%% file: paper.tex
\documentclass[a4paper,11pt]{article}
\usepackage[T1]{fontenc}
\usepackage{graphicx}
\usepackage[misc,geometry]{ifsym}
\usepackage{hyperref}
\usepackage{amsmath,amsthm}
\usepackage{amssymb}
\usepackage{ulem}
\usepackage{mathtools}
\usepackage{xcolor}
\usepackage{tikz}
\usepackage{a4wide}
\usepackage{authblk}
\usepackage{bbm}

\usetikzlibrary{arrows.meta}

\usepackage[linesnumbered,lined,boxed,commentsnumbered]{algorithm2e}

\MakeRobust{\ref}
\makeatletter
\newcommand{\labeltext}[2]{
  \@bsphack
  \csname phantomsection\endcsname
  \def\@currentlabel{#1}{\label{#2}}
  \@esphack
}
\makeatother

\newtheorem{theorem}{Theorem}[section]
\newtheorem{corollary}[theorem]{Corollary}
\newtheorem{definition}[theorem]{Definition}
\newtheorem{example}[theorem]{Example}

\newtheorem{observation}[theorem]{Observation}
\newtheorem{proposition}[theorem]{Proposition}
\newtheorem{remark}[theorem]{Remark}

\newcommand{\E}{\mathcal{E}}
\newcommand{\V}{\mathcal{V}}

\definecolor{theboss}{rgb}{1.0,0.0,0.0}

\definecolor{arifaz}{rgb}{0.7, 0.3, 0.5}
\newcommand{\af}[1]{ {#1}}
\definecolor{bluegray}{rgb}{0.4, 0.6, 0.8}
\newcommand{\jfp}[1]{ {#1}}

\def\indicator{{\mathbbm 1}}

\begin{document}

\title{Optimal transport on  gas networks  }

\author[1, \Letter]{Ariane Fazeny}
\author[1,2]{Martin Burger}
\author[3,4]{Jan-F. Pietschmann}
\affil[1]{~ {Helmholtz Imaging, Deutsches Elektronen-Synchrotron DESY, Notkestr. 85, Hamburg, 22607, Germany.} \\ \Letter ~ 
ariane.fazeny@desy.de}
\affil[2]{~ Fachbereich Mathematik, Universit\"at Hamburg, Bundesstrasse
55, Hamburg, 20146, Germany.}
\affil[3]{~ Institute of Mathematics, Universit\"{a}t Augsburg, Institut f\"ur Mathematik, Universit\"{a}tsstra\ss e 12a, 86159 Augsburg, Germany.}
\affil[4]{~ Centre for Advanced Analytics and Predictive Sciences (CAAPS), University of Augsburg,
Universit\"atsstr. 12a, 86159 Augsburg, Germany.}


\maketitle

\begin{abstract}
    \noindent Optimal transport tasks naturally arise in gas networks, which include a variety of constraints such as physical plausibility of the transport and the avoidance of extreme pressure fluctuations. To define feasible optimal transport plans, we utilize a $p$-Wasserstein metric and similar dynamic formulations minimizing the kinetic energy necessary for moving gas through the network, which we combine with suitable versions of Kirchhoff's law as coupling condition at the nodes. In contrast to existing literature, we especially focus on the non-standard case $p \neq 2$ to derive an overdamped isothermal model for gases through $p$-Wasserstein gradient flows in order to uncover and analyze underlying dynamics. We introduce different options for modeling the gas network as an oriented graph including the possibility to store gas at interior vertices and to put in or take out gas at boundary vertices. \\
    
    \noindent \textbf{Keywords} $p$-Wasserstein metric $\cdot$ Gas networks $\cdot$ PDEs on graphs $\cdot$ Isothermal Euler equations $\cdot$ Optimal control of gas $\cdot$ Optimal transport. \\

    \noindent \textbf{MSC classification} 49Q22 $\cdot$ 35R02 $\cdot$ 76N25 $\cdot$ 35Q35 $\cdot$ 60B05
\end{abstract}

\input{1_intro.tex}
\input{2_intro_gas.tex}
\input{2.1_pipe.tex}
\input{2.2_interior_vertices.tex}
\input{2.3_boundary_vertices.tex}
\input{3_OT_networks.tex}
\input{3.1_measures.tex}
\input{3.2_OT_problem.tex}
\input{3.3_feasibility.tex}
\input{3.4_solutions.tex}
\input{4_transport_cost.tex}
\input{4.1_intro_Wasserstein.tex}
\input{4.2_Wasserstein_network.tex}

\input{4.3_TODO.tex}
\input{6_gradient_flows.tex}

\input{7_numerics}

\input{8_con.tex}

\newpage

\input{appendix}

\end{document}

%% file: 1_intro.tex
\section{Introduction}

\noindent With the aim of re-purposing and extending existing natural gas infrastructure to hydrogen networks or mixed natural gas and hydrogen networks, there is a fundamental interest in understanding the dynamics of gas transport and the influence of different network topologies. \\

\noindent {In this work we model such networks as metric (or quantum) graphs where to each edge, one associates a one-dimensional interval. On these, we pose a one-dimensional partial differential equation describing the evolution of the mass density of the gas. Examples for such models are (ISO1) or (ISO3) models introduced within a whole hierarchy in \cite{domschke2021gas}. Next, we carry out a detailed modelling of possible coupling conditions at the vertices. In particular, we include the possibility of mass storage at intersections by which we extend previous approaches.
As it turns out, such models are closely related to dynamic formulations of optimal transport problems on metric graphs. This becomes evident when one only considers the conservation of mass and considers transport on networks minimizing some kinetic energy. The results of the optimal transport may give an indication of efficiency of the gas transport on the network, in particular if we extend the purely conservative framework to possible influx or outflux on some of the nodes, related to pumping in by the provider or extracting gas at the consumer. We will study the related approaches of optimal transport and also study $p-$Wasserstein metrics on the metric graphs. We mention however that those are only a metric if there is no additional in- our outflux. \\

\noindent In literature, (dynamic) optimal transport is typically restricted to bounded domains, which are a subset of $\mathbb{R}^n$ and the number of publications about  transport problems on metric graphs (or even nonconvex domains) is still quite limited and recent,  see below for a discussion. Optimal transport and Wasserstein metrics are natural tools for studying gradient flows, which we may also be interested for some of the overdamped transport models. The study of  gradient flows in Wasserstein spaces has been pioneered in \cite{erbar2021gradient} in the case $p=2$ and without storage on nodes.
Even in this case, the study of gradient flows difficult as even standard entropies are not displacement convex. 
We are able to recover (ISO3) as a gradient flow in Wasserstein spaces with $p=3$. We will also highlight some issues related to defining appropriate potentials on the metric graph, which can differ in an non-trivial way from simple definitions on single edges, since different integration constants on each edge change the interface condition. 
} \\
 

\noindent Therefore, we want to focus on how we can utilize metric graphs to model gas networks and how the $p$-Wasserstein distance allows us to tackle optimal transport problems on such graphs, where we are aiming at minimizing the necessary kinetic energy for moving gas through the network.\\



\noindent 
\noindent This paper hence analyzes how the $p$-Wasserstein distance can be defined on metric graphs encoding gas networks, in order to derive physically feasible optimal transport plans for different gas networks with varied structures, such as the possibility to store gas at interior vertices. Furthermore, we introduce dynamic formulations of the $p$-Wasserstein metric and derive $p$-Wasserstein gradient flows.

\subsection{Main contributions}
This paper generalizes the presented dynamic formulation of the $2$-Wasserstein metric in \cite{burger2023dynamic} to general $p$. Furthermore, we also present two types of coupling conditions at interior vertices (classical and generalized Kirchhoff's law) as well as time-dependent and time-independent boundary conditions. All these different types of conditions can be incorporated into the definition of our $p$-Wasserstein metric and hence solve different balanced or unbalanced optimal transport tasks on metric graphs. Furthermore, we give a detailed description of classical and weak solutions of the presented optimal transport problem as well as deriving $p$-Wasserstein gradient flows. Finally, we extend a primal-dual gradient scheme to metric graphs, with and without additional vertex dynamics.\\

\noindent The structure of the paper is as follows:
Section \ref{sec2} describes how we encode gas networks as oriented graphs, how the gas flow at individual pipes is given by the (ISO 3) model of \cite{domschke2021gas}, how we couple pipes at interior vertices and allow for storage of gas at interior vertices and how we model mass conservation at boundary vertices as well as the inflow and outflow of gas into and from the network. Section \ref{sec3} includes the introduction of masses on graphs, the formal definition of the optimal transport problem on the oriented graph together with feasibility conditions for given initial and boundary data, in particular mass conservation conditions. Section \ref{sec4} introduces the static formulation of the $p$-Wasserstein distance and different dynamic formulations depending on the coupling and boundary conditions, and discusses their basic mathematical properties. Moreover, we provide an outlook to gradient flow structures and their use for gas networks in Section \ref{sec:gf}. Finally, we discuss some numerical examples of optimal transport on networks in Section \ref{sec:num}.\\

\subsection{Related work}

\noindent The mathematical foundation of optimal transport and Wasserstein metrics, its including dynamical formulations   as well as the derivation of Wasserstein gradient flows with suitable time discretization schemes can be found in \cite{ambrosio2005gradient} and \cite{santambrogio2015optimal}. While \cite{ambrosio2005gradient} focuses on gradient flows in probability spaces, \cite{santambrogio2015optimal} gives an extensive introduction into optimal transport from the viewpoint of applied mathematics. \\

\noindent Details about the equivalence of the static and dynamic formulations of the Wasserstein distance can be found in the original publication \cite{benamou2000computational}. Dynamic formulations were extended in many directions, e.g. including non-linear mobilities \cite{Dolbeault_2008} or discrete as well as generalized graph structures \cite{MaasGradFlowEntropyFiniteMarkov2011,EPSS2021,HPS21}. Of particular interest in the context of this work is \cite{monsaingeon2021new} where two independent optimal transport problems (one in the interior of a domain and one on the boundary) are coupled so that mass exchange is possible. The idea of mass exchange between the boundary and the interior of the domain is applied to metric graphs modelling a network in \cite{burger2023dynamic}. This paper also utilizes the dynamic formulation of the $2$-Wasserstein metric on networks for optimal transport tasks, including proofs of existence of minimizers via convex duality and the derivation of $2$-Wasserstein gradient flows. The equivalence of the static and dynamic formulation of the $2$-Wasserstein metric for metric graphs is proven in \cite{erbar2021gradient}. Moreover, this paper also studies $2$-Wasserstein gradient flows of diffusion equations on metric graphs. \\

\noindent The catalogue \cite{domschke2021gas} serves as an overview of models for the transport of gas through networks and in the case of the (ISO) model hierarchy it derives from the original Euler equations different possibilities on how to model isothermal gas flow in pipes. In our paper, we will use the introduced (ISO 3) model, which based on a few assumptions on the gas flow, constitutes a simplified PDE system encoding mass and momentum conservation. Part of the model catalogue (especially the formal derivation of the (ISO 3) model) is based on the work of \cite{brouwer2011gas} and \cite{osiadacz1996different}, which use asymptotic analysis of transient gas equations to characterize different gas flow models and perform numerical simulations on examples of those models. Coupling conditions at multiple pipe connections for isothermal Euler equations are introduced in \cite{banda2005gas}. These coupling conditions resemble Kirchhoff's law but also include an equal pressure assumption at vertices of the gas network and for these coupling conditions existence of solutions at T-shaped intersections is proven. \\

\noindent In \cite{gugat2016stationary}, a method for the calculation of stationary states of gas networks is introduced. The gas flow in pipes is governed by (ISO) flow models of the previously mentioned model hierarchy and the paper specifically analyzes networks containing circles and how suitably chosen boundary conditions determine the uniqueness of stationary solutions. Furthermore, the publication \cite{egger2023stability} proves stability estimates of friction-dominated gas transport with respect to initial conditions and undertakes an asymptotic analysis of the high friction limit. In \cite{gugat2023constrained} the validity of the (ISO 2) model (a generalization of the (ISO 3) model from \cite{domschke2021gas} used in this paper) is physically valid for sufficiently low Mach numbers. The paper includes existence results of continuous solutions sufficing upper bounds for the pressure and magnitude of the Mach number of the gas flow, which are crucial parameters for physical validity of solutions.

%% file: 2_intro_gas.tex
\section{Modelling gas networks as metric graphs} \label{sec2}

Given a gas network consisting of a system of pipes, we encode it as a metric graph $\mathcal{G} = \left(\mathcal{V}, \mathcal{E}\right)$. For $m \in \mathbb{N}$, we have a set of \textbf{oriented edges}
\begin{equation*}
    \mathcal{E} = \left\{e_1, e_2, \dots e_m\right\} \subseteq \left\{e = \left\{\nu_i, \nu_j\right\} ~ \middle| ~ \nu_i, \nu_j \in \mathcal{V}\right\},
\end{equation*}
which correspond to the individual pipes with a predetermined flow direction.To each edge $e \in \mathcal{E}$, we also assign the \textbf{length}  $L_e$ and a local coordinate system $x \in (0,L_e)$, which follows the orientation.
Furthermore, we have a a set of \textbf{vertices}
\begin{equation*}
    \mathcal{V} = \left\{\nu_1, \nu_2, \dots \nu_n\right\}
\end{equation*}
for $n \in \mathbb{N}$, which encode all starting and endpoints of the individual pipes. As usual we consider the metric graph as the product space of all edges and the vertex points, where we identify the boundary points with the vertices, i.e. we take the quotient space of the closed edges subject to the identity relation of being the same vertex (cf. \cite{erbar2021gradient}).\\

\noindent Each edge $e = \left\{\nu_i, \nu_j\right\} \in \mathcal{E}$ is assigned a \textbf{start vertex} $\delta^S\left(e\right) \in \left\{\nu_i, \nu_j\right\}$ and an \textbf{end vertex} $\delta^E\left(e\right) \in \left\{\nu_i, \nu_j\right\}$, which determine the orientation of edge $e$ based on the flow direction of the corresponding pipe. Since we exclude pipes, which constitute loops from our modeling, for each edge $e = \left\{\nu_i, \nu_j\right\} \in \mathcal{E}$ it holds true that:
\begin{equation*}
    \delta^S\left(e\right) \cup \delta^E\left(e\right) = \left\{\nu_i, \nu_j\right\} \qquad \text{and} \qquad \delta^S\left(e\right) \cap \delta^E\left(e\right) = \emptyset.\\
\end{equation*}
Note that in our notation the edge coordinate system starts at $\delta^S\left(e\right)$ and ends at $\delta^E\left(e\right)$, with length $L_e$.\\

\noindent In the set of vertices, we differentiate between \textbf{boundary vertices} $\partial \mathcal{V}$, where gas enters or exits the network definitively (as supply or demand), and \textbf{interior vertices} $\mathring{\mathcal{V}}$. Hence, we first choose our boundary vertices, and then define
\begin{equation*}
    \mathring{\mathcal{V}} = \mathcal{V} \backslash \partial \mathcal{V},
\end{equation*}
which directly implies 
\begin{equation*}
    \partial \mathcal{V} \cup \mathring{\mathcal{V}} = \mathcal{V} \quad \text{and} \quad \partial \mathcal{V} \cap \mathring{\mathcal{V}} = \emptyset,
\end{equation*}
so each vertex $\nu \in \mathcal{V}$ is either a boundary vertex $\nu \in \partial \mathcal{V}$ or an interior vertex $\nu \in \mathring{\mathcal{V}}$.\\

\noindent Moreover, for the set of boundary vertices $\partial \mathcal{V}$, we differentiate between \textbf{source vertices} $\partial^+ \mathcal{V}$, which supply the network with gas, and \textbf{sink vertices} $\partial^- \mathcal{V}$, at which gas is taken out of the network to meet given demands. We will assume, that each boundary vertex $\nu \in \partial \mathcal{V}$ is either a source vertex $\nu \in \partial^+ \mathcal{V}$ or a sink vertex $\partial^- \mathcal{V}$, which implies
\begin{equation*}
    \partial^+ \mathcal{V} \cup \partial^- \mathcal{V} = \partial \mathcal{V} \quad \text{and} \quad \partial^+ \mathcal{V} \cap \partial^- \mathcal{V} = \emptyset.
\end{equation*}

\noindent In order to guarantee well-defined expressions and to simplify the notation, we assume the following properties for our underlying gas network and the resulting oriented graph.

\begin{remark}[Assumptions for the graph] ~ \\
    For our oriented graph, we will always assume that $\left\lvert\mathcal{E}\right\rvert \geq 1$ and that there are no loops $\left\{\nu, \nu\right\}$ for any $\nu \in \mathcal{V}$, which automatically implies that $\left\lvert\mathcal{V}\right\rvert \geq 2$. Moreover, we assume that our graph is connected, meaning that for every two vertices $\nu_i, \nu_j \in \mathcal{V}$ with $\nu_i \neq \nu_j$, there exists a path from $\nu_i$ to $\nu_j$ (when ignoring the orientation of edges).\\
    
    \noindent Without loss of generality, we will also assume, that the vertices are ordered in such manner, that there exist indices $o, \theta \in \left\{0, 1, \dots n\right\}$ with $o \leq \theta$ such that:
    \begin{align*}
        \partial^+ \mathcal{V} & = \left\{\nu_1, \nu_2, \dots \nu_o\right\}, \\
        \partial^- \mathcal{V} & = \left\{\nu_{o + 1}, \nu_{o + 2}, \dots \nu_{\theta}\right\}, \\
        \mathring{\mathcal{V}} & = \left\{\nu_{\theta + 1}, \nu_{\theta + 2}, \dots \nu_n\right\}.
    \end{align*}

    \noindent Furthermore, we will neglect the pipe intersection angles at interior vertices in our modelling, and only consider the length of pipes.\\
\end{remark}

\begin{example}[Gas network as oriented graph] ~ \\
    The following gas network consists of five pipes and two supply vertices and one demand vertex. Therefore, the oriented graph $\mathcal{G} = \left(\mathcal{V}, \mathcal{E}\right)$ encoding the network contains five edges
    \begin{equation*}
        \mathcal{E} = \left\{e_1, e_2, e_3, e_4, e_5\right\} = \left\{\left\{\nu_1, \nu_4\right\}, \left\{\nu_2, \nu_4\right\}, \left\{\nu_4, \nu_5\right\}, \left\{\nu_3, \nu_5\right\}, \left\{\nu_5, \nu_6\right\}\right\},
    \end{equation*}
    three boundary vertices $\partial \mathcal{V}$ and three interior vertices $\mathring{\mathcal{V}}$
    \begin{align*}
        \mathcal{V} = \partial \mathcal{V} \cup \mathring{\mathcal{V}} & = \left\{\nu_1, \nu_2, \nu_3\right\} \cup \left\{\nu_4, \nu_5, \nu_6\right\},\\
        \partial^+ \mathcal{V} & = \left\{\nu_1, \nu_2\right\},\\
        \partial^- \mathcal{V} & = \left\{\nu_3\right\}.
    \end{align*}

    \begin{minipage}{.5\textwidth}
        \begin{tikzpicture}[scale=0.6, transform shape]
		\pgfdeclarelayer{background}
		\pgfsetlayers{background,main}
		\begin{scope}[every node/.style={}]
    		\node (v1) at (0,1) {};
        	\node (v2) at (-1,-2) {};
        	\node (v5) at (2,0) {};
        	\node (v6) at (4,0) {};
                \node (v3) at (5,2) {};
        	\node (v4) at (7,-1) {};
                \node (v15) at (1.2,0.4) {};
                \node (v25) at (0.8,-0.8) {};
                \node (v56) at (3.2,0) {};
                \node (v63) at (4.65,1.3) {};
                \node (v64) at (5.8,-0.6) {};
		\end{scope}
            \begin{scope}[every node/.style={color=bluegray}]
    		\node (v1) at (0,1) {\textbf{+}};
        	\node (v2) at (-1,-2) {\textbf{+}};
                \node (v3) at (5,2) {\textbf{--}};
		\end{scope}
            \begin{scope}[>={Stealth[bluegray]},
				every edge/.style={draw=bluegray,line width=1.5pt}]
				\path [->] (v1) edge node[] {} (v15.center);
				\path [-] (1.04,0.48) edge node[] {} (v5.center);
                \path [->] (v2) edge node[] {} (v25.center);
                \path [-] (0.6,-0.93) edge node[] {} (v5.center);
                \path [->] (v5.center) edge node[] {} (v56.east);
                \path [-] (v56.west) edge node[] {} (v6.center);
                \path [->] (v6.center) edge node[] {} (v63.center);
                \path [-] (4.5,1) edge node[] {} (v3);
                \path [->] (v6.center) edge node[] {} (v64.center);
                \path [-] (5.6,-0.53) edge node[] {} (v4.center);
		\end{scope}
        \end{tikzpicture}
    \end{minipage}
    \begin{minipage}{.5\textwidth}
        \begin{tikzpicture}[scale=0.6, transform shape]
		\pgfdeclarelayer{background}
		\pgfsetlayers{background,main}
		\begin{scope}[every node/.style={draw,circle}]
                \node (v1) at (0.42,1.58) {$\nu_1$};
    	    \node (v2) at (-0.55,-2.55) {$\nu_2$};
                \node (v5) at (2,0) {$\nu_4$};
                \node (v6) at (4,0) {$\nu_5$};
                \node (v3) at (5.58,1.58) {$\nu_3$};
                \node (v4) at (6.24,-2.24) {$\nu_6$};
		\end{scope}
            \begin{scope}[>={Stealth[bluegray]},
				every edge/.style={draw=bluegray,line width=1pt}]
				\path [-] (v1) edge node[bluegray,above right] {$e_1$} (v5);
                \path [-] (v2) edge node[bluegray,above left] {$e_2$} (v5);
                \path [-] (v5) edge node[bluegray,above] {$e_3$} (v6);
                \path [-] (v6) edge node[bluegray,above left] {$e_4$} (v3);
                \path [-] (v6) edge node[bluegray,above right] {$e_5$} (v4);
		\end{scope}
        \end{tikzpicture}
    \end{minipage}
    \vspace{0cm}

    \noindent The weights of each edge $e$ are assigned according to the pipe length
    \begin{equation*}
        L_{e_1} = L_{e_4} = 45, \qquad L_{e_2} = 72, \qquad L_{e_3} = 40, \qquad L_{e_5} = 63.
    \end{equation*}

    \noindent Furthermore, the start and end vertices of each edge are given by:\newline
    \begin{minipage}{.4\textwidth}
        \begin{align*}
            \nu_1 & = \delta^S\left(e_1\right) \\
            \nu_2 & = \delta^S\left(e_2\right) \\
            \nu_3 & = \delta^E\left(e_4\right) \\
        \end{align*}
    \end{minipage}
    \begin{minipage}{.6\textwidth}
        \begin{align*}
            \nu_4 & = \delta^E\left(e_1\right) = \delta^E\left(e_2\right) = \delta^S\left(e_3\right) \\
            \nu_5 & = \delta^E\left(e_3\right) = \delta^S\left(e_4\right) = \delta^S\left(e_5\right) \\
            \nu_6 & = \delta^E\left(e_5\right) \\
        \end{align*}
    \end{minipage}
    \noindent Vertex $\nu_6$ can be considered a dead-end of our gas network, since it is only connected to one edge $e_5$, and since the vertex is neither a source nor a sink vertex.
\end{example}


%% file: 2.1_pipe.tex
\subsection{Gas flow in single pipes} \label{gasflowpipe}

\af{In a given time interval $\left[0, T\right]$ with $T > 0$, on each edge $e \in \mathcal{E}$, the physical properties of gas flow can be modelled by the \textbf{isothermal Euler equations} for compressible and inviscid fluids. In order to simplify the Euler equation system, we assume the pipe walls and the gas to have the same temperature $\mathcal{T}$, which also allows us to omit the energy conservation equation (see \cite{domschke2021gas}).} For a \textbf{mass density} $\rho_e: \left[0, L_e\right] \times \left[0, T\right] \longrightarrow \mathbb{R}_{\geq 0}$ and a \textbf{velocity} $v_e: \left[0, L_e\right] \times \left[0, T\right] \longrightarrow \mathbb{R}$, the system of equations reads as 
\begin{align*}
    \label{eqn:ISO1}
    \begin{split}
        \frac{\partial \rho_e}{\partial t} + \frac{\partial}{\partial x}\left(\rho_e v_e\right) & = 0\\\frac{\partial (\rho_e v_e)}{\partial t}
        \frac{\partial}{\partial x} \left(p_e \left(\rho_e\right)+ \rho_e v_e^2 \right) & = -\frac{\lambda_e}{2\mathcal{D}_e} \rho_e v_e \left\lvert v_e\right\rvert - g \rho_e \sin{\left(\omega_e\right)}\\
        \text{for} ~ \left(x, t\right) & \in \left(0, L_e\right) \times \left(0, T\right)
    \end{split}
    \tag{ISO1}
\end{align*}

\noindent The first equation (continuity equation) encodes conservation of mass and the second equation ensures conservation of momentum. Here, $\lambda_e \geq 0$ denotes the pipe friction coefficient, $\mathcal{D}_e > 0$ the pipe diameter, $g \approx 6.7 \cdot 10^{-11} ~ N m^2 / kg^2$ the gravitational constant, and $\omega_e \in \left[0, 2\pi\right]$ the inclination level of the pipe $e$. The \textbf{pressure} $p_e: \left[0, L_e\right] \times \left[0, T\right] \longrightarrow \mathbb{R}_{\geq 0}$ is implicitly defined by the state of real gases equation
\begin{equation}
    \label{eqn:RGE}
    p_e\left(\rho_e\right) = R \rho_e \mathcal{T} z,
    \tag{RGE}
\end{equation}  
with gas constant $R \approx 8.3 ~ J / \left(K \cdot mol\right)$, temperature $\mathcal{T} > 0$ of the gas and pipe walls, and compressibility factor $z \geq 0$, which constitutes the derivation of the pressure law for ideal gases and which we assume to be constant, but generally depends on $p_e$ and $\mathcal{T}$.\\

\noindent For each edge $e \in \mathcal{E}$, we furthermore define the  \textbf{mass flux}
\begin{equation*}
    j_e: \left[0, L_e\right] \times \left[0, T\right] \longrightarrow \mathbb{R} \quad \left(x, t\right) \mapsto j_e\left(x, t\right) := \rho_e\left(x, t\right) v_e\left(x, t\right),
\end{equation*}
and for convenience, we will use the following notation for the flux at the left (inflow) and right (outflow) boundary of each edge
\begin{equation*}
    \alpha_e\left(t\right) := j_e\left(0, t\right)   \qquad \text{and} \qquad \beta_e\left(t\right) := j_e\left(L_e, t\right), \qquad \text{for} ~ t \in \left[0, T\right].
\end{equation*}

\af{\begin{remark}[Modelling three-dimensional pipes as one-dimensional edges] ~ \\
Note that in \eqref{eqn:ISO1}, we model a three-dimensional pipe as a one-dimensional edge $e$, thus the mass density, the velocity, the pressure and the mass flux are averaged in the cross-section of the pipe to derive a one-dimensional formulation. A corresponding derivation of optimal transport from three-dimensional pipe domains to one-dimensional metric graphs constitutes an interesting question for future research.
\end{remark}}


     
   
\noindent In order to further simplify the \eqref{eqn:ISO1} model, we assume small flow rates
\begin{equation*}
    \left\lvert v_e \right\rvert \ll c_e,
\end{equation*}
which enables us to eliminate the non-linearity on the left side of the momentum equation of the  Euler equations (leading to the (ISO2) model in \cite{domschke2021gas}), as shown in \cite{osiadacz1996different}. Furthermore, we assume a friction-dominated flow, which corresponds to the friction on the pipe walls dominating the heat conduction effects, i.e. we can assume an isothermal model. If we also neglect the gravitational force in the asymptotic considerations, then we obtain, by using the scaling approaches of \cite{brouwer2011gas}, an even further simplified left side of the momentum equation and overall the model then corresponds to the  (ISO3) model in \cite{domschke2021gas}, which is given by
\begin{align*}
    \label{eqn:ISO3}
    \begin{split}
        \frac{\partial \rho_e}{\partial t} + \frac{\partial}{\partial x}\left(\rho_e v_e\right) & = 0\\
        \frac{\partial}{\partial x} \left(p_e \left(\rho_e\right)\right) & = -\frac{\lambda_e}{2\mathcal{D}_e} \rho_e v_e \left\lvert v_e\right\rvert - g \rho_e \sin{\left(\omega_e\right)}
    \end{split}
    \tag{ISO3}
\end{align*}

\noindent In this system, we can eliminate the velocity $v_e$ and write a $p-$Laplacian type equation for $\rho_e$.
Such equations can be understood as an Otto-Wasserstein gradient flow with respect to the $p$-Wasserstein metric, see for example \cite{Agueh2012_finsler}. This further motivates our study of the $p$-Wasserstein distance on metric graphs.



%% file: 2.2_interior_vertices.tex
\subsection{Gas flow at interior vertices}
For connecting the gas flow of the individual pipes to a feasible network, we need coupling conditions for the flux at vertices, which model pipe intersections. This applies to both interior vertices $\mathring{\mathcal{V}}$ and boundary vertices $\partial \mathcal{V}$, and for interior vertices $\nu \in \mathring{\mathcal{V}}$ we additionally assume that it is possible to store gas at the individual vertices.\\

\noindent In order to formulate different coupling and boundary conditions in a unified way, we introduce the \textbf{vertex excess flux}, \af{based on a generalized version of Kirchhoff's law, as}
\begin{equation*}
    \label{eqn:VF KL}
    f_{\nu}: \left[0, T\right] \longrightarrow \mathbb{R} \quad t \mapsto f_{\nu}\left(t\right) := \sum_{\substack{e \in \mathcal{E}: \\ \delta^E\left(e\right) = \nu}} \beta_e\left(t\right) - \sum_{\substack{e \in \mathcal{E}: \\ \delta^S\left(e\right) = \nu}} \alpha_e\left(t\right).\\
    \tag{VF KL}
\end{equation*}
For each vertex $\nu \in \mathcal{V}$, it catches the difference between the total amount of gas from all ingoing edges $e$ with $\delta^E\left(e\right) = \nu$ and the total amount of gas from all outgoing edges $e$ with $\delta^S\left(e\right) = \nu$.\\


\noindent To model the \af{possibility to store gas at interior vertices $\nu \in \mathring{\mathcal{V}}$}, we introduce \textbf{vertex mass densities}
\begin{equation*}
    \gamma_{\nu}: \left[0, T\right] \longrightarrow \mathbb{R}_{\geq 0} 
\end{equation*}
for all interior vertices $\nu \in \mathring{\mathcal{V}}$. At each interior vertex, we have the initial storage volume $\gamma_{\nu}\left(0\right) \geq 0$ and final storage volume $\gamma_{\nu}\left(T\right) \geq 0$, and for all $t \in \left(0, T\right)$ the vertex mass density has to satisfy
\begin{equation}\label{eq:vertex_ode}
    \frac{\partial \gamma_{\nu}}{\partial t} = f_{\nu}. \\
\end{equation}

\noindent The combination of the vertex excess flux function $f_{\nu}$ and the vertex mass density $\gamma_{\nu}$ ensures mass conservation at each interior vertex at all times, since gas can either flow through the vertex completely (indicating $f_{\nu} = 0$, as accumulated inflow equals accumulated outflow) or gas is stored at the given vertex (meaning $f_{\nu} > 0 $, as accumulated inflow is bigger than accumulated outflow) or gas is taken out of storage at the given vertex (meaning $f_{\nu} < 0 $, as accumulated inflow is smaller than accumulated outflow).\\

\noindent Not storing gas in the interior vertices can be achieved by setting $f_{\nu} \equiv 0$. Hence, we obtain the classical version of \textbf{Kirchhoff's law}
\begin{equation*}
    \label{eqn:C KL}
    0 = \sum_{\substack{e \in \mathcal{E}: \\ \delta^E\left(e\right) = \nu}} \beta_e\left(t\right) - \sum_{\substack{e \in \mathcal{E}: \\ \delta^S\left(e\right) = \nu}} \alpha_e\left(t\right).
    \tag{C KL}
\end{equation*}


\begin{observation}[No storage at interior vertices] ~ \\
    \noindent After setting $f_{\nu} \equiv 0$, the vertex mass density $\gamma_{\nu}$ can be omitted from further calculations, since $\frac{\partial \gamma_{\nu}}{\partial t} \equiv 0$ directly implies a constant amount of stored gas in each vertex $\nu \in \mathring{\mathcal{V}}$ at all times
    \begin{equation*}
        \gamma_{\nu}\left(0\right) = \gamma_{\nu}\left(T\right) \equiv \gamma_{\nu}\left(t\right), \qquad \forall t \in \left(0, T\right).\\
    \end{equation*}
\end{observation}

\noindent Note that, continuity of the velocity function $v$ or the mass density $\rho$ in the interior vertices is not covered by these coupling conditions, however it is also not expected in general.

%% file: 2.3_boundary_vertices.tex
\subsection{Gas flow at boundary vertices}

Boundary vertices $\nu \in \partial \mathcal{V}$ always either constitute source vertices or sink vertices of the network, which are responsible for the supply (inflow) and demand (outflow) of gas in the network. In contrast to interior vertices, we assume that storage of gas at boundary vertices is not possible. Therefore, we need a modified version of Kirchhoff's law, which also includes the inflow and outflow of gas to and from outside of the network, to ensure conservation of mass at boundary vertices. Hence, we introduce time-dependent and time-independent boundary conditions for boundary vertices $\nu \in \partial \mathcal{V}$.

\begin{definition}[Time-dependent boundary conditions] ~ \\
    \noindent To include time-dependent boundary conditions for the gas supply at the source vertices $\partial^+ \mathcal{V}$ as well as time-dependent gas demand conditions at the sink vertices $\partial^- \mathcal{V}$ of the gas network, we suppose the following \textbf{source vertex flux} function
    \begin{equation*}
        s_{\nu}^G: \left[0, T\right] \longrightarrow \mathbb{R}_{\leq 0} \quad t \mapsto s_{\nu}^G \left(t\right)
    \end{equation*}
    encodes the amount of gas, which has to enter the network at each source vertex $\nu \in \partial^+ \mathcal{V}$ and the following \textbf{sink vertex flux} function
    \begin{equation*}
        d_{\nu}^G: \left[0, T\right] \longrightarrow \mathbb{R}_{\geq 0} \quad t \mapsto d_{\nu}^G \left(t\right)
    \end{equation*}
    describes the given outflow of gas from the network at each sink vertex $\nu \in \partial^- \mathcal{V}$.\\
\end{definition}

\noindent Inspired by the coupling conditions at interior vertices, we can formulate these conditions as a system of equations.\\ 

\begin{observation}[Time-dependent boundary conditions] ~ \\
    \noindent With the generalized Kirchhoff's law (\ref{eqn:VF KL}), these boundary conditions can be written as
    \begin{flalign*}
        \qquad s_{\nu}^G \left(t\right) & = f_{\nu} \left(t\right) \qquad & \forall \nu \in \partial^+ \mathcal{V}, \, t \in \left[0, T\right], \qquad \\
        \qquad d_{\nu}^G \left(t\right) & = f_{\nu} \left(t\right)  \qquad & \forall \nu \in \partial^- \mathcal{V}, \, t \in \left[0, T\right], \qquad 
    \end{flalign*}
    because the source vertex flux $\left\lvert s_{\nu}^G \right\rvert$ acts as a flow into vertex $\nu \in \partial^+ \mathcal{V}$ (thus an outflow $\beta_e$ of an imaginary edge $e$) and the sink vertex flux $d_{\nu}^G$ can be seen as a flow out of vertex $\nu \in \partial^- \mathcal{V}$ (thus an inflow $\alpha_e$ of an imaginary edge $e$). With this interpretation, we obtain
    \begin{flalign*}
        \qquad 0 & = \Bigg( \left\lvert s_{\nu}^G \left(t\right) \right\rvert + \sum_{\substack{e \in \mathcal{E}: \\ \delta^E\left(e\right) = \nu}} \beta_e\left(t\right) \Bigg) - \sum_{\substack{e \in \mathcal{E}: \\ \delta^S\left(e\right) = \nu}} \alpha_e\left(t\right) \qquad & \forall \nu \in \partial^+ \mathcal{V}, \, t \in \left[0, T\right], \qquad & \\
        \qquad 0 & = \sum_{\substack{e \in \mathcal{E}: \\ \delta^E\left(e\right) = \nu}} \beta_e\left(t\right) - \Bigg( d_{\nu}^G \left(t\right) + \sum_{\substack{e \in \mathcal{E}: \\ \delta^S\left(e\right) = \nu}} \alpha_e\left(t\right) \Bigg) \qquad & \forall \nu \in \partial^- \mathcal{V}, \, t \in \left[0, T\right]. \qquad &
    \end{flalign*}
\end{observation}

\noindent Another option for conditions at the boundary vertices $\partial \mathcal{V}$ is only considering the (mandatory) total gas supply $\left\lvert S_{\nu}^G \right\rvert$ for each source vertex $\nu \in \partial^+ \mathcal{V}$ for the time period $\left[0, T\right]$, which is available starting at time $t = 0$, as well as the (mandatory) total gas demand $D_{\nu}^G$ for each sink vertex $\nu \in \partial^- \mathcal{V}$ for the time period $\left[0, T\right]$, which has to be met at the latest at time $t = T$. Since these conditions only have to be fulfilled at the end of the time period, we call them time-independent boundary conditions.\\

\begin{definition}[Time-independent boundary conditions] ~ \\
    \noindent Suppose for the time period $\left[0, T\right]$, the accumulated supplies $\left\lvert S_{\nu}^G \right\rvert$ (with $S_{\nu}^G \leq 0$) at source vertices $\nu \in \partial^+ \mathcal{V}$, and accumulated demands $D_{\nu}^G \geq 0$ at sink vertices $\nu \in \partial^- \mathcal{V}$ are given. Then our aim is to find \textbf{source vertex flux} functions
    \begin{equation*}
        s_{\nu}: \left[0, T\right] \longrightarrow \mathbb{R}_{\leq 0} \quad t \mapsto s_{\nu}\left(t\right) := f_{\nu} \left(t\right)
    \end{equation*}
    for each source vertex $\nu \in \partial^+ \mathcal{V}$, such that
    \begin{equation}
        \label{totS}
        S_{\nu}^G = \int_0^T s_{\nu}\left(t\right) \, \mathrm d t,
    \end{equation}
    and to find \textbf{sink vertex flux} functions
    \begin{equation*}
        d_{\nu}: \left[0, T\right] \longrightarrow \mathbb{R}_{\geq 0} \quad t \mapsto d_{\nu}\left(t\right) := f_{\nu} \left(t\right)
    \end{equation*}
    for each sink vertex $\nu \in \partial^- \mathcal{V}$, such that
    \begin{equation}
        \label{totD}
        D_{\nu}^G = \int_0^T d_{\nu}\left(t\right) \, \mathrm d t.\\
    \end{equation}
\end{definition}

\noindent Inspired by the coupling conditions at interior vertices, we can once again reformulate these conditions.\\ 

\begin{observation}[Time-independent boundary conditions] ~ \\
    \noindent If we define a \textbf{source vertex mass density} $S_{\nu}: \left[0, T\right] \longrightarrow \mathbb{R}_{\geq 0}$
    for all source vertices $\nu \in \partial^+ \mathcal{V}$, then constraint (\ref{totS}), can be rewritten as
    \def\filler{s_{\nu} \left(t\right)}
    \begin{flalign*}
        \qquad \frac{\partial S_{\nu} \left(t\right)}{\partial t} & = s_{\nu} \left(t\right) \qquad & \forall \nu \in \partial^+ \mathcal{V}, \, t \in \left(0, T\right), \qquad & \\
    \end{flalign*}
    with initial and final condition $S_{\nu}\left(0\right) = \left\lvert S_{\nu}^G \right\rvert$ and $S_{\nu}\left(T\right) = 0$. \\
    
    \noindent In a similar manner, by introducing a \textbf{sink vertex mass density} $D_{\nu}: \left[0, T\right] \longrightarrow \mathbb{R}_{\geq 0}$
    for each sink vertex $\nu \in \partial^- \mathcal{V}$, constraint (\ref{totD}) can be rewritten as
    \begin{flalign*}
        \qquad \frac{\partial D_{\nu} \left(t\right)}{\partial t} & = d_{\nu} \left(t\right) \qquad & \forall \nu \in \partial^- \mathcal{V}, \, t \in \left(0, T\right), \qquad &
    \end{flalign*}
    with initial and final condition $D_{\nu}\left(0\right) = 0 $ and $D_{\nu}\left(T\right) = D_{\nu}^G$. \\

\end{observation}

\noindent In the case of time-independent boundary conditions the non-positivity property of the source vertex flux and the non-negativity of the sink vertex flux are not automatically fulfilled. \\

\begin{remark}[Non-positivity of $s_{\nu}$ and non-negativity of $d_{\nu}$] ~ \\
    For time-independent boundary conditions, it is necessary to constrain $s_{\nu}\left(t\right)$ to $\mathbb{R}_{\leq 0}$ and $d_{\nu}\left(t\right)$ to $\mathbb{R}_{\geq 0}$ for all times $t \in \left[0, T\right]$. Otherwise, $s_{\nu}\left(t\right) > 0$ would indicate giving back gas supply and $d_{\nu}\left(t\right) < 0$ would correspond to taking back gas demand, which are both unrealistic behaviors for real-world applications. \\
\end{remark}

\noindent The source vertex mass and sink vertex mass densities can also be defined in the case of time-dependent boundary conditions as
\begin{equation*}
    \def\fillerb{S_{\nu}^G}
    \def\fillera{s_{\nu}^G \left(\tilde{t}\right)}
    S_{\nu} \left(t\right) := \left\lvert 
    \smash{\underbrace\fillerb_{\leq \, 0}} \vphantom\fillerb \right\rvert - \left\lvert \int_0^t \smash{\underbrace\fillera_{\leq \, 0}}\vphantom\fillera \, \mathrm d \tilde{t} \right\rvert \quad \text{with} \quad S_{\nu}^G := \int_0^T s_{\nu}^G \left(\tilde{t}\right) \, \mathrm d \tilde{t} \vspace{0.1cm}\\
\end{equation*}
for source vertices $\nu \in \partial^+ \mathcal{V}$, and as
\begin{equation*}
    D_{\nu} \left(t\right) := \int_0^t d_{\nu}^G \left(\tilde{t}\right) \, \mathrm d \tilde{t} \quad \text{with} \quad D_{\nu}^G := \int_0^T d_{\nu}^G \left(\tilde{t}\right) \, \mathrm d \tilde{t}
\end{equation*}
for sink vertices $\nu \in \partial^- \mathcal{V}$. \\

\noindent For any time $t \in \left[0, T\right]$, the source vertex mass density $S_{\nu} \left(t\right)$ gives the remaining supply at time $t$, which has not entered the network yet, and the sink vertex mass density $D_{\nu} \left(t\right)$ encodes the amount of gas which has exited the network until time $t$. \\

\begin{proposition}[Source vertex mass density] ~ \\
    The source vertex mass density $S_{\nu} \left(t\right)$ of a source vertex $\nu \in \partial^+ \mathcal{V}$ at any time $t \in \left[0, T\right]$ can be calculated by
    \begin{equation*}
        S_{\nu}\left(t\right) = \left\lvert \int_t^T s_{\nu}^G \left(\tilde{t}\right) \, \mathrm d \tilde{t} \right\rvert
    \end{equation*}
    for time-dependent boundary conditions, or with $s_{\nu}$ instead of $s_{\nu}^G$ for time-independent boundary conditions.\\
\end{proposition}

\begin{proof}
    For the source vertex mass density $S_{\nu} \left(t\right)$ of source vertex $\nu \in \partial^+ \mathcal{V}$ at time $t \in \left[0, T\right]$, the following reformulations hold true (taking into account that $s_{\nu}^G\left(t\right) \leq 0$ for all source vertices $\nu \in \partial^+ \mathcal{V}$ and all time points $t \in \left[0, T\right]$, and $S_{\nu}^G \leq 0$ for all time points $t \in \left[0, T\right]$):
    \begin{flalign*}
        S_{\nu} \left(t\right) & = \left\lvert S_{\nu}^G \right\rvert - \left\lvert \int_0^t s_{\nu}^G \left(\tilde{t}\right) \, \mathrm d \tilde{t} \right\rvert = \int_0^t s_{\nu}^G \left(\tilde{t}\right) \, \mathrm d \tilde{t} - S_{\nu}^G = \int_0^t s_{\nu}^G \left(\tilde{t}\right) \, \mathrm d \tilde{t} - \int_0^T s_{\nu}^G \left(\tilde{t}\right) \, \mathrm d \tilde{t} \\
        & = - \int_t^T s_{\nu}^G \left(\tilde{t}\right) \, \mathrm d \tilde{t} = \left\lvert \int_t^T s_{\nu}^G \left(\tilde{t}\right) \, \mathrm d \tilde{t} \right\rvert
    \end{flalign*}

    \noindent Exactly the same proof can be carried out with $s_{\nu}$ instead of $s_{\nu}^G$ in the case of time-independent boundary conditions.\\
\end{proof}

\noindent With the vertex flux functions and vertex mass densities at boundary vertices, we can also calculate the total gas supply $\left\lvert S \right\rvert$ (with $S \leq 0$) entering the network and the total gas demand $D$ exiting the network in the time period $\left[0, T\right]$. \\

\begin{observation}[Total supply and total demand] ~ \\
    In the case of time-dependent boundary conditions, the total supply and total demand of the network for the time period $\left[0, T\right]$ are given by
    \begin{flalign*}
        S & = \sum_{\nu \in \partial^+ \mathcal{V}} - S_{\nu}\left(0\right) = \sum_{\nu \in \partial^+ \mathcal{V}} S_{\nu}^G = \sum_{\nu \in \partial^+ \mathcal{V}} \int_0^T s_{\nu}^G \left(t\right) \, \mathrm d t, \\
        D & = \sum_{\nu \in \partial^- \mathcal{V}} D_{\nu} \left(T\right) = \sum_{\nu \in \partial^- \mathcal{V}} D_{\nu}^G = \sum_{\nu \in \partial^- \mathcal{V}} \int_0^T d_{\nu}^G \left(t\right) \, \mathrm d t.
    \end{flalign*}
    Alternatively, exactly the same equations hold true with $s_{\nu}$ instead of $s_{\nu}^G$ and $d_{\nu}$ instead of $d_{\nu}^G$ for time-independent boundary conditions.
\end{observation}

%% file: 3_OT_networks.tex
\section{Transport type problems on metric graphs} \label{sec3}

While the previous section was dedicated to the formulation of gas transport models on metric graphs, especially coupling conditions, we now undertake the next step towards our goal to understand such models as gradient flows. For this purpose, we will introduce optimal transport type problems on metric graphs in this section. We obtain different variants, depending on the type of coupling and boundary conditions, all of which can be understood as extensions of the dynamic formulation of the classical $p$-Wasserstein distance. However, only in the case of no boundary vertices, do we obtain a metric or distance.\\

\noindent Ultimately, in section \ref{sec:gf}, we will show that we can recover the \eqref{eqn:ISO3} model from a minimizing movements (or JKO) scheme with respect to one of these metrics. However, we start by introducing some further notation.

%% file: 3.1_measures.tex
\subsection{Measure spaces on graphs}

\noindent \af{By $\mathcal{M}_+\left(\Omega\right)$ we denote the set of Borel measures on a metric space $\left(\Omega, d\right)$ with $\Omega \subseteq \mathbb{R}^k$ and $k \in \mathbb{N}$ (and $\Sigma$ being the corresponding $\sigma$-algebra of Borel sets of $\Omega$), which are non-negative and bounded, so}
\begin{equation*}
    \forall X \subseteq \Omega: 0 \leq \mu\left(X\right) \qquad \text{and} \qquad \mu\left(\Omega\right) < \infty.
\end{equation*}
Based on this notation, we choose for each edge $e \in \mathcal{E}$ the mass density $\rho_e$ from the set of non-negative bounded measures on edge $e$ for the time period $\left[0, T\right]$ with
\begin{equation*}
    \mathcal{M}_+\left(e\right) := \mathcal{M}_+\left(\left[0, L_e\right] \times \left[0, T\right]\right).
\end{equation*}
Furthermore, we choose the vertex mass density $\gamma_{\nu}$ for each interior vertex $\nu \in \mathring{\mathcal{V}}$, the source vertex mass $S_{\nu}$ for each source vertex $\nu \in \partial^+ \mathcal{V}$, and the sink vertex mass $D_{\nu}$ for each sink vertex $\nu \in \partial^- \mathcal{V}$, from the set of non-negative bounded measures on vertex $\nu$ for the time period $\left[0, T\right]$, so we define
\begin{equation*}
    \mathcal{M}_+\left(\nu\right) := \mathcal{M}_+\left(\left\{\nu\right\} \times \left[0, T\right]\right). \\
\end{equation*}

\noindent For a fixed time $t \in \left[0, T\right]$, we define in a similar manner 
\begin{equation*}
    \mathcal{M}_+^t\left(e\right) := \mathcal{M}_+\left(\left[0, L_e\right]\right) \qquad \text{and} \qquad \mathcal{M}_+^t\left(\nu\right) := \mathcal{M}_+\left(\left\{\nu\right\}\right) = a \delta_{\nu},
\end{equation*}
for an edge $e \in \mathcal{E}$ and a vertex $\nu \in \mathcal{V}$, where $\delta_{\nu}$ denotes the corresponding Dirac measure with $a \in \mathbb{R}_{\geq 0}$. \\

\noindent For the velocity functions $v_e$ and the mass flux functions $j_e$ of edges $e \in \mathcal{E}$, and for the vertex excess flux function $f_{\nu}$ of interior vertices $\nu \in \mathring{\mathcal{V}}$, we do not need the non-negativity of the defined measure spaces. Hence, by $\mathcal{M}\left(\Omega\right)$ we denote the set of bounded Borel measures $\mu$ on a metric space $\left(\Omega, d\right)$, and thus we define
\begin{equation*}
    \mathcal{M}\left(e\right) := \mathcal{M}\left(\left[0, L_e\right] \times \left[0, T\right]\right) \quad \text{and} \quad \mathcal{M}\left(\nu\right) := \mathcal{M}\left(\left\{\nu\right\} \times \left[0, T\right]\right)
\end{equation*}
for any edge $e \in \mathcal{E}$ and any vertex $\nu \in \mathcal{V}$. \\

\noindent With the set of measures being defined on a single edge or a single vertex of the graph, we can now also introduce coupled measures on a subset of edges or a subset of vertices. These coupled measure are defined on the domain of a cartesian product of $\mathcal{M}_+\left(e\right)$, $\mathcal{M}_+^t\left(e\right)$, $\mathcal{M} \left(e\right)$, $\mathcal{M}_+\left(\nu\right)$, $\mathcal{M}_+^t\left(\nu\right)$ or $\mathcal{M} \left(\nu\right)$ for edges $e \in \mathcal{E}$ and vertices $\nu \in \mathcal{V}$. A formal definition can be found in the \ref{appendix}. \\


%% file: 3.2_OT_problem.tex
\subsection{Formulation of the optimal transport problems}

Depending on the the existence of boundary vertices and the type of boundary conditions, as well as the coupling conditions at interior vertices, we obtain different formulations of the optimal transport problem on the graph $\mathcal{G} = \left(\mathcal{V}, \mathcal{E}\right)$.\\

\noindent Note that, in this subsection we are only defining a general optimal transport task on the topology of a graph, and hence do not consider the (\ref{eqn:ISO3}) model for gas flow in the pipes yet. Instead we only consider the general continuity equation together with a version of Kirchhoff's law as a coupling condition at the vertices, which ensures mass conservation of our transport. They will serve as side constraints to an optimization problem with a (for now) general cost functional. As we will see in section \ref{sec:gf}, building a minimizing movement scheme using these optimal transport problems as distance, we will indeed recover the \eqref{eqn:ISO3} model. \\

\noindent Assuming the case with no boundary vertices and classical Kirchhoff's law as coupling condition (as in \cite{erbar2021gradient}), the corresponding optimal transport problem can formally be written as:
\begin{flalign*}
    \inf_{\rho \in \mathcal{M}_+ \left(\mathcal{E}\right), \, v \in \mathcal{M}\left(\mathcal{E}\right)} \quad c\left(\rho, v\right) \qquad \text{subject to} &&
\end{flalign*}
\vspace{-0.5cm}
\begin{alignat*}{3}
    && 0 & = \frac{\partial \rho_e}{\partial t} + \frac{\partial}{\partial x}\left(\rho_e v_e\right) & \forall e \in \mathcal{E}, \, x \in \left[0, L_e\right], \, t \in \left[0, T\right] \\ 
    && 0 & = \sum_{\substack{e \in \mathcal{E}: \\ \delta^E\left(e\right) = \nu}} \left. \left(\rho_e v_e\right) \right|_{x = L_e} - \sum_{\substack{e \in \mathcal{E}: \\ \delta^S\left(e\right) = \nu}} \left. \left(\rho_e v_e\right) \right|_{x = 0} \qquad \qquad \quad & \forall \nu \in \mathring{\mathcal{V}}, \, t \in \left[0, T\right] \\
    && \left.\rho_e\right|_{t = 0} & = \left(\rho_0\right)_e, \, \left.\rho_e\right|_{t = T} = \left(\rho_T\right)_e & \forall e \in \mathcal{E}
\end{alignat*}

\noindent If we rewrite the presented optimal transport problem with the mass flux $j_e$ instead of the product $\rho_e v_e$ for each edge $e \in \mathcal{E}$,  we obtain a convex problem with linear side constraints. Based on this form, we obtain the following extensions for different combinations of coupling- and boundary conditions. \\

\begin{definition}[Optimal transport problem]\label{def:OT_problems} ~ \\
    If there are no boundary vertices ($o = \theta = 0$, $\left\lvert \partial \mathcal{V}\right\rvert = 0$ so $\mathcal{V} = \mathring{\mathcal{V}}$) and thus no boundary conditions, we can formulate the optimal transport problem with coupling conditions based on the generalized Kirchhoff's law (\ref{eqn:VF KL}) as:
    \begin{flalign} \label{eqn:opt_ina} \tag{OT VF KL}
        \inf_{\substack{\rho \in \mathcal{M}_+ \left(\mathcal{E}\right), \, j \in \mathcal{M}\left(\mathcal{E}\right), \\ \gamma \in \mathcal{M}_+ \left(\mathring{\mathcal{V}}\right), \, f \in \mathcal{M}\left(\mathring{\mathcal{V}}\right)}} \quad c\left(\rho, j, \gamma, f\right) \qquad \text{subject to} &&
    \end{flalign}
    \vspace{-0.5cm}
    \begin{alignat*}{3}
        && 0 & = \frac{\partial \rho_e}{\partial t} + \frac{\partial j_e}{\partial x} \qquad \qquad \qquad \qquad \qquad \qquad & \forall e \in \mathcal{E}, \, x \in \left[0, L_e\right], \, t \in \left[0, T\right] \label{CE 1} \tag{CE 1} \\ 
        && f_{\nu}\left(t\right) & = \sum_{\substack{e \in \mathcal{E}: \\ \delta^E\left(e\right) = \nu}} \left. j_e \right|_{x = L_e} - \sum_{\substack{e \in \mathcal{E}: \\ \delta^S\left(e\right) = \nu}} \left. j_e \right|_{x = 0} & \forall \nu \in \mathring{\mathcal{V}}, \, t \in \left[0, T\right] \label{VF KL 1} \tag{VF KL 1} \\
        && \frac{\partial \gamma_{\nu}}{\partial t} & = f_{\nu} & \forall \nu \in \mathring{\mathcal{V}}, \, t \in \left[0, T\right] \label{VF KL 2} \tag{VF KL 2} \\
        && \left.\rho_e\right|_{t = 0} & = \left(\rho_0\right)_e, \, \left.\rho_e\right|_{t = T} = \left(\rho_T\right)_e & \forall e \in \mathcal{E} \label{CE 2} \tag{CE 2} \\
        && \left.\gamma_{\nu}\right|_{t = 0} & = \left(\gamma_0\right)_{\nu}, \, \left.\gamma_{\nu}\right|_{t = T} = \left(\gamma_T\right)_{\nu} & \forall \nu \in \mathring{\mathcal{V}} \label{VF KL 3} \tag{VF KL 3} 
    \end{alignat*}
    
    \noindent When using the classical version of Kirchhoff's law (\ref{eqn:C KL}) at the interior vertices (no storage of gas at interior vertices), then this simplifies to:
    \begin{flalign} \label{eqn:opt_erbar} \tag{OT C KL}
        \inf_{\rho \in \mathcal{M}_+ \left(\mathcal{E}\right), \, j \in \mathcal{M}\left(\mathcal{E}\right)} \quad c\left(\rho, j\right)\text{ subject to } &&
    \end{flalign}
    \vspace{-0.5cm}
    \begin{alignat}{3}
        && 0 & = \frac{\partial \rho_e}{\partial t} + \frac{\partial j_e}{\partial x} \qquad \qquad \qquad \qquad \qquad \qquad \quad & \forall e \in \mathcal{E}, \, x \in \left[0, L_e\right], \, t \in \left[0, T\right] \tag{CE 1}\label{eqn:CE1}\\ 
        && 0 & = \sum_{\substack{e \in \mathcal{E}: \\ \delta^E\left(e\right) = \nu}} \left. j_e \right|_{x = L_e} - \sum_{\substack{e \in \mathcal{E}: \\ \delta^S\left(e\right) = \nu}} \left. j_e \right|_{x = 0} & \forall \nu \in \mathring{\mathcal{V}}, \, t \in \left[0, T\right] \label{C KL 1} \tag{C KL 1} \\
        && \left.\rho_e\right|_{t = 0} & = \left(\rho_0\right)_e, \, \left.\rho_e\right|_{t = T} = \left(\rho_T\right)_e & \forall e \in \mathcal{E} \tag{CE 2}
    \end{alignat}
\end{definition}

\noindent The optimal transport problem \eqref{eqn:opt_ina} is exactly the same as the model studied in \cite{burger2023dynamic}, while \eqref{eqn:opt_erbar} was analyzed in \cite{erbar2021gradient}. In both cases, no boundary vertices are present and thus, mass transport into or out of the network is not possible. In view of our application, we hence extend these problem by the following boundary conditions.\\

\begin{remark}[Extension with boundary conditions] ~ \\
    \noindent If boundary vertices are present ($0 < \theta$, $\left\lvert \partial \mathcal{V}\right\rvert > 0$ so $\mathcal{V} = \partial \mathcal{V} \cup \mathring{\mathcal{V}} \supset \mathring{\mathcal{V}}$), then depending on the type of boundary conditions, different constraints need to be added to the optimal transport problem. Assuming time-dependent boundary conditions, then the following constraints need to be included:
    \begin{alignat}{3}
        && 0 & = \Bigg( \left\lvert s_{\nu}^G \right\rvert + \sum_{\substack{e \in \mathcal{E}: \\ \delta^E\left(e\right) = \nu}} \left. j_e \right|_{x = L_e} \Bigg) - \sum_{\substack{e \in \mathcal{E}: \\ \delta^S\left(e\right) = \nu}} \left. j_e \right|_{x = 0} \qquad \qquad & \forall \nu \in \partial^+ \mathcal{V}, \, t \in \left[0, T\right] \label{TD BC 1} \tag{TD BC 1} \\ 
        && 0 & = \sum_{\substack{e \in \mathcal{E}: \\ \delta^E\left(e\right) = \nu}} \left. j_e \right|_{x = L_e} - \Bigg( d_{\nu}^G + \sum_{\substack{e \in \mathcal{E}: \\ \delta^S\left(e\right) = \nu}} \left. j_e \right|_{x = 0} \Bigg) & \forall \nu \in \partial^- \mathcal{V}, \, t \in \left[0, T\right] \label{TD BC 2} \tag{TD BC 2}
    \end{alignat}

    \noindent Similarly, for time-independent boundary conditions, the following constraints, as well as initial and final conditions need to be included in the optimal transport problem instead:
    \begin{alignat}{3}
        && s_{\nu} & = \sum_{\substack{e \in \mathcal{E}: \\ \delta^S\left(e\right) = \nu}} \left. j_e \right|_{x = 0} - \sum_{\substack{e \in \mathcal{E}: \\ \delta^E\left(e\right) = \nu}} \left. j_e \right|_{x = L_e} & \qquad \qquad \qquad \quad \forall \nu \in \partial^+ \mathcal{V}, \, t \in \left[0, T\right] \label{TI BC 1} \tag{TI BC 1} \\
        && d_{\nu} & = \sum_{\substack{e \in \mathcal{E}: \\ \delta^E\left(e\right) = \nu}} \left. j_e \right|_{x = L_e} - \sum_{\substack{e \in \mathcal{E}: \\ \delta^S\left(e\right) = \nu}} \left. j_e \right|_{x = 0} & \forall \nu \in \partial^- \mathcal{V}, \, t \in \left[0, T\right] \label{TI BC 2} \tag{TI BC 2} \\
        && \frac{\partial S_{\nu}}{\partial t} & = s_{\nu} & \forall \nu \in \partial^+ \mathcal{V}, \, t \in \left[0, T\right] \label{TI BC 3} \tag{TI BC 3} \\
        && \frac{\partial D_{\nu}}{\partial t} & = d_{\nu} & \forall \nu \in \partial^- \mathcal{V}, \, t \in \left[0, T\right] \label{TI BC 4} \tag{TI BC 4} \\
        && \left.S_{\nu}\right|_{t = 0} & = \left\lvert S_{\nu}^G \right\rvert, \, \left.S_{\nu}\right|_{t = T} = 0 & \forall \nu \in \partial^+ \mathcal{V} \label{TI BC 5} \tag{TI BC 5}\\
        && \left.D_{\nu}\right|_{t = 0} & = 0, \, \left.D_{\nu}\right|_{t = T} = D_{\nu}^G & \forall \nu \in \partial^- \mathcal{V} \label{TI BC 6} \tag{TI BC 6}
    \end{alignat}
    
    \noindent In order to avoid $s_{\nu}\left(t\right) = \frac{\partial S_{\nu} \left(t\right)}{\partial t} > 0$ or $d_{\nu}\left(t\right) = \frac{\partial D_{\nu} \left(t\right)}{\partial t} < 0$ at all times $t \in \left[0, T\right]$, we require $-s \in \mathcal{M}_+ \left(\partial^+ \mathcal{V}\right)$ and $d \in \mathcal{M}_+ \left(\partial^- \mathcal{V}\right)$. Furthermore, for time-independent boundary conditions, the cost functional changes to
    \begin{equation*}
        \inf_{\substack{\rho \in \mathcal{M}_+ \left(\mathcal{E}\right), \, j \in \mathcal{M}\left(\mathcal{E}\right), \\ \gamma \in \mathcal{M}_+ \left(\mathring{\mathcal{V}}\right), \, f \in \mathcal{M}\left(\mathring{\mathcal{V}}\right), \\ S, \, -s \in \mathcal{M}_+ \left(\partial^+ \mathcal{V}\right), \\ D, \, d \in \mathcal{M}_+ \left(\partial^- \mathcal{V}\right)}} c\left(\rho, j, \gamma, f, S, s, D, d\right), \qquad \text{or} \qquad \inf_{\substack{\rho \in \mathcal{M}_+ \left(\mathcal{E}\right), \, j \in \mathcal{M}\left(\mathcal{E}\right), \\ S, \, -s \in \mathcal{M}_+ \left(\partial^+ \mathcal{V}\right), \\ D, \, d \in \mathcal{M}_+ \left(\partial^- \mathcal{V}\right)}} c\left(\rho, j, S, s, D, d\right),
    \end{equation*}
    depending on the type of coupling conditions at the interior vertices.\\
\end{remark}

\noindent \af{Examples for possible (action or) \textbf{cost functionals} will be given in  section \ref{sec:p_Wasserstein_networks}, aiming at calculating the necessary kinetic energy for the transport, which is sensible for gas networks, but is also the type of action functional typically used for dynamic optimal transport.}

%% file: 3.3_feasibility.tex
\subsection{Feasibility of the optimal transport problem}

Depending on whether we allow storage of gas at interior vertices $\mathring{\mathcal{V}}$ or not, and depending on the existence of boundary vertices and the type of boundary conditions, necessary conditions for the feasibility of the optimal transport problem differ slightly. \\

\noindent If we write the space of the graph $\mathcal{G} = \left(\mathcal{V}, \mathcal{E}\right)$ as
\begin{equation*}
    \Omega_{\mathcal{G}} := \left(\bigtimes_{\nu \in \mathcal{V}} \left\{\nu\right\}\right) \times \left(\bigtimes_{e \in \mathcal{E}} \left[0, L_e\right]\right),
\end{equation*}
then we can introduce a \textbf{total density} $\sigma$ for the entire graph.\\

\begin{definition}[Total density] ~ \\
    In the case of generalized Kirchhoff's law as coupling condition at interior vertices, and time-dependent or time-independent boundary conditions, we define the total density on graph $\mathcal{G} = \left(\mathcal{V}, \mathcal{E}\right)$ as
    \begin{equation*}
        \sigma: \Omega_{\mathcal{G}} \times \left[0, T\right] \longrightarrow \mathbb{R}_{\geq 0} 
    \end{equation*}
    \begin{equation*}
        \left(\nu_1, \nu_2, \dots \nu_n, x_{e_1}, x_{e_2}, \dots x_{e_m}, t\right) \mapsto \sigma \left(t\right) := \sum_{\nu \in \mathring{\mathcal{V}}} \gamma_{\nu} \left(t\right) + \sum_{\nu \in \partial^+ \mathcal{V}} S_{\nu} \left(t\right) + \sum_{\nu \in \partial^- \mathcal{V}} D_{\nu} \left(t\right) + \sum_{e \in \mathcal{E}} \rho_e \left(x_e, t\right),
    \end{equation*}
    \af{where $x_{e_i} \in \left[0, L_{e_i}\right]$ for all edges $e_i \in \mathcal{E}$.}\\
    
    \noindent Depending on the existence of boundary vertices and the coupling conditions at interior vertices, this definition simplifies to:
    \begin{itemize}
        \item Classical Kirchhoff's law and no boundary vertices: \\ $\sigma \left(t\right) = \sum_{e \in \mathcal{E}} \rho_e \left(x_e, t\right)$ 
        \item Classical Kirchhoff's law and time-dependent or time-independent boundary conditions: \\ $\sigma \left(t\right) = \sum_{\nu \in \partial^+ \mathcal{V}} S_{\nu} \left(t\right) + \sum_{\nu \in \partial^- \mathcal{V}} D_{\nu} \left(t\right) + \sum_{e \in \mathcal{E}} \rho_e \left(x_e, t\right)$ 
        \item Generalized Kirchhoff's law and no boundary vertices: \\
        $\sigma \left(t\right) = \sum_{\nu \in \mathring{\mathcal{V}}} \gamma_{\nu} \left(t\right) + \sum_{e \in \mathcal{E}} \rho_e \left(x_e, t\right)$
    \end{itemize}
\end{definition}

\noindent Note that for time-dependent boundary conditions, the source vertex mass density $S_{\nu}$ and sink vertex mass density $D_{\nu}$ at time $t \in \left[0, T\right]$ can easily be calculated by utilizing the given source vertex flux $s_{\nu}^G$ and sink vertex flux $d_{\nu}^G$ functions with
\begin{equation*}
    S_{\nu} \left(t\right) = \left\lvert \int_t^T s_{\nu}^G (\tilde{t}\,) \, \mathrm d \tilde{t} \right\rvert, \qquad \text{and} \qquad D_{\nu} (t) = \int_0^t d_{\nu}^G \left(\tilde{t\,}\right) \, \mathrm d \tilde{t}.
\end{equation*}

\begin{remark}[Total density a probability measure] ~ \\
    For the measures $\rho_e$, $\gamma_{\nu}$, $S_{\nu}$, and $D_{\nu}$ we will always assume that they are chosen in such a way, that for all times $t \in \left[0, T\right]$ the total density $\sigma \left(t\right)$ is a probability measure in $\mathbb{P}\left(\Omega_{\mathcal{G}}\right)$. This means that the total mass of gas in the network is constant for all $t \in \left[0, T\right]$, so
    \begin{equation*}
        \Sigma \left(t\right) := \sum_{\nu \in \mathring{\mathcal{V}}} \gamma_{\nu} \left(t\right) + \sum_{\nu \in \partial^+ \mathcal{V}} S_{\nu} \left(t\right) + \sum_{\nu \in \partial^- \mathcal{V}} D_{\nu} \left(t\right) + \sum_{e \in \mathcal{E}} \int_0^{L_e} \rho_e \left(x, t\right) \, \mathrm d x \equiv C \in \mathbb{R}_{> 0},
    \end{equation*}
    in the case of generalized Kirchhoff's law as coupling condition and the existence of boundary vertices. For the initial and final conditions of the density functions $\gamma_{\nu}$, $S_{\nu}$, $D_{\nu}$, and $\rho_e$ we will specifically demand that they are given in such way that $\sigma \left(0\right), \sigma \left(T\right) \in \mathbb{P}\left(\Omega_{\mathcal{G}}\right)$.\\
\end{remark}

\noindent The condition of global mass conservation of the gas can be expressed utilizing the total density.\\

\begin{definition}[Global mass conservation] ~ \\
    The total mass in a graph $\mathcal{G} = \left(\mathcal{V}, \mathcal{E}\right)$ includes the gas stored at interior vertices, the supply, which has not entered the network yet, the demand, which has already exited the network, and the gas, which is currently in the pipes of the network. Therefore, the global mass conservation condition is given by
    \begin{equation*}
        \frac{\mathrm d \Sigma \left(t\right)}{\mathrm d t} = 0 \qquad \forall t \in \left(0, T\right).
    \end{equation*}

    \noindent Depending on the type of coupling conditions and whether there are any boundary vertices, the total mass is given by:
    \begin{itemize}
        \item Classical Kirchhoff's law and no boundary vertices: \\ $\Sigma \left(t\right) = \sum_{e \in \mathcal{E}} \int_0^{L_e} \rho_e \left(x, t\right) \, \mathrm d x$ 
        \item Classical Kirchhoff's law and time-dependent or time-independent boundary conditions: \\ $\Sigma \left(t\right) = \sum_{\nu \in \partial^+ \mathcal{V}} S_{\nu} \left(t\right) + \sum_{\nu \in \partial^- \mathcal{V}} D_{\nu} \left(t\right) + \sum_{e \in \mathcal{E}} \int_0^{L_e} \rho_e \left(x, t\right) \, \mathrm d x$ 
        \item Generalized Kirchhoff's law and no boundary vertices: \\
        $\Sigma \left(t\right) = \sum_{\nu \in \mathring{\mathcal{V}}} \gamma_{\nu} \left(t\right) + \sum_{e \in \mathcal{E}} \int_0^{L_e} \rho_e \left(x, t\right) \, \mathrm d x$ 
        \item Generalized Kirchhoff's law and time-dependent or time-independent boundary conditions: \\
        $\Sigma \left(t\right) = \sum_{\nu \in \mathring{\mathcal{V}}} \gamma_{\nu} \left(t\right) + \sum_{\nu \in \partial^+ \mathcal{V}} S_{\nu} \left(t\right) + \sum_{\nu \in \partial^- \mathcal{V}} D_{\nu} \left(t\right) + \sum_{e \in \mathcal{E}} \int_0^{L_e} \rho_e \left(x, t\right) \, \mathrm d x$
    \end{itemize}

    \noindent As before, in the case of time-dependent boundary conditions, the source vertex mass density $S_{\nu}$ and sink vertex mass density $D_{\nu}$ can easily be calculated.\\
\end{definition}

\noindent The global mass conservation conditions can be used to define feasible upper bounds for the total demand of all sink vertices for the time period $\left[0, T\right]$. \\

\begin{proposition}[Upper bound for the accumulated demand] ~ \\
    In the case of generalized Kirchhoff's law at interior vertices and the existence of boundary vertices, the accumulated demand $D = \sum_{\nu \in \partial^- \mathcal{V}} D_{\nu} \left(T\right)$ of the time period $\left[0, T\right]$ is bounded from above through the inequality
    \begin{equation*}
        \sum_{\nu \in \mathring{\mathcal{V}}} \gamma_{\nu} \left(0\right) + \left\lvert S \right\rvert + \sum_{e \in \mathcal{E}} \int_0^{L_e} \rho_e \left(x, 0\right) \, \mathrm d x \geq D,
    \end{equation*}
    where $\left\lvert S \right\rvert = \sum_{\nu \in \partial^- \mathcal{V}} S_{\nu} \left(0\right)$ corresponds to the accumulated supply for the time period $\left[0, T\right]$. \\
\end{proposition}

\begin{proof}
    \noindent If we compare the total mass $\Sigma \left(t\right)$ for any $t \in \left[0, T\right]$ to the initial total mass $\Sigma \left(0\right)$, then we obtain the following equation
    \begin{align*}
        & \sum_{\nu \in \mathring{\mathcal{V}}} \gamma_{\nu} \left(0\right) + \sum_{\nu \in \partial^+ \mathcal{V}} \underbrace{S_{\nu} \left(0\right)}_{= \, \left\lvert S_{\nu}^G \right\rvert} + \sum_{\nu \in \partial^- \mathcal{V}} \underbrace{D_{\nu} \left(0\right)}_{= \, 0} + \sum_{e \in \mathcal{E}} \int_0^{L_e} \rho_e \left(x, 0\right) \, \mathrm d x = \\
        & \sum_{\nu \in \mathring{\mathcal{V}}} \gamma_{\nu} \left(t\right) + \sum_{\nu \in \partial^+ \mathcal{V}} S_{\nu} \left(t\right) + \sum_{\nu \in \partial^- \mathcal{V}} D_{\nu} \left(t\right) + \sum_{e \in \mathcal{E}} \int_0^{L_e} \rho_e \left(x, t\right) \, \mathrm d x,
    \end{align*}
    which can be rewritten to
    \begin{align*}
        \label{eqn:GMC}
        \begin{split}
            & \sum_{\nu \in \mathring{\mathcal{V}}} \gamma_{\nu} \left(0\right) + \sum_{\nu \in \partial^+ \mathcal{V}} \underbrace{\left(\left\lvert S_{\nu}^G \right\rvert - S_{\nu} \left(t\right)\right)}_{= \, \left\lvert \int_0^t s_{\nu} \left(\tilde{t}\right) \, \mathrm d \tilde{t} \right\rvert} + \sum_{e \in \mathcal{E}} \int_0^{L_e} \rho_e \left(x, 0\right) \, \mathrm d x = \\
            & \sum_{\nu \in \mathring{\mathcal{V}}} \gamma_{\nu} \left(t\right) + \sum_{\nu \in \partial^- \mathcal{V}} \underbrace{D_{\nu} \left(t\right)}_{= \, \int_0^t d_{\nu} \left(\tilde{t}\right) \, \mathrm d \tilde{t}} + \sum_{e \in \mathcal{E}} \int_0^{L_e} \rho_e \left(x, t\right) \, \mathrm d x.
        \end{split}
        \tag{GMC}
    \end{align*}

    \noindent The left side of the global mass conservation equation calculates the amount of gas stored in the network at time $t = 0$ (in interior vertices or in the pipes) plus the gas, which has entered the network through source vertices in the time period $\left[0, t\right]$. The right side encodes the amount of gas in the interior vertices or the pipes at time $t$ together with the gas, which has exited the network through sink vertices in the time period $\left[0, t\right]$.\\
    
    \noindent For $t = T$ we obtain from the global mass conservation equation (\ref{eqn:GMC})
    \begin{flalign*}
        & \sum_{\nu \in \mathring{\mathcal{V}}} \gamma_{\nu} \left(0\right) + \sum_{\nu \in \partial^+ \mathcal{V}} \left\lvert \int_0^T s_{\nu} \left(t\right) \, \mathrm d t \right\rvert + \sum_{e \in \mathcal{E}} \int_0^{L_e} \rho_e \left(x, 0\right) \, \mathrm d x = \\
        & \sum_{\nu \in \mathring{\mathcal{V}}} \gamma_{\nu} \left(T\right) + \sum_{\nu \in \partial^- \mathcal{V}} \int_0^T d_{\nu} \left(t\right) \, \mathrm d t + \sum_{e \in \mathcal{E}} \int_0^{L_e} \rho_e \left(x, T\right) \, \mathrm d x,
    \end{flalign*}
    which is equivalent to the following equation
    \begin{equation*}
        \sum_{\nu \in \mathring{\mathcal{V}}} \gamma_{\nu} \left(0\right) + \left\lvert S \right\rvert + \sum_{e \in \mathcal{E}} \int_0^{L_e} \rho_e \left(x, 0\right) \, \mathrm d x = \sum_{\nu \in \mathring{\mathcal{V}}} \gamma_{\nu} \left(T\right) + D + \sum_{e \in \mathcal{E}} \int_0^{L_e} \rho_e \left(x, T\right) \, \mathrm d x.
    \end{equation*}
    By utilizing the non-negativity of the vertex mass density $\gamma$ and of the mass density $\rho$
    \begin{equation*}
        \gamma_{\nu} \in \mathcal{M}_+ \left(\mathring{\mathcal{V}}\right), \, \rho_e \in \mathcal{M}_+ \left(\mathcal{E}\right) \quad \Longrightarrow \quad \sum_{\nu \in \mathring{\mathcal{V}}} \gamma_{\nu} \left(T\right), \, \sum_{e \in \mathcal{E}} \int_0^{L_e} \rho_e \left(x, T\right) \, \mathrm d x \geq 0,
    \end{equation*}
    we obtain the upper bound for the accumulated demand of all sink vertices $\nu \in \partial^- \mathcal{V}$
    \begin{equation*}
        \sum_{\nu \in \mathring{\mathcal{V}}} \gamma_{\nu} \left(0\right) + \left\lvert S \right\rvert + \sum_{e \in \mathcal{E}} \int_0^{L_e} \rho_e \left(x, 0\right) \, \mathrm d x \geq D.
    \end{equation*}
\end{proof}

{\corollary[Upper bound for the demand with time-dependent boundary conditions] ~ \\
    If the graph has boundary vertices with time-dependent boundary conditions and uses generalized Kirchhoff's law as coupling condition at interior vertices, the accumulated demand at all times $t \in \left[0, T\right]$ is bounded above by
    \begin{equation*}
        \sum_{\nu \in \mathring{\mathcal{V}}} \gamma_{\nu} \left(0\right) + \sum_{\nu \in \partial^+ \mathcal{V}} \left\lvert \int_0^t s_{\nu}^G \left(\tilde{t}\right) \, \mathrm d \tilde{t}\right\rvert + \sum_{e \in \mathcal{E}} \int_0^{L_e} \rho_e \left(x, 0\right) \, \mathrm d x \geq \sum_{\nu \partial^- \mathcal{V}} \int_0^t d_{\nu}^G \left(\tilde{t}\right) \, \mathrm d \tilde{t},
    \end{equation*}
    which can be used to test the given source vertex flux $s_{\nu}^G$ and sink vertex flux $d_{\nu}^G$ for plausibility. \\
}

\begin{proof}
    \noindent \noindent As the vertex mass density $\gamma$ and the mass density $\rho$ on the left side of the global mass conservation equality (\ref{eqn:GMC}) do not depend on $t \in \left[0, T\right]$, and the non-negativity
    \begin{equation*}
        \gamma_{\nu} \in \mathcal{M}_+ \left(\mathring{\mathcal{V}}\right), \, \rho_e \in \mathcal{M}_+ \left(\mathcal{E}\right) \quad \Longrightarrow \quad \sum_{\nu \in \mathring{\mathcal{V}}} \gamma_{\nu} \left(t\right), \, \sum_{e \in \mathcal{E}} \int_0^{L_e} \rho_e \left(x, t\right) \, \mathrm d x \geq 0,
    \end{equation*}
    holds true for all $t \in \left[0, T\right]$, we obtain the inequality
    \begin{equation*}
        \sum_{\nu \in \mathring{\mathcal{V}}} \gamma_{\nu} \left(0\right) + \sum_{\nu \in \partial^+ \mathcal{V}} \left(\left\lvert S_{\nu}^G \right\rvert - S_{\nu} \left(t\right)\right) + \sum_{e \in \mathcal{E}} \int_0^{L_e} \rho_e \left(x, 0\right) \, \mathrm d x \geq \sum_{\nu \partial^- \mathcal{V}} D_{\nu} \left(t\right).
    \end{equation*}
    Plugging in the definitions of the source vertex mass densities $S_{\nu}^G$ and $S_{\nu}$, as well as the sink vertex mass densities $D_{\nu}$, yields the inequality
    \begin{equation*}
        \sum_{\nu \in \mathring{\mathcal{V}}} \gamma_{\nu} \left(0\right) + \sum_{\nu \in \partial^+ \mathcal{V}} \left\lvert \int_0^t s_{\nu}^G \left(\tilde{t}\right) \, \mathrm d \tilde{t}\right\rvert + \sum_{e \in \mathcal{E}} \int_0^{L_e} \rho_e \left(x, 0\right) \, \mathrm d x \geq \sum_{\nu \partial^- \mathcal{V}} \int_0^t d_{\nu}^G \left(\tilde{t}\right) \, \mathrm d \tilde{t}.
    \end{equation*}
\end{proof}

\noindent The global mass conservation equation (\ref{eqn:GMC}) together with the upper bound for the accumulated demand can be used to test given initial, final and boundary conditions for possibly leading to an infeasible optimal transport problem. However, note that these tests do not detect instances which for example require unfeasibly fast velocities $v$. \\


%% file: 3.4_solutions.tex
\subsection{Strong and weak solutions of the continuity equation constraints}

\noindent Assuming feasibility of the introduced optimal transport problem, we now want to define strong and weak solutions of the \af{continuity equation constraints (\ref{CE 1}), (\ref{CE 2}) combined with the constraints from the coupling conditions and possibly the boundary conditions.} 

\begin{definition}[Strong solutions] ~ \\
    A strong solution of (\ref{CE 1}), (\ref{CE 2}) together with (\ref{C KL 1}) or (\ref{VF KL 1}) - (\ref{VF KL 3}), and possibly also (\ref{TD BC 1}), (\ref{TD BC 2}) or (\ref{TI BC 1}) - (\ref{TI BC 6}), consists of a tuple of functions $\left(\rho, j\right)$, $\left(\rho, j, S, s, D, d\right)$, $\left(\rho, j, \gamma, f\right)$ or $\left(\rho, j, \gamma, f, S, s, D, d\right)$ fulfilling a subset of the following conditions:
    \begin{itemize}
        \item[(S1)] $\forall e \in \mathcal{E}: \left(x, t\right) \mapsto \rho_e \left(x, t\right)$ continuous on $\left[0, L_e\right] \times \left[0, T\right]$\\
        $\forall e \in \mathcal{E}, \, \forall x \in \left[0, L_e\right]: t \mapsto \rho_e \left(x, t\right)$ continuously differentiable on $\left(0, T\right)$
        \item[(S2)] $\forall e \in \mathcal{E}: \left(x, t\right) \mapsto j_e \left(x, t\right)$ continuous on $\left[0, L_e\right] \times \left[0, T\right]$\\
        $\forall e \in \mathcal{E}, \, \forall t \in \left[0, T\right]: x \mapsto j_e \left(x, t\right)$ continuously differentiable on $\left(0, L_e\right)$
        \item[(S3)]
        \begin{enumerate}
           \item[i)]$\rho \in \mathcal{M}_+ \left(\mathcal{E}\right)$, $j \in \mathcal{M} \left(\mathcal{E}\right)$ fulfill the side constraints (\ref{CE 1}), (\ref{CE 2}) and (\ref{C KL 1}) 
            \item[ii)] $\rho \in \mathcal{M}_+ \left(\mathcal{E}\right)$, $j \in \mathcal{M} \left(\mathcal{E}\right)$, $\gamma \in \mathcal{M}_+ \left(\mathring{\mathcal{V}}\right)$ and $f \in \mathcal{M} \left(\mathring{\mathcal{V}}\right)$ fulfill the side constraints (\ref{CE 1}), (\ref{CE 2}) and (\ref{VF KL 1}) - (\ref{VF KL 3}) 
        \end{enumerate}
        \item[(S4)] $\forall \nu \in \mathring{\mathcal{V}}: t \mapsto \gamma_{\nu} \left(t\right)$ continuous on $\left[0, T\right]$, and continuously differentiable on $\left(0, T\right)$\\
        (which directly implies $\forall \nu \in \mathring{\mathcal{V}}: t \mapsto f_{\nu} \left(t\right)$ continuous on $\left(0, T\right)$)
        \item[(S5)] $\forall \nu \in \partial^+ \mathcal{V}: t \mapsto S_{\nu} \left(t\right)$ continuous on $\left[0, T\right]$, and continuously differentiable on $\left(0, T\right)$\\
        (which directly implies $\forall \nu \in \partial^+ \mathcal{V}: t \mapsto s_{\nu} \left(t\right)$ continuous on $\left(0, T\right)$)
        \item[(S6)] $\forall \nu \in \partial^- \mathcal{V}: t \mapsto D_{\nu} \left(t\right)$ continuous on $\left[0, T\right]$, and continuously differentiable on $\left(0, T\right)$\\
        (which directly implies $\forall \nu \in \partial^- \mathcal{V}: t \mapsto d_{\nu} \left(t\right)$ continuous on $\left(0, T\right)$)
        \item[(S7)] $S, -s \in \mathcal{M}_+\left(\partial^+ \mathcal{V}\right)$ and $D, d \in \mathcal{M}_+\left(\partial^- \mathcal{V}\right)$ fulfill the constraints (\ref{TD BC 1}), (\ref{TD BC 2}) or (\ref{TI BC 1}) - (\ref{TI BC 6})  
    \end{itemize}

    \noindent Strong solutions are hence given by:
    \begin{itemize}
        \item Classical Kirchhoff's law and no boundary vertices: \\ 
        $\left(\rho, j\right)$ fulfilling (S1), (S2) and (S3) i)
        \item Classical Kirchhoff's law and time-dependent or time-independent boundary conditions: \\ 
        $\left(\rho, j, S, s, D, d\right)$ fulfilling (S1), (S2), (S3) i), (S5), (S6) and (S7)
        \item Generalized Kirchhoff's law and no boundary vertices: \\
        $\left(\rho, j, \gamma, f\right)$ fulfilling (S1), (S2), (S3) ii) and (S4)
        \item Generalized Kirchhoff's law and time-dependent or time-independent boundary conditions: \\
        $\left(\rho, j, \gamma, f, S, s, D, d\right)$ fulfilling (S1), (S2), (S3) ii), (S4), (S5), (S6) and (S7)
    \end{itemize}

    \noindent Note that, the vertex excess flux $f$ can be calculated directly from $j$ based on the generalized Kirchhoff's law at the interior vertices and since the vertex flux is not restricted to a non-negative sign, it can also be omitted from the definition of strong solutions. As the source vertex flux $s$ and the sink vertex flux $d$ need to fulfill non-positivity and non-negativity conditions respectively, they can only be omitted from the definition of strong solutions, if corresponding constraints for the sign of the time derivatives of $S$ and $D$ are included. 
\end{definition}
\noindent For the definition of weak solutions, we need suitable test functions defined on the domain of the graph.\begin{definition}[Test function]   We call a continuous function $\varphi: \Omega_{\mathcal{G}} \times \left[0, T\right] \longrightarrow \mathbb{R}$ a test function, if for all edges $e \in \mathcal{E}$ its restriction $\varphi_e :=\left. \varphi \right|_{e}$ is continuously differentiable on $\left(0, L_e\right) \times \left(0, T\right)$.
    In addition, for any vertex $\nu \in \mathcal{V}$ at any time $t \in \left[0, T\right]$ we use the notation  $\varphi_{\nu} := \left.\varphi\right|_{v}$. 
\end{definition}

%
\noindent The following definitions of weak solutions can be derived from the presented strong solutions. 
\begin{definition}[Weak solutions] ~ \\
    The tuples $\left(\rho, j\right)$, $\left(\rho, j, S, s, D, d\right)$, $\left(\rho, j, \gamma, f\right)$ or $\left(\rho, j, \gamma, f, S, s, D, d\right)$ constitute weak solutions of (\ref{CE 1}), (\ref{CE 2}) combined with (\ref{C KL 1}) or (\ref{VF KL 1}) - (\ref{VF KL 3}), and possibly also (\ref{TD BC 1}), (\ref{TD BC 2}) or (\ref{TI BC 1}) - (\ref{TI BC 6}), if they fulfill a particular subset of the following conditions:
    \begin{itemize}
        \item[(S1)] $\rho \in \mathcal{M}_+ \left(\mathcal{E}\right)$ and $j \in \mathcal{M} \left(\mathcal{E}\right)$
        \item[(S2)] $\forall e \in \mathcal{E}: \left.\rho_e\right|_{t = 0} = \left(\rho_0\right)_e, \, \left.\rho_e\right|_{t = T} = \left(\rho_T\right)_e$
        \item[(S3)]
        \begin{enumerate}
            \item[i)]
        \begin{equation*}
            \int_0^T \sum_{e \in \mathcal{E}} \int_0^{L_e} \left\lvert j_e \left(x, t\right) \right\rvert \, \mathrm d x \, \mathrm d t < \infty
        \end{equation*}
        \item[ii)] \begin{equation*}
            \int_0^T \Bigg(\sum_{e \in \mathcal{E}} \int_0^{L_e} \left\lvert j_e \left(x, t\right) \right\rvert \, \mathrm d x - \sum_{\nu \in \partial^+ \mathcal{V}} s_{\nu} \left(t\right) + \sum_{\nu \in \partial^- \mathcal{V}} d_{\nu} \left(t\right)\Bigg) \, \mathrm d t < \infty
        \end{equation*}
        \item[iii)] \begin{equation*}
            \int_0^T \Bigg( \sum_{e \in \mathcal{E}} \int_0^{L_e} \left\lvert j_e \left(x, t\right) \right\rvert \, \mathrm d x + \sum_{\nu \in \mathring{\mathcal{V}}} \left\lvert f_{\nu} \left(t\right) \right\rvert \Bigg)\, \mathrm d t < \infty
        \end{equation*}
        \item[iv)] \begin{equation*}
            \int_0^T \Bigg(\sum_{e \in \mathcal{E}} \int_0^{L_e} \left\lvert j_e \left(x, t\right) \right\rvert \, \mathrm d x + \sum_{\nu \in \mathring{\mathcal{V}}} \left\lvert f_{\nu} \left(t\right) \right\rvert - \sum_{\nu \in \partial^+ \mathcal{V}} s_{\nu} \left(t\right) + \sum_{\nu \in \partial^- \mathcal{V}} d_{\nu} \left(t\right)\Bigg) \, \mathrm d t < \infty
        \end{equation*}
        \end{enumerate}
        \item[(S4)] $\gamma \in \mathcal{M}_+ \left(\mathring{\mathcal{V}}\right)$ and $f \in \mathcal{M} \left(\mathring{\mathcal{V}}\right)$
        \item[(S5)] $\forall \nu \in \mathring{\mathcal{V}}: \left.\gamma_{\nu}\right|_{t = 0} = \left(\gamma_0\right)_{\nu}, \, \left.\gamma_{\nu}\right|_{t = T} = \left(\gamma_T\right)_{\nu}$
        \item[(S6)] $S, -s \in \mathcal{M}_+ \left(\partial^+ \mathcal{V}\right)$
        \item[(S7)] $D, d \in \mathcal{M}_+ \left(\partial^- \mathcal{V}\right)$
        \item[(S8)] $\forall \nu \in \partial^+ \mathcal{V}: \, \left.S_{\nu}\right|_{t = 0} = \left\lvert S_{\nu}^G \right\rvert, \, \left.S_{\nu}\right|_{t = T} = 0$ \vspace{0.2cm} \\
        $\forall \nu \in \partial^+ \mathcal{V}$ and for $\mathcal{L}$-a.e. $t \in \left(0, T\right): \, \frac{\partial S_{\nu} \left(t\right)}{\partial t} = s_{\nu} \left(t\right)$ \vspace{0.4cm} \\
        $\forall \nu \in \partial^- \mathcal{V}: \, \left.D_{\nu}\right|_{t = 0} = 0, \, \left.D_{\nu}\right|_{t = T} = D_{\nu}^G$ \vspace{0.2cm} \\
        $\forall \nu \in \partial^- \mathcal{V}$ and for $\mathcal{L}$-a.e. $t \in \left(0, T\right): \, \frac{\partial D_{\nu} \left(t\right)}{\partial t} = d_{\nu} \left(t\right)$
        \item[(S9)] 
        \begin{enumerate}
            \item[i)]
        $\forall$ test functions $\varphi$ and $\mathcal{L}$-a.e. $t \in \left(0, T\right):$
        \begin{flalign*}
            & \sum_{e \in \mathcal{E}} \int_0^{L_e} \frac{\partial \rho_e \left(x, t\right)}{\partial t} \varphi_e \left(x, t\right) \, \mathrm d x = \sum_{e \in \mathcal{E}} \int_0^{Le} j_e \left(x, t\right) \frac{\partial \varphi_e \left(x, t\right)}{\partial x} \, \mathrm d x
        \end{flalign*}
        \item[ii)] $\forall$ test functions $\varphi$ and $\mathcal{L}$-a.e. $t \in \left(0, T\right):$
        \begin{flalign*}
            & \sum_{e \in \mathcal{E}} \int_0^{L_e} \frac{\partial \rho_e \left(x, t\right)}{\partial t} \varphi_e \left(x, t\right) \, \mathrm d x + \sum_{\nu \in \partial^+ \mathcal{V}} s_{\nu} \left(t\right) \varphi_{\nu} \left(t\right) + \sum_{\nu \in \partial^- \mathcal{V}} d_{\nu} \left(t\right) \varphi_{\nu} \left(t\right) = \\
            & \sum_{e \in \mathcal{E}} \int_0^{Le} j_e \left(x, t\right) \frac{\partial \varphi_e \left(x, t\right)}{\partial x} \, \mathrm d x
        \end{flalign*}
        \item[iii)] $\forall$ test functions $\varphi$ and $\mathcal{L}$-a.e. $t \in \left(0, T\right):$
        \begin{equation*}
            \sum_{e \in \mathcal{E}} \int_0^{L_e} \frac{\partial \rho_e \left(x, t\right)}{\partial t} \varphi_e \left(x, t\right) \, \mathrm d x + \sum_{\nu \in \mathring{\mathcal{V}}} f_{\nu} \left(t\right) \varphi_{\nu} \left(t\right) = \sum_{e \in \mathcal{E}} \int_0^{Le} j_e \left(x, t\right) \frac{\partial \varphi_e \left(x, t\right)}{\partial x} \, \mathrm d x
        \end{equation*}
        \item[iv)] $\forall$ test functions $\varphi$ and $\mathcal{L}$-a.e. $t \in \left(0, T\right):$
        \begin{flalign*}
            & \sum_{e \in \mathcal{E}} \int_0^{L_e} \frac{\partial \rho_e \left(x, t\right)}{\partial t} \varphi_e \left(x, t\right) \, \mathrm d x + \sum_{\nu \in \mathring{\mathcal{V}}} f_{\nu} \left(t\right) \varphi_{\nu} \left(t\right) + \sum_{\nu \in \partial^+ \mathcal{V}} s_{\nu} \left(t\right) \varphi_{\nu} \left(t\right) \\
            &+ \sum_{\nu \in \partial^- \mathcal{V}} d_{\nu} \left(t\right) \varphi_{\nu} \left(t\right) =  \sum_{e \in \mathcal{E}} \int_0^{Le} j_e \left(x, t\right) \frac{\partial \varphi_e \left(x, t\right)}{\partial x} \, \mathrm d x
        \end{flalign*}
        \end{enumerate}
    \end{itemize}

    \noindent Here, $\mathcal{L}$ denotes the standard $1$-dimensional Lebesgue measure on $\left[0, T\right]$. \\

    \noindent Weak solutions are thus given by:
    \begin{itemize}
        \item Classical Kirchhoff's law and no boundary vertices: \\ 
        $\left(\rho, j\right)$ fulfilling (S1), (S2), (S3) i) and (S9) i)
        \item Classical Kirchhoff's law and time-dependent or time-independent boundary conditions: \\ 
        $\left(\rho, j, S, s, D, d\right)$ fulfilling (S1), (S2), (S3) ii), (S6), (S7), (S8) and (S9) ii)
        \item Generalized Kirchhoff's law and no boundary vertices: \\
        $\left(\rho, j, \gamma, f\right)$ fulfilling (S1), (S2), (S3) iii), (S4), (S5) and (S9) iii)
        \item Generalized Kirchhoff's law and time-dependent or time-independent boundary conditions: \\
        $\left(\rho, j, \gamma, f, S, s, D, d\right)$ fulfilling (S1), (S2), (S3) iv), (S4), (S5), (S6), (S7), (S8) and (S9) iv)
    \end{itemize}
\end{definition}


%% file: 4_transport_cost.tex
\section{$p$-Wasserstein distance as cost of transport} \label{sec4}

In order to find solutions of the optimal transport task on a given gas network, we first need to assign suitable costs to each feasible transport plan. Thus, we will first review different formulations of the classical $p$-Wasserstein distances, and secondly extend those to our case.

%% file: 4.1_intro_Wasserstein.tex
\subsection{General $p$-Wasserstein distance}

\af{Given a metric space $\left(\Omega, d\right)$, where we assume the domain $\Omega \subseteq \mathbb{R}^k$ with $k \in \mathbb{N}$ to be bounded, the $p$-Wasserstein distance can be defined in the form of a static formulation based on the distance function $d$, or it can also be given in a dynamic formulation based on absolutely continuous curves
\begin{equation*}
    \mu: \left(0, T\right) \longrightarrow \mathcal{M}_+ \left(\Omega_{\mathcal{G}}\right) \quad t \mapsto \mu\left(t\right).
\end{equation*}}

\noindent \af{In literature, the study of the $p$-Wasserstein distance mostly utilizes the distance function
\begin{equation*}
    d \left(x, y\right) := \left\lvert x - y\right\rvert,
\end{equation*}
where $\left\lvert \, \cdot \, \right\rvert$ denotes the vector norm in $\mathbb{R}^k$. The following static formulation of the $p$-Wasserstein distance could analogously be formulated for a general distance function $d$, however for the dynamic formulation of the $p$-Wasserstein distance, $d \left(x, y\right) = \left\lvert x - y\right\rvert$ is a necessary requirement.} \\


\begin{definition}[$p$-Wasserstein distance \cite{ambrosio2005gradient}] ~ \\
    For two probability measures $\mu_0, \mu_T \in \mathbb{P}\left(\Omega\right)$, the \textbf{static formulation} of the $p$-Wasserstein distance for $p \in \left[1, \infty \right)$ is defined as
    \begin{equation*}
        W_p \left(\mu_0, \mu_T\right) = \min_{\pi \in \Pi\left(\mu_0, \mu_T\right)} \left\{\int_{\Omega \times \Omega} \left\lvert x - y\right\rvert^p \, \mathrm d \pi\left(x, y\right)\right\}^{\frac{1}{p}},
    \end{equation*}
    with $\Pi\left(\mu_0, \mu_T\right)$ being the set of all joint probability distributions $\pi$ on $\Omega \times \Omega$, with the respective marginals $\mu_0$ and $\mu_T$, i.e.
    \begin{equation*}
        \mu_0\left(x\right) = \int_{\Omega} \pi\left(x, y\right) \, \mathrm d y \qquad \text{and} \qquad \mu_T\left(y\right) = \int_{\Omega} \pi\left(x, y\right) \, \mathrm d x.
    \end{equation*}
    
    \noindent \af{Note that, as we assume $\Omega$ to be bounded, the finiteness of the $p$-th moments of $\mu_0$ and $\mu_T$ are not required. However, for non-negative bounded measures $\mu_0, \mu_t \in \mathcal{M}_+\left(\Omega\right)$ instead of probability measures, the $p$-Wasserstein distance is only able to assign costs to balanced optimal transport, thus $\mu_0$ and $\mu_T$ need to have the same total mass
    \begin{equation*}
        \int_{\Omega} \mu_0\left(z\right) \, \mathrm d z = \int_{\Omega} \mu_T\left(z\right) \, \mathrm d z.\\
    \end{equation*}}
    
    
\end{definition}


\noindent In their seminal work \cite{benamou2000computational}, Benamou and Brenier introduced a dynamical formulation of the $2$-Wasserstein distance that was also extended to the case $p\neq 2$. In this formulation the task is to minimize an action functional, which corresponds to the kinetic energy of curves connecting the initial measure $\mu_0$ to the final measure $\mu_T$. This minimization is subject to the initial and final conditions, $\left.\mu\right|_{t = 0} = \mu_0$ and $\left.\mu\right|_{t = T} = \mu_T$, as well as the continuity equation $\frac{\partial \mu}{\partial t} + \frac{\partial j}{\partial x} = 0$ and therefore closely resembles the models introduced in section \ref{sec2}. \\


\noindent \jfp{Note that, to increase readability, we will always denote the densities of measures with respect to a given reference measure such as the Lebesgue measure by the same symbol as the original measures themselves.}  

\noindent \af{We directly present the dynamic formulation of $p$-Wasserstein distance in terms of mass density and mass flux (instead of mass density and velocity), and for general time intervals $\left[0,T\right]$ (instead of $\left[0, 1\right]$). Note that, we write $\mu$ instead of $\rho$, since later on, the measure $\mu$ will not only consist of the mass density $\rho$, but also the vertex mass $\gamma$, and optionally the source vertex mass $S$ and the sink vertex mass $D$.} \\

\noindent Furthermore, we will restrict the set of feasible mass densities to $\mu \in \mathcal{M}_+ \left(\Omega \times \left[0, T\right]\right)$ and in order to avoid a positivity constraint for $\mu$, when reformulating $j := v \mu$ to $v = \frac{j}{\mu}$, we define the helper function
\begin{equation*}
    \label{eqn:helper}
    h: \mathbb{R} \times \mathbb{R} \longrightarrow \mathbb{R}_{\geq 0} \quad \left(a, b\right) \mapsto h\left(a, b\right) = \left\{ \begin{array}{lr}
        \frac{\left\lvert a \right\rvert^p}{b^{p - 1}} \quad & b > 0 \vspace{0.1cm}\\
        0 \quad & b = a = 0\\
        \infty \quad & b < 0 ~ \vee ~ b = 0, \, a \neq 0
    \end{array} \right. .
\end{equation*}


\noindent \af{This allows us to write the dynamic formulation of the $p$-Wasserstein distance in terms of mass densities and mass fluxes, which results in the minimization problem becoming convex with linear side constraints.} \\

\begin{observation}[Dynamic formulation of the $p$-Wasserstein distance] ~ \\
    On a convex and compact domain $\Omega$, the $p$-Wasserstein distance for two probability measures $\mu_0, \mu_1 \in \mathbb{P} \left(\Omega\right)$ can be calculated as
    \begin{equation*}
        \label{eqn:DYN-T} \tag{DYN-T}
        W_p^p \left(\mu_0, \mu_T\right) =
    \end{equation*}
    \begin{equation*}
        \min_{\substack{\mu \in \mathcal{M}_+ \left(\Omega \times \left[0, T\right]\right), \\ j \in \mathcal{M} \left(\Omega \times \left[0, T\right]\right)}} \left\{T^{p - 1} \int_0^T \int_{\Omega} h\left(j\left(x, t\right), \mu\left(x, t\right)\right) \, \mathrm d\eta  \, \mathrm d t ~ \middle| ~ \begin{array}{l}
            \frac{\partial \mu}{\partial t} + \frac{\partial j}{\partial x} = 0, \vspace{0.1cm} \\
            \left. \mu \right|_{t = 0} = \mu_0, \, \left. \mu \right|_{t = T} = \mu_T
        \end{array} \right\},
    \end{equation*}
    where $\eta$ is some reference measure.
%
    Note that, this formulation of the $p$-Wasserstein distance is restricted to physically feasible solutions as $W_p^p \left(\mu_0, \mu_T\right) < \infty$ is only possible, if for $\mathcal{L}$-a.e. $x \in \Omega$ and $t \in \left[0, T\right]$ either
    \begin{equation*}
        \mu \left(x, t\right) > 0 \quad \text{or} \quad \mu \left(x, t\right) = j \left(x, t\right) = 0,
    \end{equation*}
    so a non-zero flux is only possible if the mass density is positive, and a negative mass density always constitutes an infeasible transport. \\


    \noindent Here, $\mathcal{L}$ denotes the $k + 1$-dimensional Lebesgue measure on $\Omega \times \left[0, T\right] \subseteq \mathbb{R}^{k + 1}$.
\end{observation}

%% file: 4.2_Wasserstein_network.tex
\subsection{$p$-Wasserstein distance on networks} \label{sec:p_Wasserstein_networks}

\noindent In order to utilize the $p$-Wasserstein distance as the cost functional for the transport along an edge $e \in \mathcal{E}$, we can set $\Omega := \left[0, L_e\right]$ as well as $k := 1$, and with the notation of the previously introduced optimal transport problems we can thus write for a given initial mass distribution of gas $\left(\rho_0\right)_e \in \mathbb{P} \left(\left[0, L_e\right]\right)$, a given final distribution $\left(\rho_T\right)_e \in \mathbb{P} \left(\left[0, L_e\right]\right)$ and the mass flux $j_e$ along the edge $e$ 
\begin{equation*}
    W_p^p \left(\left(\rho_0\right)_e, \left(\rho_T\right)_e\right) =
\end{equation*}
\begin{equation*}
    \min_{\substack{\rho_e \in \mathcal{M}_+ \left(e\right), \\ j_e \in \mathcal{M} \left(e\right)}} \left\{T^{p - 1} \int_0^T \int_0^{L_e} h\left(j_e\left(x, t\right), \rho_e\left(x, t\right)\right) \, \mathrm d \eta \, \mathrm d t ~ \middle| ~ \begin{array}{l}
        \frac{\partial \rho_e}{\partial t} + \frac{\partial j_e}{\partial x} = 0, \vspace{0.1cm} \\
        \left. \rho_e \right|_{t = 0} = \left(\rho_0\right)_e, \, \left. \rho_e \right|_{t = T} = \left(\rho_T\right)_e
    \end{array} \right\} , \\
\end{equation*}
for a reference measure $\eta$.
\noindent For a single interior vertex $\nu \in \mathring{\mathcal{V}}$, in the case of generalized Kirchhoff's law as coupling condition, we could also write the cost of transport from an initial vertex mass density $\left(\gamma_0\right)_{\nu} \in \mathbb{P} \left(\left\{\nu\right\}\right)$ to a final vertex mass density $\left(\gamma_T\right)_{\nu} \in \mathbb{P} \left(\left\{\nu\right\}\right)$ with vertex excess flux $f_{\nu}$ as
\begin{equation*}
    W_p^p \left(\left(\gamma_0\right)_{\nu}, \left(\gamma_T\right)_{\nu}\right) =
\end{equation*}
\begin{equation*}
    \min_{\substack{\gamma_{\nu} \in \mathcal{M}_+ \left(\nu\right), \\ f_{\nu} \in \mathcal{M} \left(\nu\right)}} \left\{T^{p - 1} \int_0^T h\left(f_{\nu}\left(t\right), \gamma_{\nu}\left(t\right)\right) \, \mathrm d t ~ \middle| ~ \begin{array}{l}
        \frac{\partial \gamma_{\nu}}{\partial t} = f_{\nu}, \vspace{0.1cm} \\
        \left. \gamma_{\nu} \right|_{t = 0} = \left(\gamma_0\right)_{\nu}, \, \left. \gamma_{\nu} \right|_{t = T} = \left(\gamma_T\right)_{\nu}
    \end{array} \right\} ,
\end{equation*}
where instead of the standard continuity equation, we have $\frac{\partial \gamma_{\nu}}{\partial t} = f_{\nu}$ to ensure mass conservation in the vertex. \\

\begin{observation}[$p$-Wasserstein distance at interior vertex] ~ \\
    The application of the dynamic formulation of the $p$-Wasserstein distance to an interior vertex $\nu \in \mathring{\mathcal{V}}$ is not particularly interesting as $\left(\gamma_0\right)_{\nu}, \left(\gamma_T\right)_{\nu} \in \mathbb{P} \left(\left\{\nu\right\}\right)$ directly implies
    \begin{equation*}
        \left(\gamma_0\right)_{\nu} = \left(\gamma_T\right)_{\nu} = \delta_{\nu}
    \end{equation*}
    and as $h\left(f_{\nu}\left(t\right), \gamma_{\nu}\left(t\right)\right) \geq 0$ for all $\gamma_{\nu} \in \mathcal{M}_+ \left(\left\{\nu\right\} \times \left[0, T\right]\right)$ and $f_{\nu} \in \mathcal{M} \left(\left\{\nu\right\} \times \left[0, T\right]\right)$, the unique strong solution of the dynamic formulation of the $p$-Wasserstein distance is given by
    \begin{equation*}
        \gamma_{\nu} \left(t\right) \equiv \left(\gamma_0\right)_{\nu} = \delta_{\nu} = \left(\gamma_T\right)_{\nu} \qquad \text{and} \qquad f_{\nu} \left(t\right) \equiv 0,
    \end{equation*}
    because $f_{\nu} \left(t\right) \neq 0$ for some $t \in \left[0, T\right]$ would lead to 
    \begin{equation*}
        h\left(f_{\nu}\left(t\right), \gamma_{\nu}\left(t\right)\right) \in \mathbb{R}_{> 0} ~ \text{for} ~ \gamma_{\nu}\left(t\right) > 0 \qquad \text{and} \qquad h\left(f_{\nu}\left(t\right), \gamma_{\nu}\left(t\right)\right) = \infty ~ \text{for} ~ \gamma_{\nu}\left(t\right) = 0.
    \end{equation*}
    
     \noindent The continuity of $t \mapsto \gamma_{\nu} \left(t\right)$ and $t \mapsto f_{\nu}$ thus would lead to $T^{p - 1} \int_0^T h\left(f_{\nu}\left(t\right), \gamma_{\nu}\left(t\right)\right) \, \mathrm d t > 0$ if there exists any $t \in \left[0, T\right]$ with $f_{\nu} \left(t\right) \neq 0$.
\end{observation}

\noindent By extending the previous remarks about the dynamic formulation of the $p$-Wasserstein distance on a single edge or an interior vertex, we can define the $p$-Wasserstein metric on general oriented graphs modeling network structures. \\

\begin{definition}[$p$-Wasserstein distance on networks without boundary vertices] ~ \\
    For an oriented graph modeling a network, the dynamic formulation of the $p$-Wasserstein distance, depending on the coupling condition at interior vertices, is given by:
    \begin{itemize}
        \item Classical Kirchhoff's law and no boundary vertices, with
        $\mu_0 := \rho_0$ and $\mu_T := \rho_T$:
        \begin{flalign*}
            & W_p^p \left(\mu_0, \mu_T\right) = \min_{\substack{\rho \in \mathcal{M}_+ \left(\mathcal{E}\right), \\ j \in \mathcal{M} \left(\mathcal{E}\right)}} \left\{T^{p - 1} \int_0^T \sum_{e \in \mathcal{E}} \left(\int_0^{L_e} h\left(j_e\left(x, t\right), \rho_e\left(x, t\right)\right) \, \mathrm d \eta\right) \, \mathrm d t ~  
            \middle| \begin{array}{l}
                \forall e \in \mathcal{E}: \, 1), \, 2) \vspace{0.1cm} \\
                \forall \nu \in \mathcal{V}: \, 3 \, i)
            \end{array} \right\}
        \end{flalign*}
        \item Generalized Kirchhoff's law and no boundary vertices, with $\mu_0 := \left(\rho_0, \gamma_0\right)$ and $\mu_T := \left(\rho_T, \gamma_T\right)$:
        \begin{flalign*}
            W_p^p \left(\mu_0, \mu_T\right) = \min_{\substack{\rho \in \mathcal{M}_+ \left(\mathcal{E}\right), \, j \in \mathcal{M} \left(\mathcal{E}\right), \\ \gamma \in \mathcal{M}_+ \left(\mathcal{V}\right), \, f \in \mathcal{M} \left(\mathcal{V}\right)}} \Biggl\{ & T^{p - 1} \int_0^T \sum_{e \in \mathcal{E}} \left(\int_0^{L_e} h\left(j_e\left(x, t\right), \rho_e\left(x, t\right)\right) \, \mathrm d \eta\right) \\
            & + \sum_{\nu \in \mathcal{V}} h\left(f_{\nu}\left(t\right), \gamma_{\nu}\left(t\right)\right) \, \mathrm d t 
            \left. \middle| ~ \begin{array}{l}
                \forall e \in \mathcal{E} : \, 1), \, 2) \vspace{0.1cm} \\
                \forall \nu \in \mathcal{V} : \, 3 \, ii), \, 4), \, 5)
            \end{array} \right\}
        \end{flalign*}
    \end{itemize}
    
    \noindent The respective side constraints are given by ($\eta$ again being a reference measure): \\
    
    \begin{minipage}{0.5\textwidth}
        \begin{itemize}
            \item[1)] $\frac{\partial \rho_e}{\partial t} + \frac{\partial j_e}{\partial x} = 0$
            \item[2)] $\left. \rho_e \right|_{t = 0} = \left(\rho_0\right)_e, \, \left. \rho_e \right|_{t = T} = \left(\rho_T\right)_e$
            \item[3 i)] $0 = \sum_{\substack{e \in \mathcal{E}: \\ \delta^E\left(e\right) = \nu}} \left. j_e \right|_{x = L_e} - \sum_{\substack{e \in \mathcal{E}: \\ \delta^S\left(e\right) = \nu}} \left. j_e \right|_{x = 0}$
        \end{itemize}
    \end{minipage}
    \begin{minipage}{0.5\textwidth}
        \begin{itemize}
            \item[3 ii)] $f_{\nu} = \sum_{\substack{e \in \mathcal{E}: \\ \delta^E\left(e\right) = \nu}} \left. j_e \right|_{x = L_e} - \sum_{\substack{e \in \mathcal{E}: \\ \delta^S\left(e\right) = \nu}} \left. j_e \right|_{x = 0}$
            \item[4)] $\frac{\partial \gamma_{\nu}}{\partial t} = f_{\nu}$
            \item[5)] $\left. \gamma_{\nu} \right|_{t = 0} = \left(\gamma_0\right)_{\nu}, \, \left. \gamma_{\nu} \right|_{t = T} = \left(\gamma_T\right)_{\nu}$
        \end{itemize}
    \end{minipage}
\end{definition}

\noindent \af{This definition of the $p$-Wasserstein distance allows for a specific type of unbalanced optimal transport, as mass can be moved from edges (the network) to vertices (storage opportunity), and vice versa, for additional costs. This concept of unbalanced optimal transport can be further generalized to also include boundary vertices and therefore the possibility to move gas from interior vertices and edges to boundary vertices (supply and demand).} \\

\begin{remark}[$p$-Wasserstein distance on networks with boundary vertices] ~ \\
    For an oriented graph with boundary vertices, the dynamic formulation of the $p$-Wasserstein distance, depending on the coupling condition at interior vertices, is given by:
    \begin{itemize}
        \item Classical Kirchhoff's law and time-dependent or time-independent boundary conditions, with $\mu_0 := \left(\rho_0, S_0. D_0\right)$ and $\mu_T := \left(\rho_T, S_T, D_T\right)$:
        \begin{flalign*}
            & W_p^p \left(\mu_0, \mu_T\right) = \min_{\substack{\rho \in \mathcal{M}_+ \left(\mathcal{E}\right), \, j \in \mathcal{M} \left(\mathcal{E}\right), \\ S, -s \in \mathcal{M}_+ \left(\partial^+ \mathcal{V}\right), \\ D, d \in \mathcal{M}_+ \left(\partial^- \mathcal{V}\right)}} \Biggl\{T^{p - 1} \int_0^T \sum_{e \in \mathcal{E}} \left(\int_0^{L_e} h\left(j_e\left(x, t\right), \rho_e\left(x, t\right)\right) \, \mathrm d \eta\right) \\
            & + \sum_{\nu \in \partial^+ \mathcal{V}} h\left(s_{\nu}\left(t\right), S_{\nu}\left(t\right)\right) + \sum_{\nu \in \partial^- \mathcal{V}} h\left(d_{\nu}\left(t\right), D_{\nu}\left(t\right)\right) \, \mathrm d t ~
            \left. \middle| ~ 
            \begin{array}{l}
                \forall e \in \mathcal{E} : \, 1), \, 2) \vspace{0.1cm} \\
                \forall \nu \in \mathring{\mathcal{V}}: \, 3 \, i) \vspace{0.1cm} \\
                \forall \nu \in \partial^+ \mathcal{V}: \, 6), \, 7) , \, 8) \vspace{0.1cm} \\
                \forall \nu \in \partial^- \mathcal{V} : \, 9), \, 10), \, 11)
            \end{array} \right\} .
        \end{flalign*}
        \item Generalized Kirchhoff's law and time-dependent or time-independent boundary conditions, with $\mu_0 := \left(\rho_0, \gamma_0, S_0, D_0\right)$ and $\mu_T := \left(\rho_T, \gamma_T, S_T, D_T\right)$:
        \begin{flalign*}
            & W_p^p \left(\mu_0, \mu_T\right) = \\
            & \min_{\substack{\rho \in \mathcal{M}_+ \left(\mathcal{E}\right), \, j \in \mathcal{M} \left(\mathcal{E}\right), \\ \gamma \in \mathcal{M}_+ \left(\mathring{\mathcal{V}}\right), \, f \in \mathcal{M} \left(\mathring{\mathcal{V}}\right), \\ S, -s \in \mathcal{M}_+ \left(\partial^+ \mathcal{V}\right), \\ D, d \in \mathcal{M}_+ \left(\partial^- \mathcal{V}\right)}} \Biggl\{T^{p - 1} \int_0^T \sum_{e \in \mathcal{E}} \left(\int_0^{L_e} h\left(j_e\left(x, t\right), \rho_e\left(x, t\right)\right) \, \mathrm d \eta\right) + \sum_{\nu \in \mathring{\mathcal{V}}} h\left(f_{\nu}\left(t\right), \gamma_{\nu}\left(t\right)\right) \\
            & + \sum_{\nu \in \partial^+ \mathcal{V}} h\left(s_{\nu}\left(t\right), S_{\nu}\left(t\right)\right) + \sum_{\nu \in \partial^- \mathcal{V}} h\left(d_{\nu}\left(t\right), D_{\nu}\left(t\right)\right) \, \mathrm d t ~
            \left. \middle| ~ \begin{array}{l}
                \forall e \in \mathcal{E} : \, 1), \, 2) \vspace{0.1cm} \\
                \forall \nu \in \mathring{\mathcal{V}}: \, 3 \, ii), \, 4), \, 5) \vspace{0.1cm} \\
                \forall \nu \in \partial^+ \mathcal{V}: \, 6), \, 7), \, 8) \vspace{0.1cm} \\
                \forall \nu \in \partial^- \mathcal{V} : \, 9), \, 10), \, 11)
            \end{array} \right\} .
        \end{flalign*}
    \end{itemize}
    
    \noindent Here, $\eta$ is a reference measure and the additional constraints are given by: \\

    \begin{minipage}{0.5\textwidth}
        \begin{itemize}
            \item[6)] $s_{\nu} = \sum_{\substack{e \in \mathcal{E}: \\ \delta^E\left(e\right) = \nu}} \left. j_e \right|_{x = L_e} - \sum_{\substack{e \in \mathcal{E}: \\ \delta^S\left(e\right) = \nu}} \left. j_e \right|_{x = 0}$
            \item[7)] $\frac{\partial S_{\nu}}{\partial t} = s_{\nu}$
            \item[8)] $\left. S_{\nu} \right|_{t = 0} = \left\lvert S_{\nu}^G \right\rvert, \, \left. S_{\nu} \right|_{t = T} = 0$
        \end{itemize}
    \end{minipage}
    \begin{minipage}{0.5\textwidth}
        \begin{itemize}
            \item[9)] $d_{\nu} = \sum_{\substack{e \in \mathcal{E}: \\ \delta^E\left(e\right) = \nu}} \left. j_e \right|_{x = L_e} - \sum_{\substack{e \in \mathcal{E}: \\ \delta^S\left(e\right) = \nu}} \left. j_e \right|_{x = 0}$
            \item[10)] $\frac{\partial D_{\nu}}{\partial t} = d_{\nu}$
            \item[11)] $\left. D_{\nu} \right|_{t = 0} = 0, \, \left. D_{\nu} \right|_{t = T} = D_{\nu}^G$
        \end{itemize}
    \end{minipage}
\end{remark}

%% file: 4.3_TODO.tex
\subsection{Characterization of absolutely continuous curves}

Before investigating the connection between $p$-absolutely continuous curves and velocity fields in this subsections, we need to firstly introduce some further concepts and notations. \\

\begin{definition}[$p$-absolutely continuous curves \cite{ambrosio2005gradient}, \cite{erbar2021gradient}] ~ \\
    For the complete metric space $\left(\mathcal{X}, \tilde{d}\right)$, we call a curve
    \begin{equation*}
        \mu: \left(0, T\right) \longrightarrow \mathcal{X} \quad t \mapsto \mu\left(t\right),
    \end{equation*}
    $p$-absolutely continuous for $p \in \left[1, \infty\right)$, if there exists a function $g \in L^p \left(\left(0, T\right), \mathbb{R}\right)$ such that
    \begin{equation}
        \label{abscont}
        \forall 0 < r \leq s < T: \, \tilde{d} \left(\mu \left(r\right), \mu \left(s\right)\right) \leq \int_r^s g \left(t\right) \, \mathrm d t .
    \end{equation}
    Here, $L^p \left(\left(0, T\right), \mathbb{R}\right)$ denotes the $L^p$-space of $\mathcal{L}$-measurable and $\mathbb{R}$-valued functions on the interval $\left(0, T\right)$.
%
    By $AC^p \left(\left(0, T\right), \mathcal{X}\right)$ we will denote the corresponding space of $p$-absolutely continuous curves on the complete metric space $\left(\mathcal{X}, \tilde{d}\right)$. \\
\end{definition}

\noindent Among all possible function choices of $g$ in (\ref{abscont}), we are interested in the minimal one, whose existence is given by the following proposition. \\

\begin{proposition}[Metric derivative of $p$-absolutely continuous curves \cite{ambrosio2005gradient}] ~ \\
    For every $p$-absolutely continuous curve $\mu \in AC^p \left(\left(0, T\right), \mathcal{X}\right)$, the limit
    \begin{equation*}
        \left\lvert \mu' \left(t\right) \right\rvert := \lim_{r \rightarrow s} \frac{\tilde{d} \left(\mu\left(r\right), \mu\left(s\right)\right)}{\left\lvert r - s \right\rvert}
    \end{equation*}
    exists for $\mathcal{L}$-a.e. $s \in \left(0, T\right)$ and the function $t \mapsto \left\lvert \mu' \left(t\right) \right\rvert$ belongs to $L^p \left(\left(0, T\right), \mathbb{R}\right)$. Furthermore, $\left\lvert \mu' \left(t\right) \right\rvert$ is an admissible choice for $g$ and minimal in the sense of
    \begin{equation*}
        \left\lvert \mu' \left(t\right) \right\rvert \leq \tilde{g}\left(t\right) ~ \text{for} ~ \mathcal{L} \text{-a.e.} ~ t \in \left(0, T\right)
    \end{equation*}
    for each function $\tilde{g} \in L^p \left(\left(0, T\right), \mathbb{R}\right)$ satisfying (\ref{abscont}).
\end{proposition}

\begin{proof}
    See for example Theorem 1.1.2 in \cite{ambrosio2005gradient}.
\end{proof}

\noindent As we are interested in the previously introduced network setting, we will set $\mathcal{X} := \mathcal{M}_+ \left(\Omega_{\mathcal{G}}\right)$ (or for simpler examples we can also set $\mathcal{X} := \mathcal{M}_+ \left(\left[0, L_e\right]\right)$ for an edge $e \in \mathcal{E}$, or $\mathcal{X} := \mathcal{M}_+ \left(\left\{\nu\right\}\right)$ for a vertex $\nu \in \mathcal{V}$) and $\tilde{d} := W_p$. \\

\begin{theorem}[Connection $p$-absolutely continuous curves and velocity fields] ~ \\
    Let
    \begin{equation*}
        \mu: \left(0, T\right) \longrightarrow \mathcal{M}_+ \left(\Omega_{\mathcal{G}}\right) \quad t \mapsto \mu\left(t\right)
    \end{equation*}
    be a curve and 
    \begin{equation*}
        V: \Omega_{\mathcal{G}} \times \left(0, T\right) \longrightarrow \mathbb{R}^{k + 1} \quad \left(x, t\right) \mapsto V\left(x, t\right)
    \end{equation*}
    a vector field, then the following two statements hold true:
    \begin{itemize}
        \item[1)] $\mu \in AC^p \left(\left(0, T\right), \mathcal{M}_+ \left(\Omega_{\mathcal{G}}\right)\right) \quad \Longrightarrow \quad$ for $\mathcal{L}$-a.e. $\tilde{t} \in \left(0, T\right): \, \exists \left. V \right|_{t = \tilde{t}} \in L^p\left(\Omega_{\mathcal{G}}, \mathbb{R}^{k + 1}, \mu\right): \, \left(\int_{\Omega_{\mathcal{G}}} \left\lvert V \left(x, \tilde{t}\right) \right\rvert^p \, \mathrm d \mu \left(x, \tilde{t}\right) \right)^{\frac{1}{p}} \leq \left\lvert \mu' \left(\tilde{t}\right) \right\rvert$ and $\left(\mu, \mu V\right)$ is a weak solution of (\ref{CE 1})
        \item[2)] $\left(\mu, \mu V\right)$ is a weak solution of (\ref{CE 1}) with $\int_0^T \left(\int_{\Omega_{\mathcal{G}}} \left\lvert V \left(x, t\right) \right\rvert^p \, \mathrm d \mu \left(x,t\right) \right)^{\frac{1}{p}} \, \mathrm d t < \infty \quad \Longrightarrow \quad \mu \in AC^p \left(\left(0, T\right), \mathcal{M}_+ \left(\Omega_{\mathcal{G}}\right)\right)$ and for $\mathcal{L}$-a.e. $\tilde{t} \in \left(0, T\right): \, \left\lvert \mu' \left(\tilde{t}\right) \right\rvert \leq \left(\int_{\Omega_{\mathcal{G}}} \left\lvert V \left(x, \tilde{t}\right) \right\rvert^p \, \mathrm d \mu \left(x, \tilde{t}\right) \right)^{\frac{1}{p}}$
    \end{itemize}
    where $L^p\left(\Omega_{\mathcal{G}}, \mathbb{R}^{k + 1}, \mu\right)$ denotes the $L^p$-space $L^p\left(\Omega_{\mathcal{G}}, \mathbb{R}^{k + 1}\right)$ with respect to the measure $\mu$.
\end{theorem}

\begin{proof}
    For the case of $p = 2$ this theorem has been proven in \cite{erbar2021gradient}, theorem 3.8. A proof for the presented more generalized framework is postponed to future work.
\end{proof}

%% file: 6_gradient_flows.tex
\section{$p$-Wasserstein gradient flows} \label{sec:gf}

We finally have a closer look at gradient flows in the $p$-Wasserstein metric on networks. As usual in the metric setting (cf. \cite{ambrosio2005gradient}), they can be derived in the limit $\tau \rightarrow 0$, for $\tau > 0$, from the minimizing movement scheme 
\begin{equation} \label{eq:minimizingmovements}
    \mu^\tau_{\left(l + 1\right) \tau} = \text{arg}\min_{\mu} \left( 
    \frac{1}{p \tau^{p - 1}} d \left(\mu, \mu^{\tau}_{l \tau}\right)^p + E \left(\mu\right) \right)
\end{equation}
with appropriate (variational) interpolation. The analysis of gradient flows in such a setting for the standard Wasserstein metric has already been considered in \cite{agueh,ambrosio2005gradient}. \\
 
\noindent In this section we will build on the rigorous definitions from the previous parts, but mainly take a formal approach to clarify two issues: one is the derivation of consistent interface conditions for the gradient flow, the second is to rewrite an existing standard model as a gradient flow on the network. 

\subsection{$p$-Wasserstein gradient flow on networks}

In the following, let us derive a more concrete structure for Wasserstein gradient flows on networks (for the case of classical Kirchhoff's law at interior vertices and no boundary vertices). For this sake we consider energies of the natural form
\begin{equation*}
    E(\mu) = \sum_{e\in \mathcal{E}} E_e \left(\rho_e\right),
\end{equation*}
i.e., an additive composition of functionals on each edge $e \in \mathcal{E}$. \\

\noindent As in the case of Wasserstein gradient flows in domains, we expect that 
\begin{flalign} \label{eq:gfp1}
    \frac{\partial \rho_e}{\partial t} + \frac{\partial j_e}{\partial x} & = 0, \\ \label{eq:gfp2}
    j_e \left\lvert j_e \right\rvert^{p - 2} & = - \rho_e^{p-1} \frac{\partial}{\partial x} \left(E'\left(\rho_e\right)\right),
\end{flalign}
together with the classical Kirchhoff condition for all $\nu \in \mathcal{V}$
\begin{equation} \label{eq:kirchhoff}
  0 = \sum_{\substack{e \in \mathcal{E}: \\ \delta^E\left(e\right) = \nu}} \left. j_e \right|_{x = L_e} - \sum_{\substack{e \in \mathcal{E}: \\ \delta^S\left(e\right) = \nu}} \left. j_e \right|_{x = 0}.
\end{equation}
However, to complete the system, we need additional coupling conditions, which will arise automatically from the variational formulation. \\

\noindent To see this, let us derive the Hamilton-Jacobi equation for the dual variable in the Wasserstein metric, which will  be related to the first variation of the energy. A simple formal derivation can be obtained from assuming the existence of a saddle-point of the associated Lagrange functional, whose first variations we can be set to zero
\begin{equation*}
   \frac{1}{p \tau^{p-1}}W_p^p \left(\left(\rho_t\right)_e, \left(\rho_{t+\tau}\right)_e\right) =
    \min_{\rho_e  ,j_e }  \max_{\phi_e}
    \sum_{e \in \mathcal{E}} \int_t^{t+\tau} \int_0^{L_e} 
   \left( \frac{\left\lvert j_e \right\rvert^p}{p \rho_e^{p-1}} +  \frac{\partial \rho_e}{\partial t} \phi_e + \frac{\partial j_e}{\partial x} \phi_e \right) \, \mathrm d x \, \mathrm d t 
\end{equation*}
for $\tau > 0$, $t \in \left[0, T - \tau\right]$ and test functions $\phi_e$, which, in case of the test function being a maximizing argument, correspond to the Lagrange multipliers. Here, the minimization is carried out among those $\rho_e, j_e$, which fulfill the given initial and terminal values and satisfy \eqref{eq:kirchhoff}. Thus, the variation with respect to $j_e$ in direction $k_e$ yields 
\begin{equation*}
    \sum_{e \in \mathcal{E}} \int_t^{t+\tau} \int_0^{L_e} \left( \frac{j_e \left\lvert j_e\right\rvert^{p - 2}}{\rho_e^{p-1}} k_e + \frac{\partial k_e}{\partial x} \phi_e \right) \, \mathrm d x \, \mathrm d t = 0.
\end{equation*}
Now, let us integrate by parts and use that each feasible variation $k_e$ still needs to satisfy Kirchhoff's law, which yields
\begin{align*}  0 =& \sum_{\nu \in \mathcal{V}} \left(\sum_{\substack{e \in \mathcal{E}: \\ \delta^E\left(e\right) = \nu}} \left. (k_e \phi_e) \right|_{x = L_e} - \sum_{\substack{e \in \mathcal{E}: \\ \delta^S\left(e\right) = \nu}} \left. (k_e \phi_e) \right|_{x = 0}\right) + \\&
\sum_{e \in \mathcal{E}}\int_t^{t+\tau} \int_0^{L_e} 
   \left( \frac{j_e |j_e|^{p-2}}{ \rho_e^{p-1}}  -\frac{\partial \phi_e}{\partial x}  \right) k_e\, \mathrm d x \, \mathrm d t \\
   =& \sum_{\nu \in \mathcal{V}} \left(\sum_{\substack{e \in \mathcal{E}: \\ \delta^E\left(e\right) = \nu}} \left. (k_e (\phi_e-\Phi_\nu)) \right|_{x = L_e} - \sum_{\substack{e \in \mathcal{E}: \\ \delta^S\left(e\right) = \nu}} \left. (k_e (\phi_e-\Phi_\nu)) \right|_{x = 0}\right) + \\&
\sum_{e \in \mathcal{E}}\int_t^{t+\tau} \int_0^{L_e} 
   \left( \frac{j_e |j_e|^{p-2}}{ \rho_e^{p-1}}  -\frac{\partial \phi_e}{\partial x}  \right) k_e\, \mathrm d x \, \mathrm d t
\end{align*}
with the nodal mean
\begin{equation*}
    \Phi_{\nu} := \frac{1}{\text{deg}(\nu)} \left(\sum_{\substack{e \in \mathcal{E}: \\ \delta^E\left(e\right) = \nu}} \phi_e +   \sum_{\substack{e \in \mathcal{E}: \\ \delta^S\left(e\right) = \nu}} \phi_e \right),
\end{equation*}
\af{and the degree function
\begin{equation*}
    \text{deg} : \mathcal{V} \longrightarrow \mathbb{N}_0 \quad \nu \mapsto \text{deg}(\nu) := \left\lvert\left\{e \in \mathcal{E} ~ \middle| ~ \nu \in e\right\}\right\rvert,
\end{equation*}
where for a connected graph it holds true that $\text{deg}\left(\nu\right) \geq 1$ for all vertices $\nu \in \mathcal{V}$.} \\
   
\noindent Choosing arbitrary fluxes $k_e$ with compact support inside a single edge immediately yields
\begin{equation*}
    j_e \left\lvert j_e \right\rvert^{p - 2} = \rho_e^{p - 1} \frac{\partial \phi_e}{\partial x} 
\end{equation*}
for each $e \in \mathcal{E}$. Choosing further arbitrary $k_e$ with support only close to a single vertex $\nu \in \mathcal{V}$ (up to satisfying \eqref{eq:kirchhoff}), we see that further the condition
\begin{equation*}
    \phi_e = \Phi_\nu \qquad \forall e \in \mathcal{E} ~ \text{with} ~ \delta^S\left(e\right) = \nu \text{ or } \delta^E\left(e\right) = \nu
\end{equation*}
holds for the vertex $\nu$. \af{For dead-end vertices, i.e. deg$\left(\nu\right) = 1$, this yields no additional conditions, while for all other vertices this implies continuity of the dual variable.} \\

\noindent To relate the dual variable to gradient flows, we can compute the variation of $\frac{1}{p \tau^{p-1}}W_p^p \left(\left(\rho_t\right)_e, \left(\rho_{t+\tau}\right)_e\right)$ with respect to $\left(\rho_{t+\tau}\right)_e$, which is indeed given by $\phi_e\left(t+\tau\right)$ as a standard computation analogous to the Wasserstein distance on domains, it shows. Thus, the optimality condition in \eqref{eq:minimizingmovements} becomes
\begin{equation*}
    \phi_e \left(t+\tau\right) + E_e'\left(\rho_e\left(t+\tau\right)\right) = 0. 
\end{equation*}

\noindent In the limit $\tau \downarrow 0$, we hence obtain the gradient flow \eqref{eq:gfp1},\eqref{eq:gfp2} supplemented by \eqref{eq:kirchhoff}
and the following interface conditions: there exists a nodal function $\Phi: \mathcal{V} \rightarrow \mathbb{R}$ such that
\begin{equation} \label{eq:interface}
   E_e'(\rho_e)=\Phi_\nu \qquad \forall e \in \mathcal{E} ~ \text{with} ~ \delta^S\left(e\right) = \nu \text{ or } \delta^E\left(e\right) = \nu,
\end{equation}
for vertices $\nu \in \mathcal{V}$.

\subsection{The (ISO3) model for gas networks}

In the following, we turn our attention to the (\ref{eqn:ISO3}) model of transport in gas networks, which is derived and embedded in a full model hierarchy in \cite{domschke2021gas} (see also Section \ref{gasflowpipe}). Based on mass density $\rho_e$ and velocity $v_e$ functions defined on edges, the (IS03) model is given by
\begin{flalign*}
    \frac{\partial \rho_e}{\partial t} + \frac{\partial  }{\partial x} ( \rho_e v_e) & = 0, \\
    \frac{\lambda_e}{2 \mathcal{D}_e} \rho_e v_e \left\lvert v_e \right\rvert & = - \frac{\partial}{\partial x} \left(p_e\left(\rho_e\right)\right) - g \rho_e \sin \left(\omega_e\right),
\end{flalign*}
together with the classical Kirchhoff condition (\ref{eqn:C KL}) on the vertices of the network. \\

\noindent We start by reformulating the model in terms of the mass density variable $\rho_e$ and the flux variable $j_e$ as:
\begin{flalign*}
    \frac{\partial \rho_e}{\partial t} + \frac{\partial j_e}{\partial x} & = 0 \\
    j_e \left\lvert j_e \right\rvert & = - \frac{2 \mathcal{D}_e}{\lambda_e} \rho_e \frac{\partial}{\partial x} \left(p_e\left(\rho_e\right)\right) - \rho_e^2 \frac{2 \mathcal{D}_e g}{\lambda_e} \sin \left(\omega_e\right) \\
    &= - \rho_e^2 \left( \frac{2 \mathcal{D}_e}{\lambda_e \rho_e }\frac{\partial}{\partial x} \left(p_e\left(\rho_e\right)\right) + \frac{2 \mathcal{D}_e g}{\lambda_e} \sin \left(\omega_e\right)\right)
\end{flalign*}
Specifying the general form of the Wasserstein gradient flow to $p = 3$ implies 
\begin{equation*}
    \frac{\partial}{\partial x} E'(\mu) = \frac{2 \mathcal{D}_e}{\lambda_e \rho_e }\frac{\partial}{\partial x} \left(p_e\left(\rho_e\right)\right) + \frac{2 \mathcal{D}_e g}{\lambda_e} \sin \left(\omega_e\right)
\end{equation*}
on each edge $e \in \mathcal{E}$. As usual, we need to introduce an entropy functional, which  - due to the dependence on diameter $\mathcal{D}_e$ and friction coefficient $\lambda_e$ - may vary with the edge. The derivative of the entropy is defined as  
\begin{equation*}
    F_e''(s) := \frac{2 \mathcal{D}_e}{\lambda_e} \frac{p_e'\left(s\right)}s,
\end{equation*}
and we define $F_e$ as a second primitive with $F_e \left(0\right) = 0$. Then we see that
\begin{equation*}
    \frac{\partial}{\partial x} F_e'\left(\rho_e\right) = \frac{2 \mathcal{D}_e}{\lambda_e \rho_e }\frac{\partial}{\partial x} \left(p_e\left(\rho_e\right)\right)
\end{equation*}
holds. With the short-hand notation 
\begin{equation*}
    c_e:=\frac{2 \mathcal{D}_e g}{\lambda_e} \sin \left(\omega_e\right),
\end{equation*}
we obtain the energy functional
\begin{equation*}
    E(\mu) = \sum_ e\int_0^{L_e} \left(F_e\left(\rho_e\right) +  c_e x \rho_e + d_e \rho_e\right) \, \mathrm d x.
\end{equation*}
Note that $d_e$ is an arbitrary integration constant, which however impacts the gradient flow as soon as there is more than one edge, since the mass per edge is not necessarily conserved then. Since a global change by a constant will not affect the result, due to the overall mass conservation on the network, we can however use the condition
\begin{equation} \label{eq:dnormalization}
    \sum_{e \in \mathcal{E}} d_e = 0.
\end{equation}

\noindent As the derivation above shows, the gradient flow formulation immediately yields the appropriate interface conditions \eqref{eq:interface}, in addition to the classical Kirchhoff's law  for \eqref{eqn:ISO3}, which are given by
\begin{align*} 
   F_e'(\rho_e) + d_e&=\Phi_\nu \qquad \forall e \in \mathcal{E} ~ \text{with} ~ \delta^S\left(e\right) = \nu, \\
   F_e'(\rho_e)+ c_e L_e + d_e \ &=\Phi_\nu \qquad \forall e \in \mathcal{E} ~ \text{with} ~ \delta^E\left(e\right) = \nu,
\end{align*}
for vertices $\nu \in \mathcal{V}$. \\

\noindent The interface condition effectively means that the potential $F_e'(\rho_e) + c_e x + d_e$ is constant across vertices. This may seem quite arbitrary, since a potential is usually specified only up to a constant. However, we have to understand the potential rather as a global quantity on the network, so clearly changes on single edges affect the global structure of the potential. Since there is quite some freedom to determine the constants $d_e$, we could in turn specify certain interface conditions and determine the associated potential. Again, it is reasonable that different jump conditions on the vertices have the same effect as changing a potential. 

\subsection{Vanishing Diffusion Limit}

As mentioned above, there are several potentials and thus interface conditions that lead to gradient flows. We may want to specify however a potential such that in the limit of vanishing diffusion ($F_e=0$), we obtain a consistent system. The problem of vanishing diffusion has been investigated for linear problems on networks in \cite{eggerphilippi} with a focus on incompressibility rather than a gradient structure. \\

\noindent Consistency of the vector field with the vanishing diffusion limit means
\begin{align}  \label{eq:dequationS}
     d_e&=\Phi_\nu \qquad \forall e \in \mathcal{E} ~ \text{with} ~ \delta^S\left(e\right) = \nu, \\ \label{eq:dequationE}
    c_e L_e + d_e \ &=\Phi_\nu \qquad \forall e \in \mathcal{E} ~ \text{with} ~ \delta^E\left(e\right) = \nu,
\end{align}
for $\nu \in \mathcal{V}$, which - together with \eqref{eq:dnormalization} - can be considered a system of linear equations for the variables $d_e$, or can even be extended to a linear system for the variables $d_e$ and $\Phi_\nu$. Hence, it is natural to investigate the (unique) solvability of this system. \\

\noindent The number of unknowns is $\left\lvert \mathcal{E}\right\rvert + \left\lvert\mathcal{V}\right\rvert$, while the number of equations is $\sum_{\nu \in \mathcal{V}}$ deg$\left(\nu\right) + 1$. By the sum formula \cite[Theorem 1.1]{bondy}, we have $\sum_{\nu \in \mathcal{V}}$ deg$\left(\nu\right) + 1 = 2 \left\lvert\mathcal{E}\right\rvert + 1$. Thus, the number of equations and unknowns coincides if $\left\lvert\mathcal{E}\right\rvert + 1 = \left\lvert \mathcal{V}\right\rvert$.\\

\noindent \af{It is a simple exercise to see that $\left\lvert \mathcal{E}\right\rvert \geq \left\lvert \mathcal{V}\right\rvert$ implies the existence of a cycle (cf. \cite[p.15]{bondy}) in a so-called simple graph (a graph without loops and parallel edges). Hence, in a simple connected graph without a cycle (which seems a suitable condition for gas networks), we have $\left\lvert \mathcal{E} \right\rvert + 1 \leq \left\lvert \mathcal{V} \right\rvert$, which is however also an upper bound for the number of edges.} Thus, in this case, the number of equations and unknowns coincides and we can see that the linear system has full rank.

\begin{proposition}
    Let the graph $\mathcal{G} = \left(\mathcal{V},\mathcal{E}\right)$ be simple, connected, and contain no cycles. Then there exists a unique solution $(d_e,\Phi_\nu)$ of the linear system \eqref{eq:dnormalization}, \eqref{eq:dequationS} and \eqref{eq:dequationE}. 
\end{proposition}
\begin{proof}
    Since under the above conditions on the graph, we have verified that the number of equations and unknowns in the linear system coincides, we only need to show that the homogeneous problem (corresponding to $c_e = 0$) has a unique trivial solution. \\
    
    \noindent We first notice that in the homogeneous case, the orientation of the graph plays no role and that \eqref{eq:dequationS} and \eqref{eq:dequationE} imply that $d_e = d_{e'}$ if $e$ and $e'$ are adjacent to a common vertex. Since the graph is connected, there exists a path that connects two arbitrary edges, hence by transitivity we obtain $d_e = d_{e'}$ for any two edges $e,e' \in \mathcal{E}$. \\
    
    \noindent Since all $d_e$ are equal, \eqref{eq:dnormalization} implies $d_e = 0$ for all $e \in \mathcal{E}$, which further implies $\Phi_\nu = 0$ for all $\nu \in \mathcal{V}$, since each vertex is connected to at least on edge.
\end{proof}

\noindent We finally illustrate the choice of the constants with a simple example:

\begin{example} ~ \\
    Let us consider a very simple tree network consisting of the vertices $\mathcal{V} = \left\{\nu_1, \nu_2, \nu_3, \nu_4\right\}$, and the edges $\mathcal{E} = \left\{e_1, e_2, e_3\right\}$ with
    \begin{align*}
        & \delta^S\left(e_1\right) = \nu_1, \qquad \delta^S\left(e_2\right) = \nu_2, \qquad \delta^S\left(e_3\right) = \nu_2, \\
        & \delta^E\left(e_1\right) = \nu_2, \qquad \delta^E\left(e_2\right) = \nu_3, \qquad \delta^E\left(e_3\right) = \nu_4.
    \end{align*}
    Thus, we obtain the following linear system: \\

    \begin{minipage}{.3\textwidth}
        \begin{center}
            \begin{tikzpicture}[scale=0.6, transform shape]
    		\pgfdeclarelayer{background}
    		\pgfsetlayers{background,main}
    		\begin{scope}[every node/.style={draw,circle}]
                    \node (v1) at (0,0) {$\nu_1$};
        	    \node (v2) at (3,0) {$\nu_2$};
                    \node (v3) at (6,1.5) {$\nu_3$};
                    \node (v4) at (6,-1.5) {$\nu_4$};
    		\end{scope}
                \begin{scope}[>={Stealth[bluegray]},
    				every edge/.style={draw=bluegray,line width=1pt}]
    				\path [-] (v1) edge node[bluegray,below] {$e_1$} (v2);
                    \path [-] (v2) edge node[bluegray,above left] {$e_2$} (v3);
                    \path [-] (v2) edge node[bluegray,below left] {$e_3$} (v4);
    		\end{scope}
            \end{tikzpicture}
        \end{center}
    \end{minipage}
    \begin{minipage}{.6\textwidth}
        \begin{align*}
            d_1 - \Phi_1 &= 0 &&
            d_1 - \Phi_2 = - c_1 L_1 \\
            d_2 - \Phi_2 &= 0 &&
            d_3 - \Phi_2 = 0 \\
            d_2 - \Phi_3 &= - c_2 L_2 &&
            d_3 - \Phi_4 = - c_3 L_3 \\
            d_1+d_2+d_3 &= 0
        \end{align*}
    \end{minipage}

    \vspace{0.5cm}
    
    \noindent Solving this system yields $d_1=-\frac{2}3 c_1 L_1$ and $d_2 = d_3 = \frac{1}3 c_2 L_2$. 
\end{example}

\noindent In the example we see that the constants $d_e$ one edges $e_2$ and $e_3$ only depend on the properties in edge $e_1$, but not in $c_2, c_3, L_2, L_3$. Indeed it is easy to see that independence of the edge parameters is always the case for an edge leading to an outgoing boundary of the network.

%% file: 7_numerics.tex
\section{Numerical examples} \label{sec:num}

This section presents several numerical examples of solutions to the transport problems introduced in section \ref{sec4} above. We examine different types of boundary conditions, always for the case $p=2$. \\

\noindent Our numerical algorithm is based on a variant of the primal-dual optimization algorithm applied to optimal transport that was introduced in \cite{Carrillo2022_PrimalDual}. We work with a slightly different discretziation of the constraints as presented in \cite{pietschmann2022} (and recently extended to metric graphs in \cite{krautz2024}), yet only in the case of no-flux boundary conditions. \\

\noindent The basic idea is to relax the constraints in the optimization problem by only asking them to be satisfied up to a given, fixed precision. The main advantage of this approach is that, when adding the relaxed constraints to the objective functional, this results in a simple proximal operator (i.e. an shrinkage operation) which can be implemented rather efficiently. \\

\noindent For concreteness, we focus on (\ref{TD BC 1}) and (\ref{TD BC 2}) to explain the algorithm in detail. 
\subsection{Discretizing the constraints}
We start by discretizing the interval $[0,L_e]$ associated to every edge into an equidistant grid with $N_x^e$ points. The time interval $[0,T]$ is similarly discretized into an equidistant grid with $N_t$ points. With $\delta_x^e = L_e/N_x^e$ and $\delta_t = T/N_t$ we obtain the grid points 
\begin{align*}
    x_i^e & = (i-1)\delta_x^e \qquad \text{for} ~ e \in E, \, i \in \{1, 2, \ldots, N_x^i+1\}, \\
    t_k & = (k-1)\delta_t \qquad \text{for} ~ k \in \{1, 2, \ldots, N_t+1\}.
\end{align*}
%
Note that in this section, we will slightly deviate from the notation used before, as we will use superscript $e$ and superscript $\nu$ in order to denote the corresponding edge or vertex, respectively, and not to interfere too much with the additional indices coming from the space-time discretization. \\

\noindent To any continuous function $v:[0,L]\times[0,T] \to \mathbb{R}$ we associate a grid function $v_h=(v_{i,k})_{i,k}$ via
\begin{align*}
    v_{i,k}=v(x_i,t_k) \qquad \text{for} ~ i\in\{1,2, \ldots, N_x+1\}, \, k\in\{1, 2, \ldots,N_t+1\}.
\end{align*}
Similar notation is used for grid functions associated to the spatial and temporal partitions, respectively.
On these functions, we define a scalar product which we obtain by discretizing the integrals involved using the trapezoidal rule. Given grid functions $u_h=(\rho_{i,k}^u,m_{i,k}^u)$, $v_h=(\rho_{i,k}^v,m_{i,k}^v)$ we have
\begin{align}\label{eq:def_norm}
    \langle u_h,v_h\rangle = \sum_{k=1}^{N_k+1}\sum_{i=1}^{N_x+1} w_k^tw_j^x\left( \rho_{i,k}^u\rho_{i,k}^v + m_{i,k}^u m_{i,k}^v\right),\quad \|v_h\|=\sqrt{\langle v_h,v_h\rangle},
\end{align}
with 
\begin{align*}
    w_i^x =\begin{cases} \frac{\delta_x}{2} &i\in\{1,N_x+1\},\\
                        \delta_x   &\text{otherwise}.
                        \end{cases}
\end{align*}
We denote by $w_k^t$ a similar weight function for the composite trapezoidal rule associated with the temporal partition, and note that $\sum_{k=1}^{N_t+1}w_k^t=1$.
This renders the linear space of grid functions a Hilbert space. \\

\noindent We will only present the discretization of the constraints for the variant used in (\ref{TD BC 1}) and (\ref{TD BC 2}), with obvious modifications for the other cases. We discretize the continuity equation  \eqref{eqn:CE1} using a  centred differencing in space and a backward differencing in time for the interior grid points $i\in\{2, 3, \ldots N_x\},\ k\in\{2, 3, \ldots,N_t+1\}$, i.e.,
\begin{align}\label{eq:FD_int}
    \frac{\rho^e_{i,k}-\rho^e_{i,k-1}}{\delta_t}+\frac{ j^e_{i+1,k}-j^e_{i-1,k}}{2\delta_x}=0.
\end{align}
For the boundary we use a one-sided finite difference approximation to approximate $\partial_x j^e$, i.e., for $i\in\{1,N_x+1\}$ we use
\begin{align}
    \frac{\rho^e_{1,k}-\rho^e_{1,k-1}}{\delta_t}+\frac{ j^e_{2,k}-j^e_{1,k}}{\delta_x}&=0,\quad\text{and}\label{eq:FD_bdry_left}\\
    \frac{\rho^e_{N_x+1,k}-\rho^e_{N_x+1,k-1}}{\delta_t}+\frac{ j^e_{N_x+1,k}-j^e_{N_x,k}}{\delta_x}&=0.\label{eq:FD_bdry_right}
\end{align}
The coupling conditions (\ref{TD BC 1}) and (\ref{TD BC 2}) at the boundary vertices become
 \begin{alignat*}{3}
        && 0 & = \Bigg( \left\lvert s^{G,\nu}_k  \right\rvert + \sum_{\substack{e \in \mathcal{E}: \\ \delta^E\left(e\right) = \nu}}  j^e_{N_x+1,k}   \Bigg) - \sum_{\substack{e \in \mathcal{E}: \\ \delta^S\left(e\right) = \nu}}  j^e_{1,k}  \qquad \qquad \qquad \quad & \forall \nu \in \partial^+ \mathcal{V}, \\ 
        \qquad && 0 & = \sum_{\substack{e \in \mathcal{E}: \\ \delta^E\left(e\right) = \nu}} j^e_{N_x+1,k}  - \Bigg( d_k^{G,\nu}  + \sum_{\substack{e \in \mathcal{E}: \\ \delta^S\left(e\right) = \nu}} j^e_{1,k} \Bigg) & \forall \nu \in \partial^- \mathcal{V}, 
    \end{alignat*}
for all $k\in\{1, 2, \ldots N_t+1\}$, while the vertex ODEs at the interior vertices, \eqref{eq:vertex_ode},  become
\begin{align}
\frac{\gamma^\nu_{k}-\gamma^\nu_{k-1}}{\delta_t} = f^\nu_k, \quad \forall \nu \in \mathring{\mathcal{V}}.
\end{align}
Since the discretization does not depend on $m_{i, 1}$, we set $m_{i, 1} = 1$ for $i= 1, 2, \ldots N_x + 1$.
Furthermore, the initial- and final conditions of both $\rho^e$ and $\gamma^{\nu}$ are taken into account in the following way

\begin{align}
    \rho^e_{i, 1} = \rho_0^e \left(x_i\right) \quad \text{and} \quad  \rho^e_{i, N_t+1} = \rho_T^e \left(x_i\right) & \qquad \qquad \text{for} ~ i \in \{1, 2, \ldots N_x + 1\}, \, e \in \E, \label{eq:FD_initial} \\
    \gamma^{\nu}_1 = \gamma_0^{\nu} \quad \text{and} \quad \gamma^{\nu}_{N_t + 1} = \gamma_1^{\nu} & \qquad \qquad \text{for} ~ \nu \in \mathring{\V}.
\end{align}
Finally, since the continuity equation is only satisfied approximately, we explicitly add a constraint on the total mass of the graph being a probability density
\begin{align}
    \sum_{e\in\E}\sum_{i=1}^{N_x+1} w_i^x \left(\rho^e_{i,k}-\rho^e_0(x_i)\right) + \sum_{\nu\in\V} \left(\gamma^{\nu}_k - \gamma^{\nu}_0\right)&=0 \qquad \text{for} ~ k\in\{2, 3, \ldots N_t+1\}.\label{eq:FD_mass}
\end{align}
\subsection{Primal-dual algorithm}
We now sketch the main ingredients of the primal-dual algorithm, referring to \cite{Carrillo2022_PrimalDual,pietschmann2022} for details.
As mentioned above, the key idea is to enforce the constraints \eqref{eq:FD_int}--\eqref{eq:FD_mass} only up to a finite precision, i.e. using the norm \eqref{eq:def_norm}, instead of \eqref{eq:FD_int}, we only ask for
\begin{align*}
     \sum_{k=2}^{N_t+1}w_k^t \Bigg( w_1^x \left(\frac{\rho^e_{i,k}-\rho^e_{i,k-1}}{\delta_t}+\frac{ j^e_{i+1,k}-j^e_{i-1,k}}{2\delta_x}\right)^2\Bigg) &\leq \delta_1^2\Bigg. 
\end{align*}
These weakened constraints are quadratic and can be written in the form
\begin{align*} 
   Au \in \mathcal{C}_\delta = \left\{x ~ \middle| ~ \|x_j -b_j\|_2 \leq \delta_j, \, j = 1, \ldots, 8\right\}, 
\end{align*}
where the vector $u$ contains the coefficient of the grid functions $(\rho_h,m_h)$. We have a total of $8$ constraints since we include, in addition to the continuity equation on the edges, the initial- and final data for $\rho^e$, for $\gamma^{\nu}$, the coupling conditions at the vertices, vertex dynamics and an additional mass constraint. Note that the weights $w_i^x$ and $w_k^t$ are included in the definition of $A_j$ and $b_j$, respectively, and the vectors $x_j$ are slices of the vector $x$ corresponding to the number of rows in $A_j$.
We define the matrix $A$ by vertically concatenating the matrices $A_j$ for $j = 1, \ldots, 8$. We note that $A$ is the matrix of a linear map from the Hilbert space of grid functions, with the inner product being defined in \eqref{eq:def_norm}, to a Euclidean space.

\noindent We then aim to solve the fully discretized,  problem
\begin{align*}
    \inf_{(\rho_h,m_h,\gamma_h, f_h)} \left(\sum_{e \in \E}\sum_{k=1}^{N_t+1}\sum_{i=1}^{N_x+1} w^x_i w_k^t h(\rho^e_{i,k},j^e_{i,k}) + \sum_{\nu \in \mathring{\V}}\sum_{k=1}^{N_t+1}h(\gamma^{\nu}_{k},f^{\nu}_k) \right)+ \mathfrak{i}_{\mathcal{C}_\delta}(A u),
\end{align*}
where $\mathfrak{i}$ denotes the convex indicator function and with $h$ as defined in \eqref{eqn:helper} (for $p=2$). We now have an unconstrained optimization problem involving a convex, yet non-differentiable, functional and thus we can employ a primal-dual algorithm in the variant given in \ref{alg:JKOStep}.

\IncMargin{1em}
\begin{algorithm}
\SetKwInOut{Input}{Input}\SetKwInOut{Output}{Output}
\Input{$u^{(0)},\, \phi^{(0)},\, it_{\max},\, \lambda,\, \sigma,\, A,\, b,\, \delta$}
\Output{$u^*$, $\phi^*$}
\BlankLine
Initialize $\bar u^{(0)}=u^{(0)}$ and $l=0$\;
\For{$i = 0$ \KwTo $it_{\max}$}{
   {$\phi^{(i+1)}={\rm prox}_{\sigma \mathfrak{i}^*_\delta}(\phi^{(i)} +\sigma A \bar u^{(i)})$}\;
  {$u^{(i+1)} = {\rm prox}_{\lambda\Phi}( u^{(i)}-\lambda A^* \phi^{(i+1)})$}\;
  {$\bar u^{(i+1)}= 2 u^{(i+1)}-u^{(i)}$}\;
  \If{convergence}{
  {$u^*=u^{(i+1)}$}\;
  {$\phi^*=\phi^{(i+1)}$}\;
  break}{}
}
\caption{Primal-Dual algorithm for the solution of the optimal transport problems. 
.\label{alg:JKOStep}}
\end{algorithm}\DecMargin{1em}

\noindent In all examples below and unless explicitly stated differently, we make the following parameter choices: $N_x = 150$, $N_t = 75$, $T=1$, $L_e=1$ for all edges $e \in \E$.

\subsection{Example 1: Geodesics and branching}
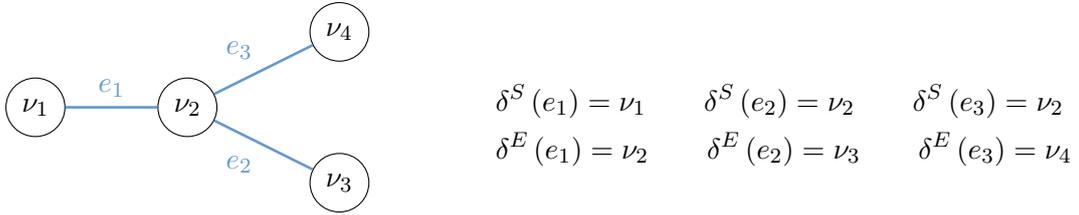
\begin{figure}[ht]
    \centering
    \begin{minipage}{.4\textwidth}
    \begin{tikzpicture}[scale=1, transform shape]
    \pgfdeclarelayer{background}
    \pgfsetlayers{background,main}
    \begin{scope}[every node/.style={draw,circle}]
            \node (v1) at (-2,0) {$\nu_1$};
            \node (v2) at (0,0) {$\nu_2$};
            \node (v3) at (2,-1) {$\nu_3$};
            \node (v4) at (2,1) {$\nu_4$};
    \end{scope}
        \begin{scope}[>={Stealth[bluegray]},
            every edge/.style={draw=bluegray,line width=1pt}]
            \path [-] (v1) edge node[bluegray,above] {$e_1$} (v2);
            \path [-] (v2) edge node[bluegray,below left] {$e_2$} (v3);
            \path [-] (v2) edge node[bluegray,above left] {$e_3$} (v4);
    \end{scope}
    \end{tikzpicture}
    \end{minipage}
    \begin{minipage}{.5\textwidth}
        \begin{align*}
            & \delta^S\left(e_1\right) = \nu_1 \qquad \delta^S\left(e_2\right) = \nu_2 \qquad \delta^S\left(e_3\right) = \nu_2 \\
            & \delta^E\left(e_1\right) = \nu_2 \qquad  \delta^E\left(e_2\right) = \nu_3 \qquad \delta^E\left(e_3\right) = \nu_4
        \end{align*}
    \end{minipage}
    \caption{Sketch of the graph used in the first example for branching geodesics. Here, no in- or outflux via the boundary is assumed (i.e. $\partial^+ \V = \partial^- \V = \emptyset$).}
    \label{fig:graph_erbar_example}
\end{figure}

\begin{figure}[ht]
    \centering
    \includegraphics[width=.32\textwidth]{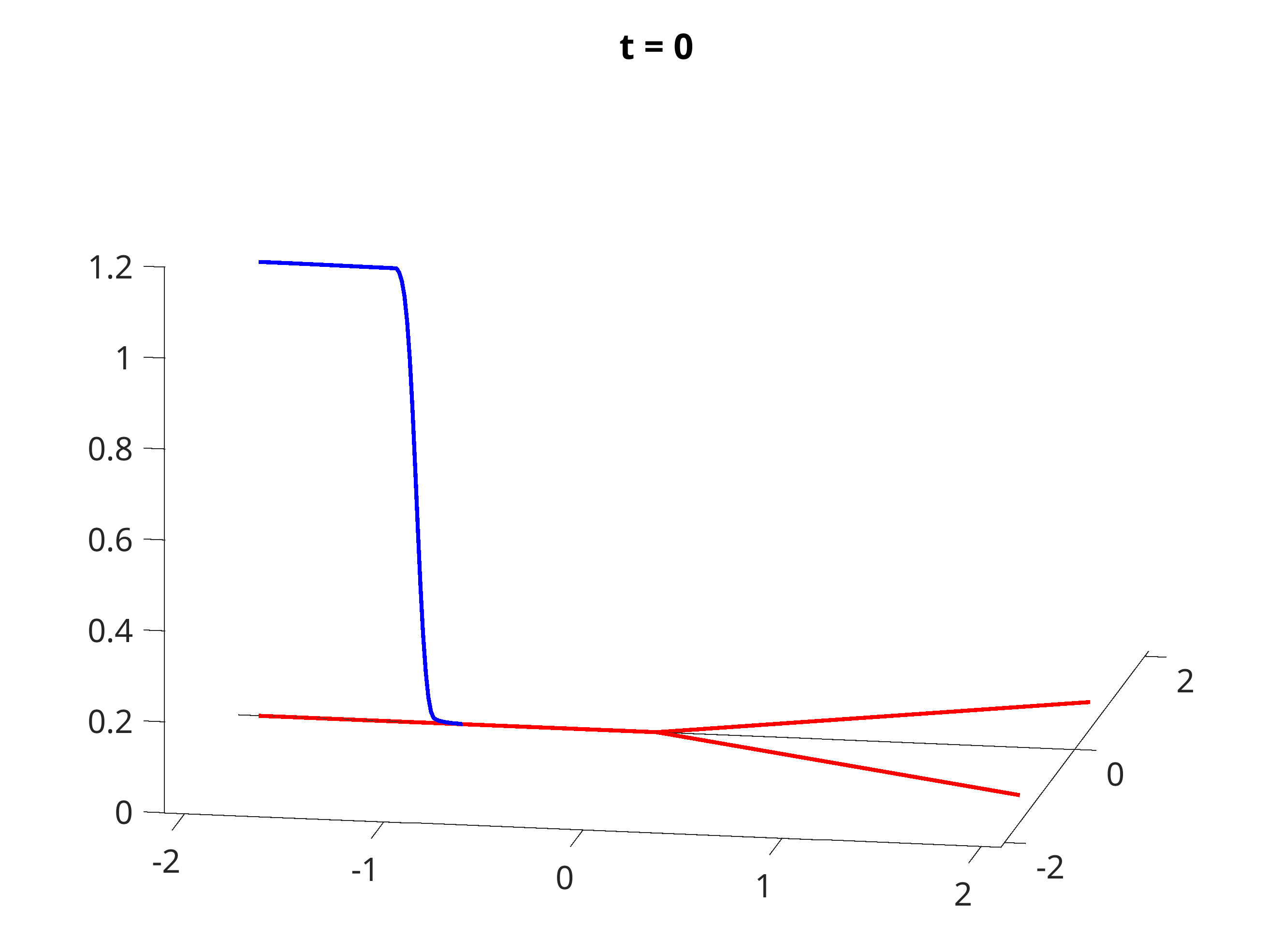}
    \includegraphics[width=.32\textwidth]{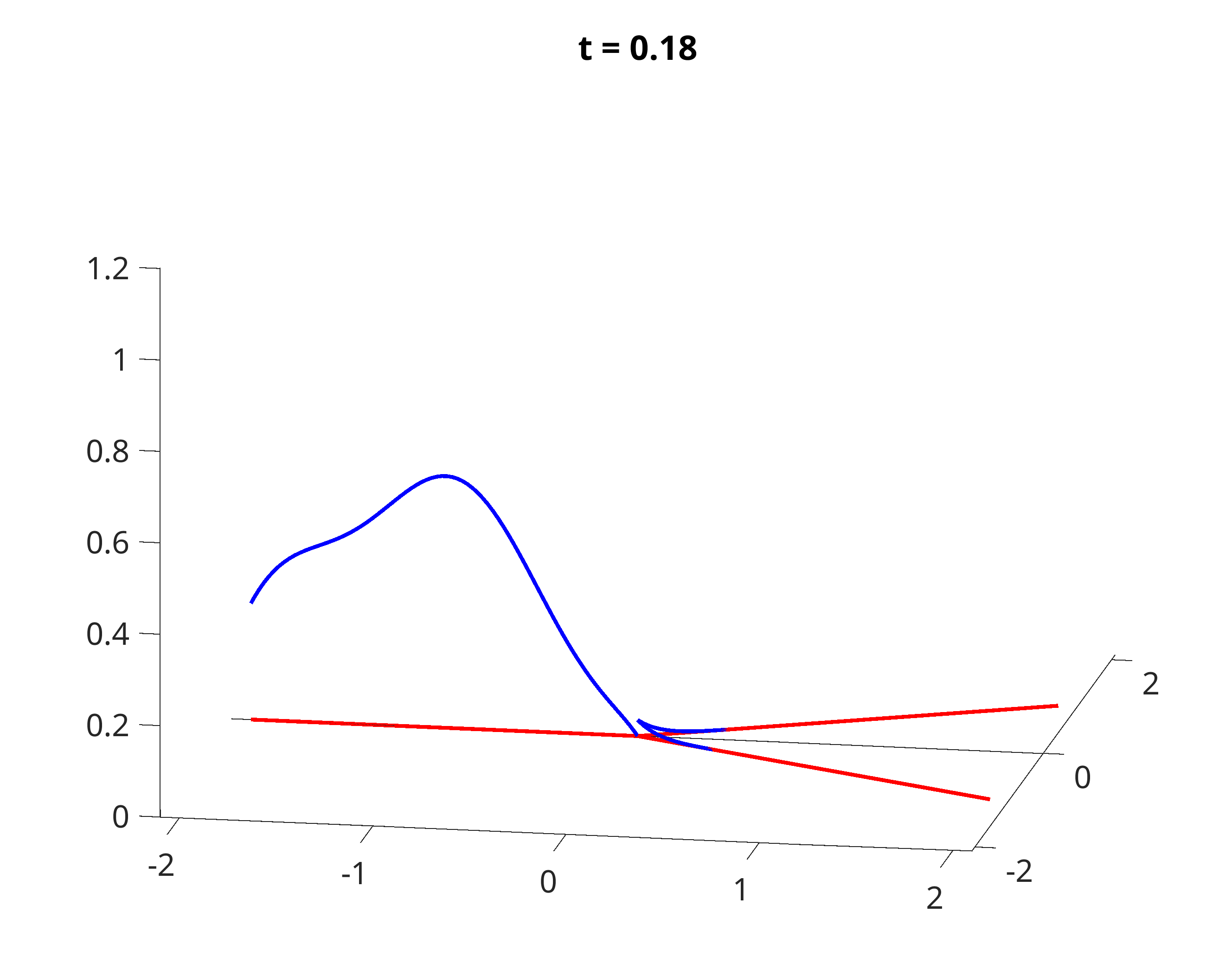}
    \includegraphics[width=.32\textwidth]{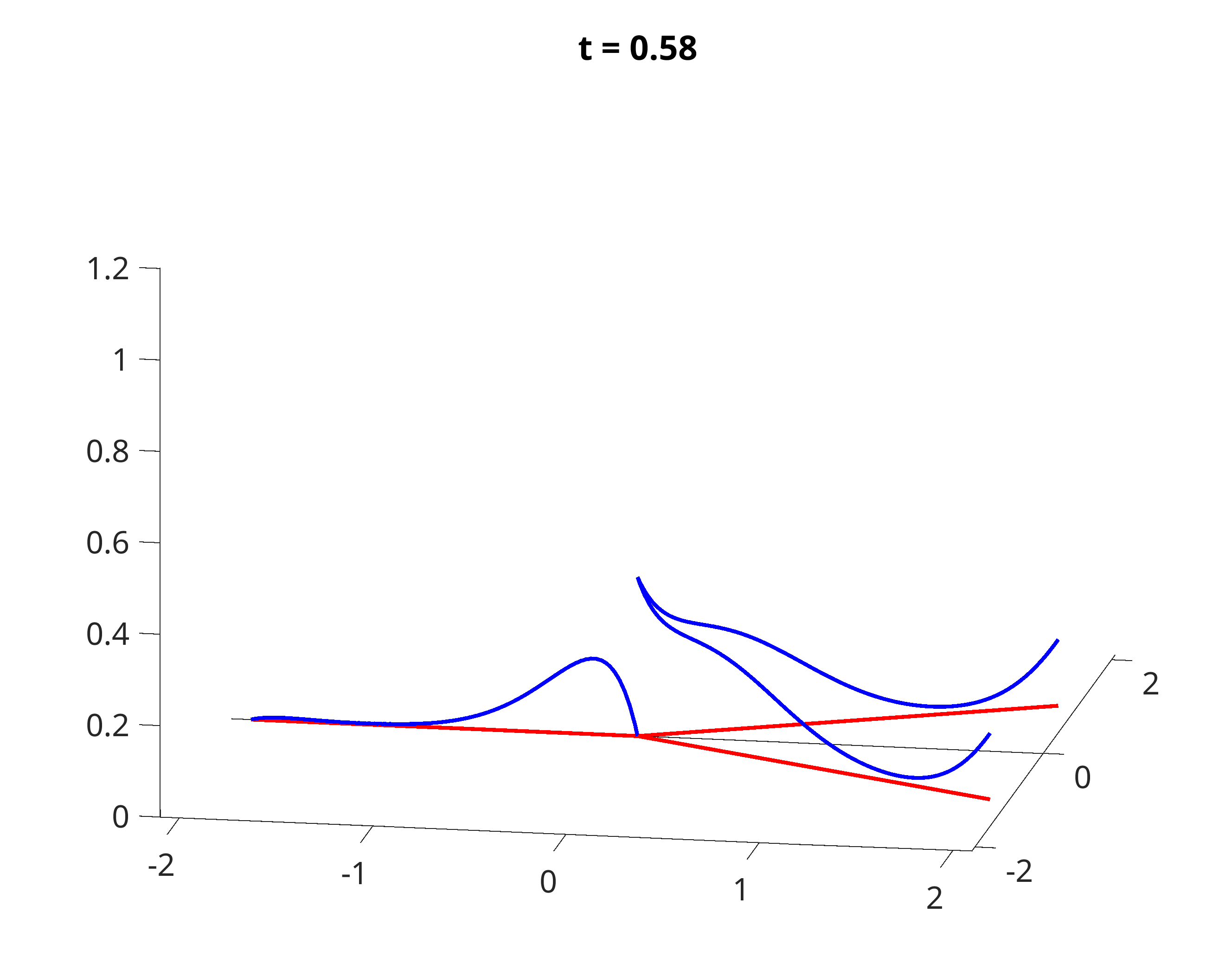}\\
    \includegraphics[width=.32\textwidth]{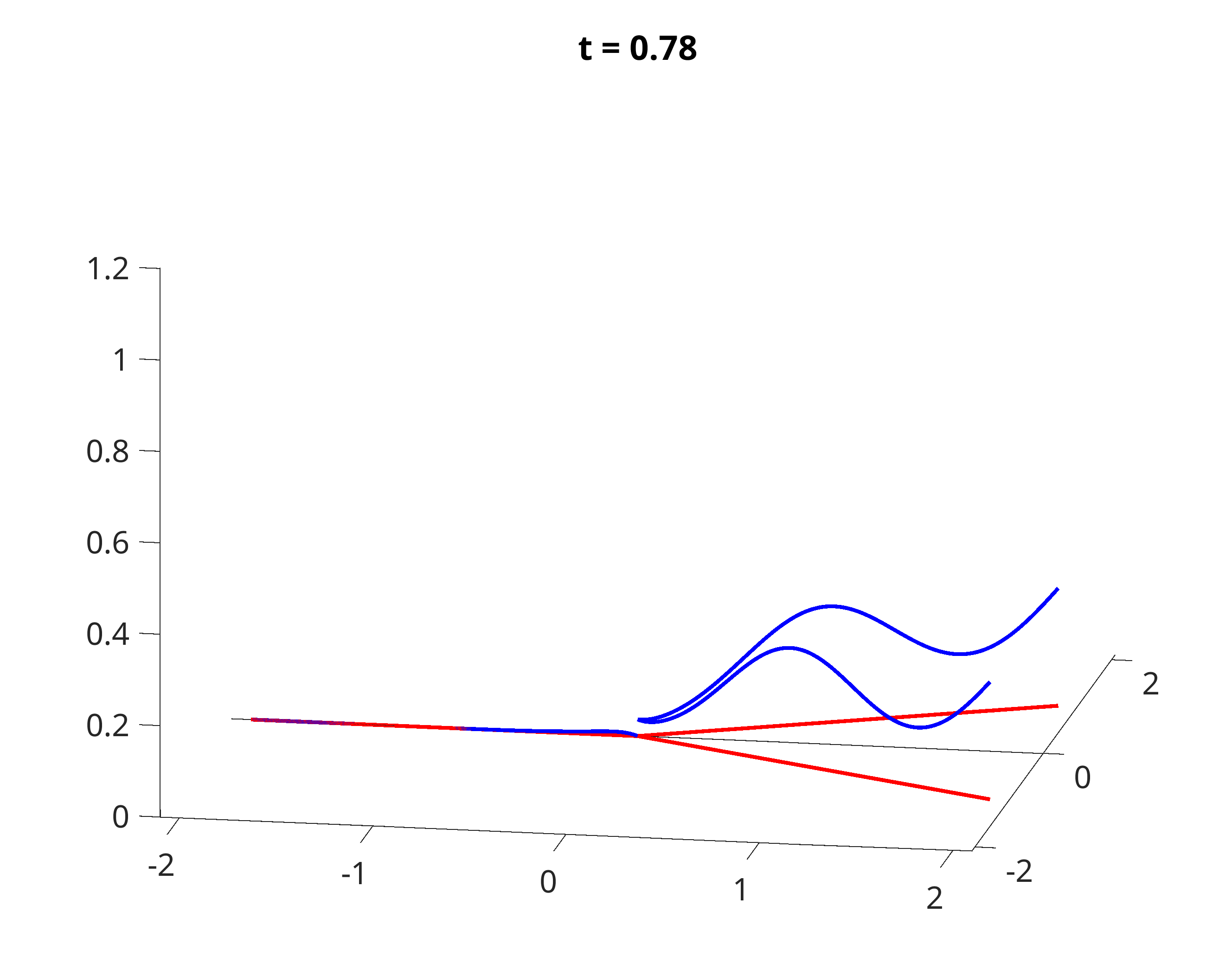}
    \includegraphics[width=.32\textwidth]{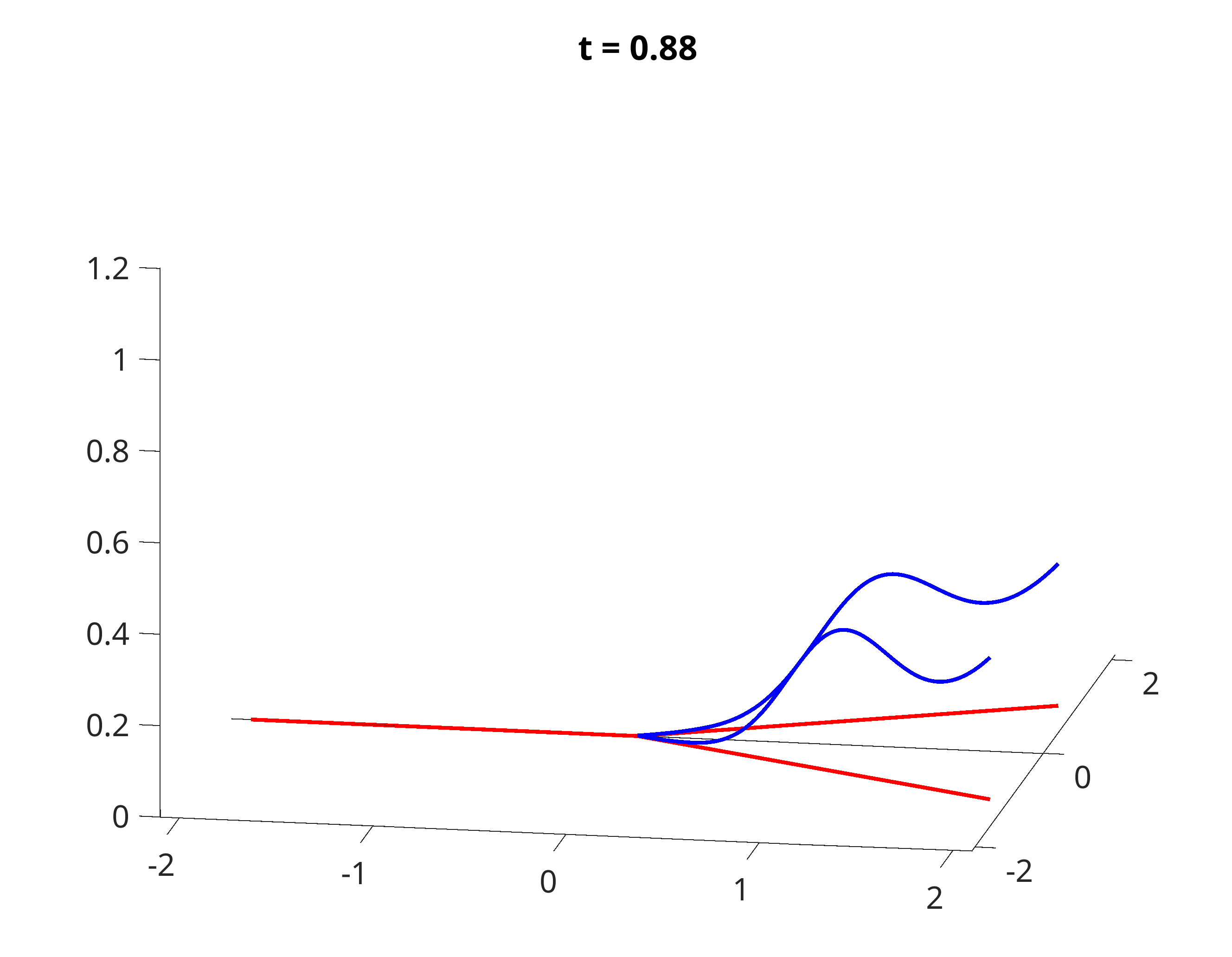}
    \includegraphics[width=.32\textwidth]{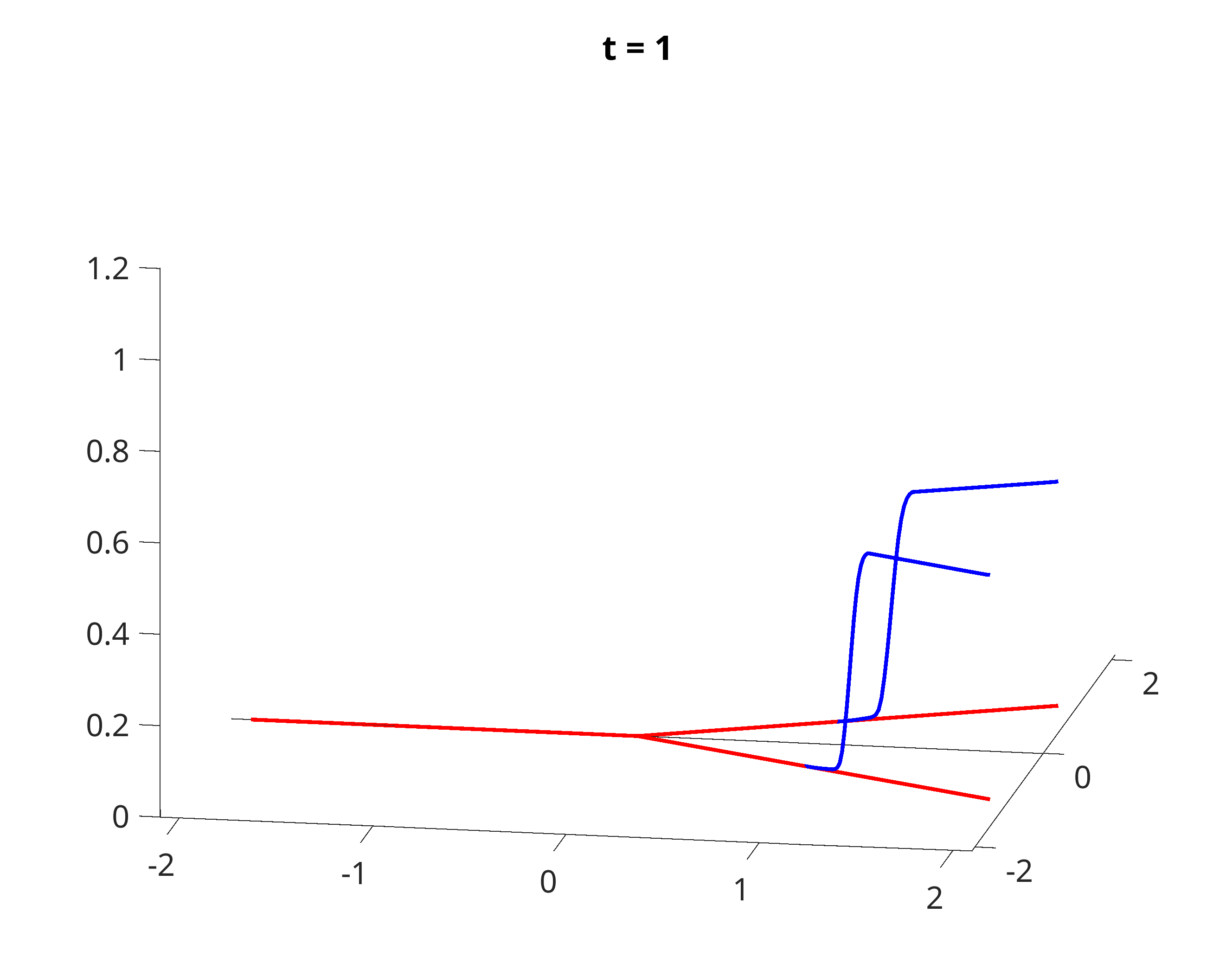}
    \caption{Branching geodesic without vertex dynamic: snapshots of the dynamics of the densities $\rho_e$  at different times.}
    \label{fig:num_ex_erbar_no_vertices}
\end{figure}

\begin{figure}[ht]
    \centering
    \includegraphics[width=.32\textwidth]{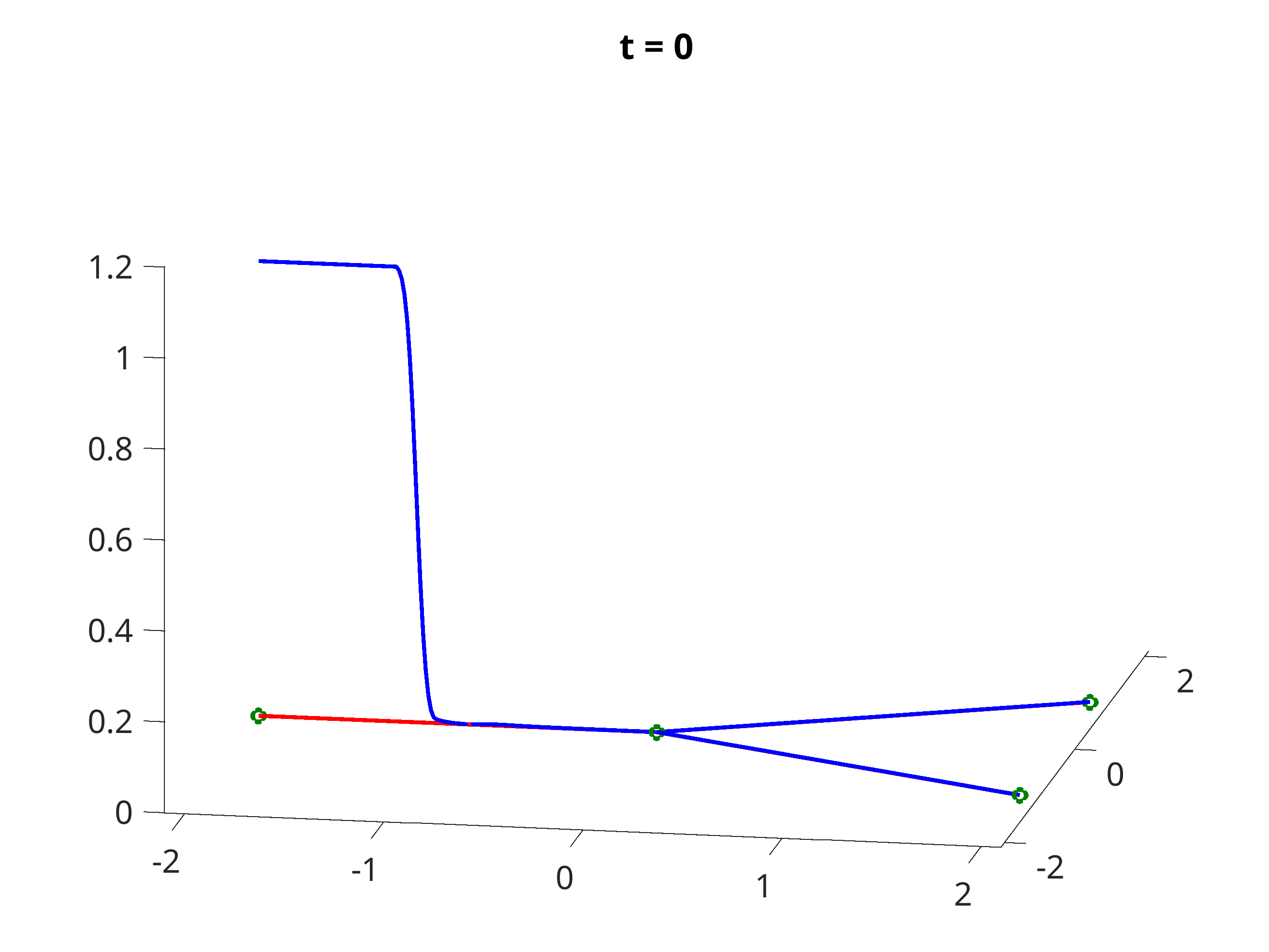}
    \includegraphics[width=.32\textwidth]{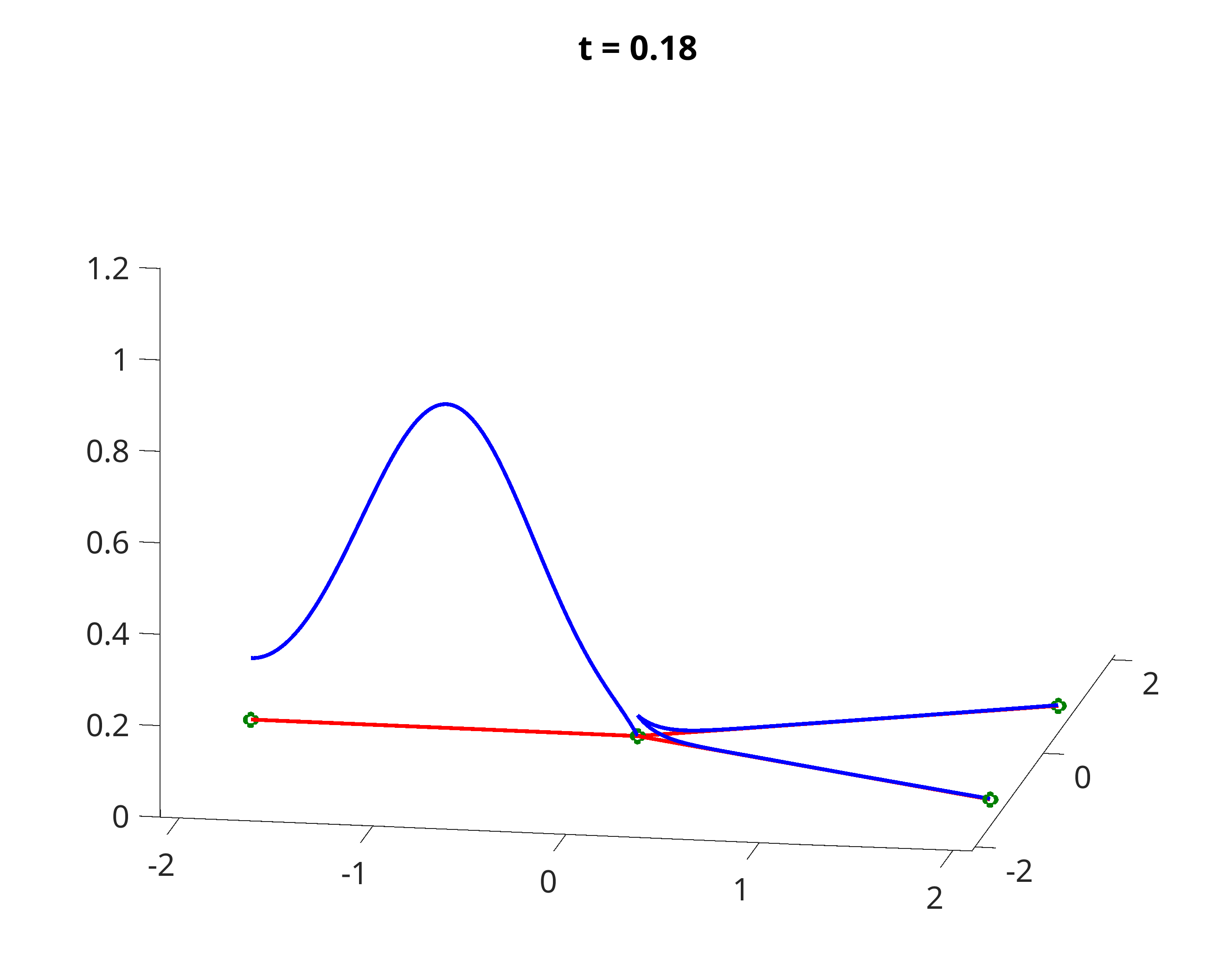}
    \includegraphics[width=.32\textwidth]{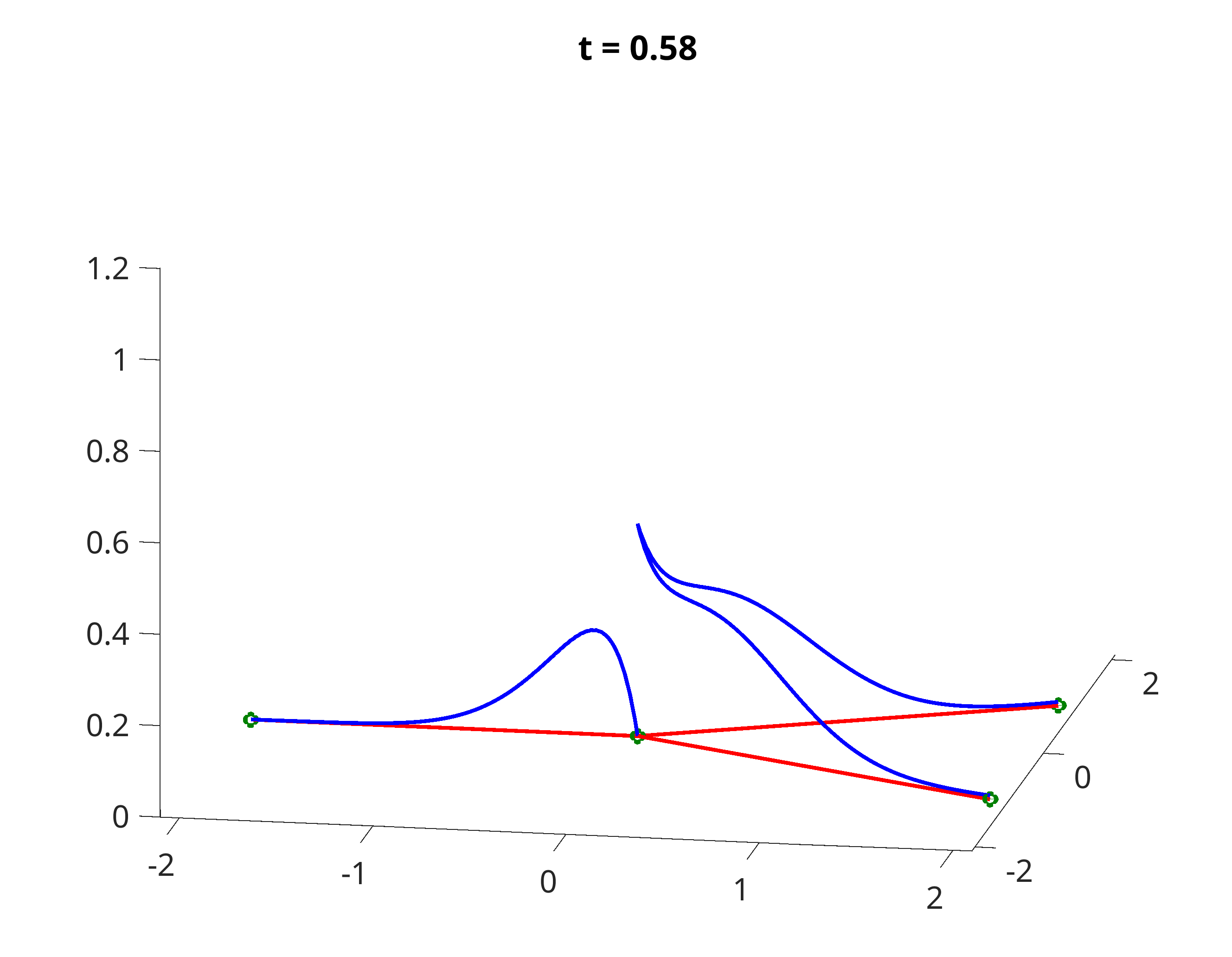}\\
    \includegraphics[width=.32\textwidth]{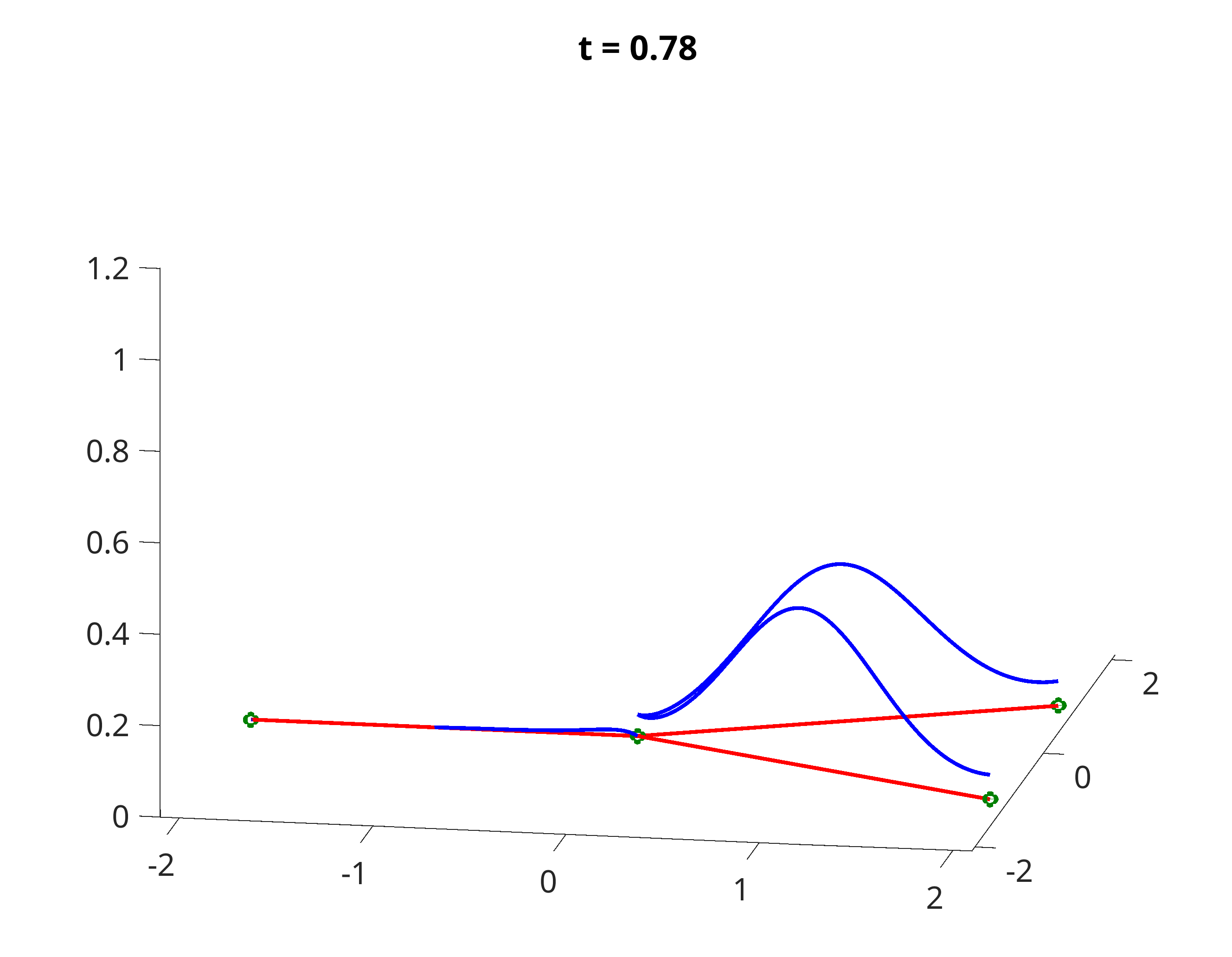}
    \includegraphics[width=.32\textwidth]{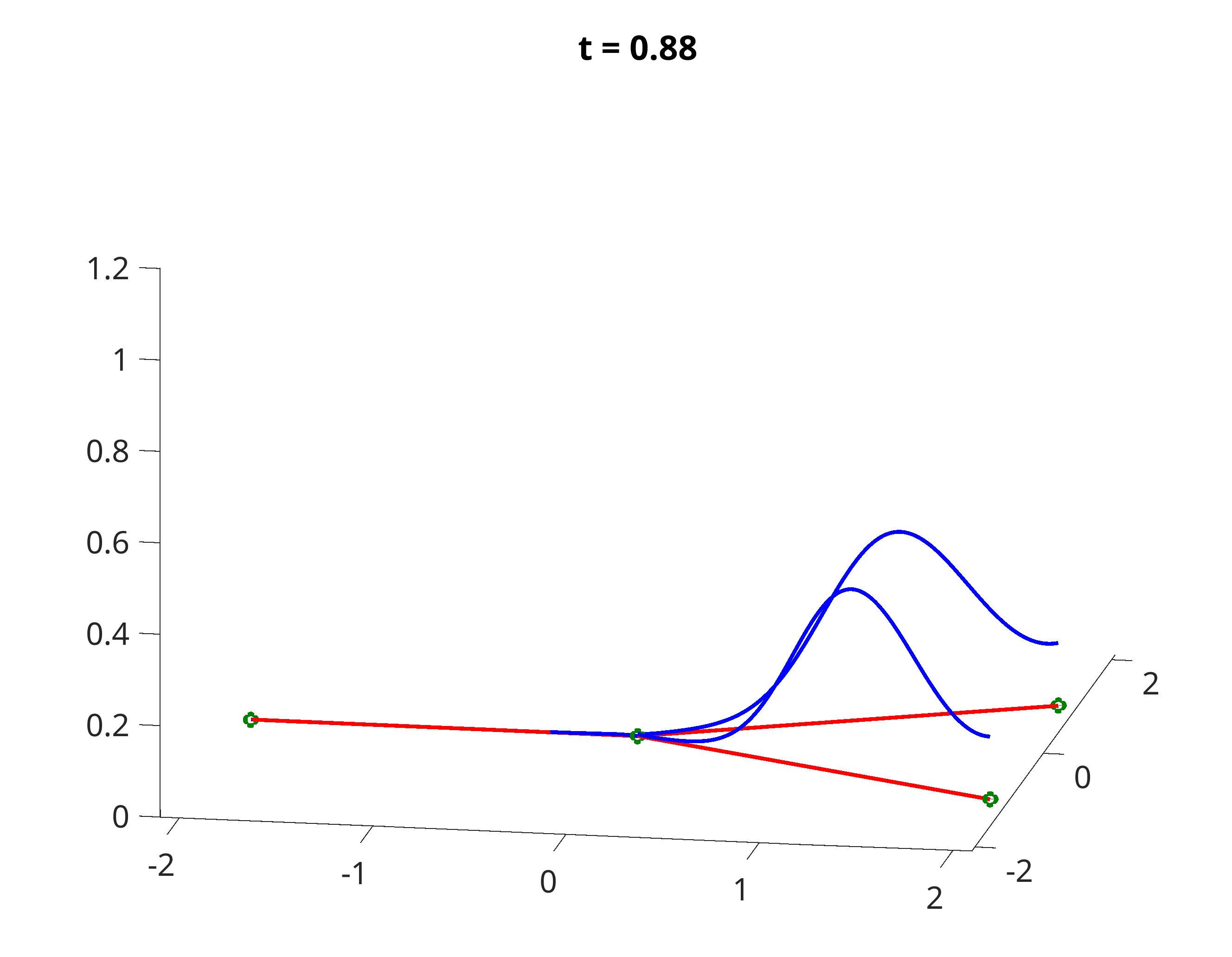}
    \includegraphics[width=.32\textwidth]{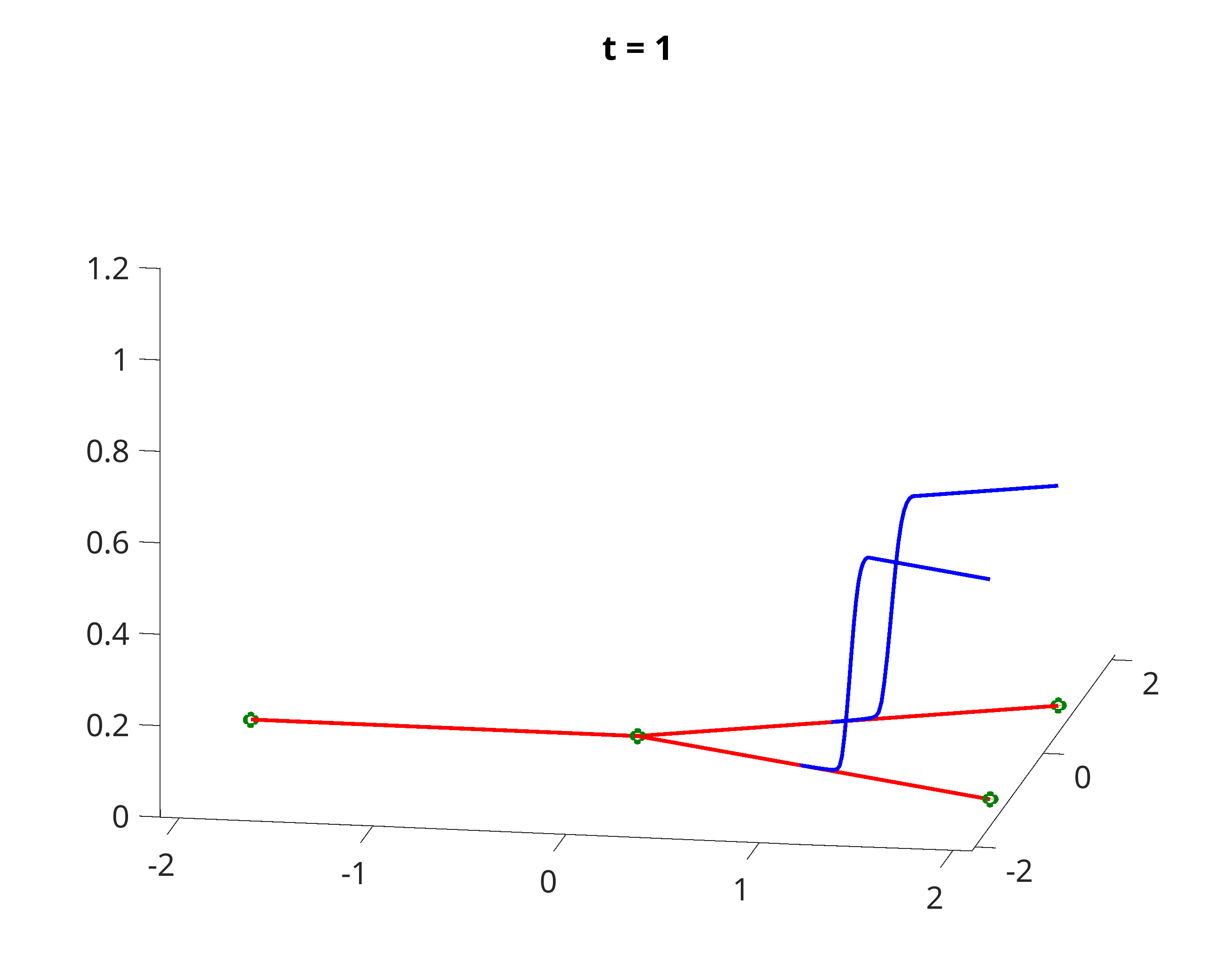}
    \caption{Branching geodesic with vertex dynamic: snapshots of the dynamics of the densities $\rho_e$ and $\gamma_{\nu}$ at different times.}
    \label{fig:num_ex_erbar_vertices}
\end{figure}

One important observation made in \cite{erbar2021gradient} is the fact that the logarithmic entropy on metric graphs is not geodesically convex, due to the fact that geodesics may branch. This is illustrated by an explicit example, which we reproduce here numerically, both for the case with and without vertex dynamics, and for the graph topology depicted in Figure \ref{fig:graph_erbar_example}.
As initial and final data we chose
\begin{equation*}
        \rho_e^0(x):= \begin{cases}
        \indicator_{[0,0.4]}(x), &e =e_1 ,\\ 
        0, & e=e_2, \, e_3,\end{cases}        
                \quad \text{ and }\quad
        \rho_e^1(x) := \begin{cases}0, & e= e_1 ,\\ 
        \frac{1}{2} \indicator_{[1-0.4,1]} (x), & e=e_2,\,e_3, \end{cases}
\end{equation*}
and, in the case when we allow for a vertex dynamic,
\begin{align*}
\gamma_{\nu}^0 = \gamma_{\nu}^1 = 0 \qquad \qquad \forall \nu \in \mathcal{V}.
\end{align*}
The results are depicted in Figures \ref{fig:num_ex_erbar_no_vertices} and \ref{fig:num_ex_erbar_vertices} for the case without and with vertex dynamics, respectively. \jfp{One can observe that with explicit vertex dynamics, the transport over the vertex is slowed down, delaying the whole transport.}


\subsection{Example 2: time-dependent in- and outflow}
\begin{figure}
    \centering
    \includegraphics[width=.32\textwidth]{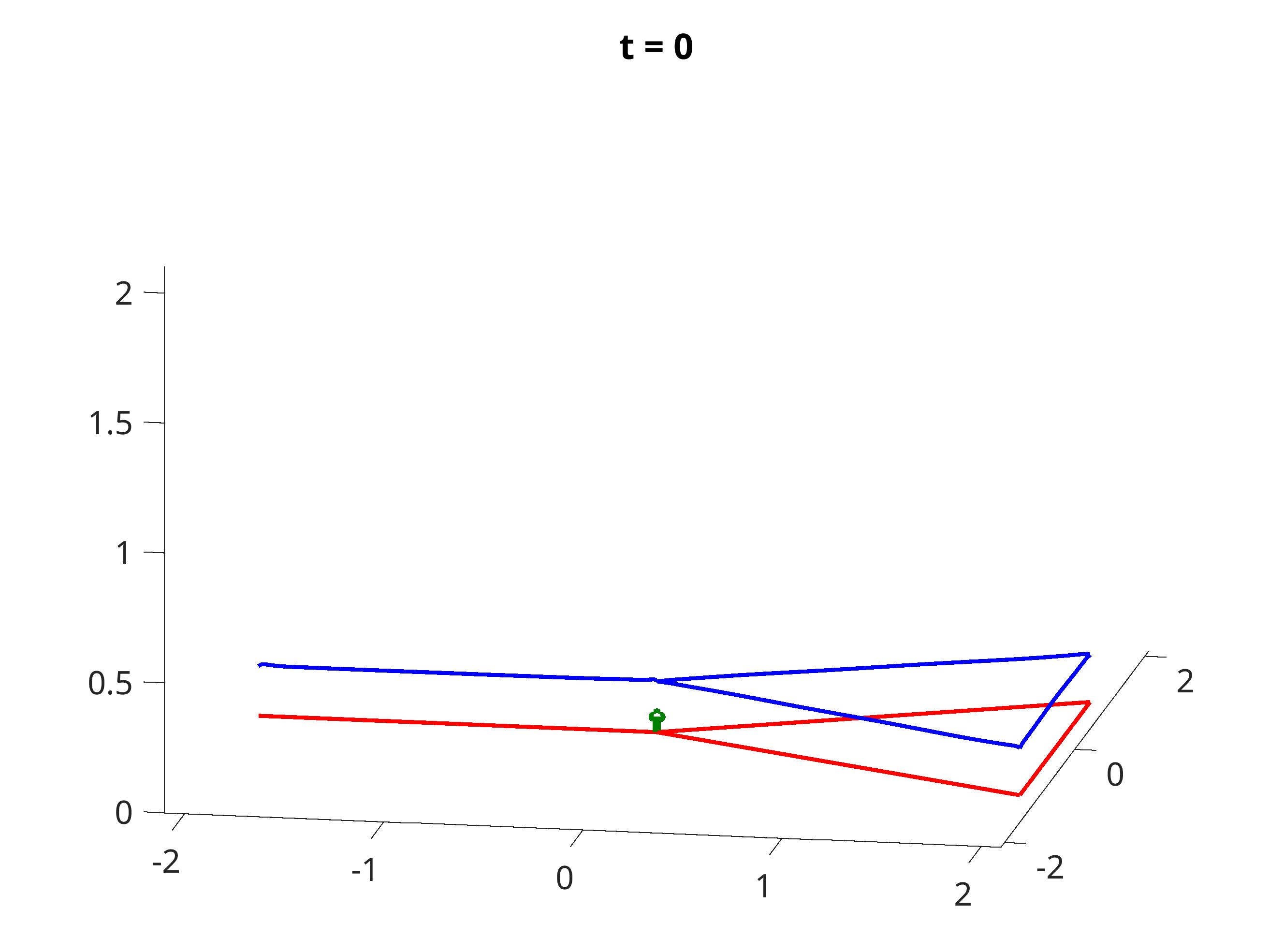}
    \includegraphics[width=.32\textwidth]{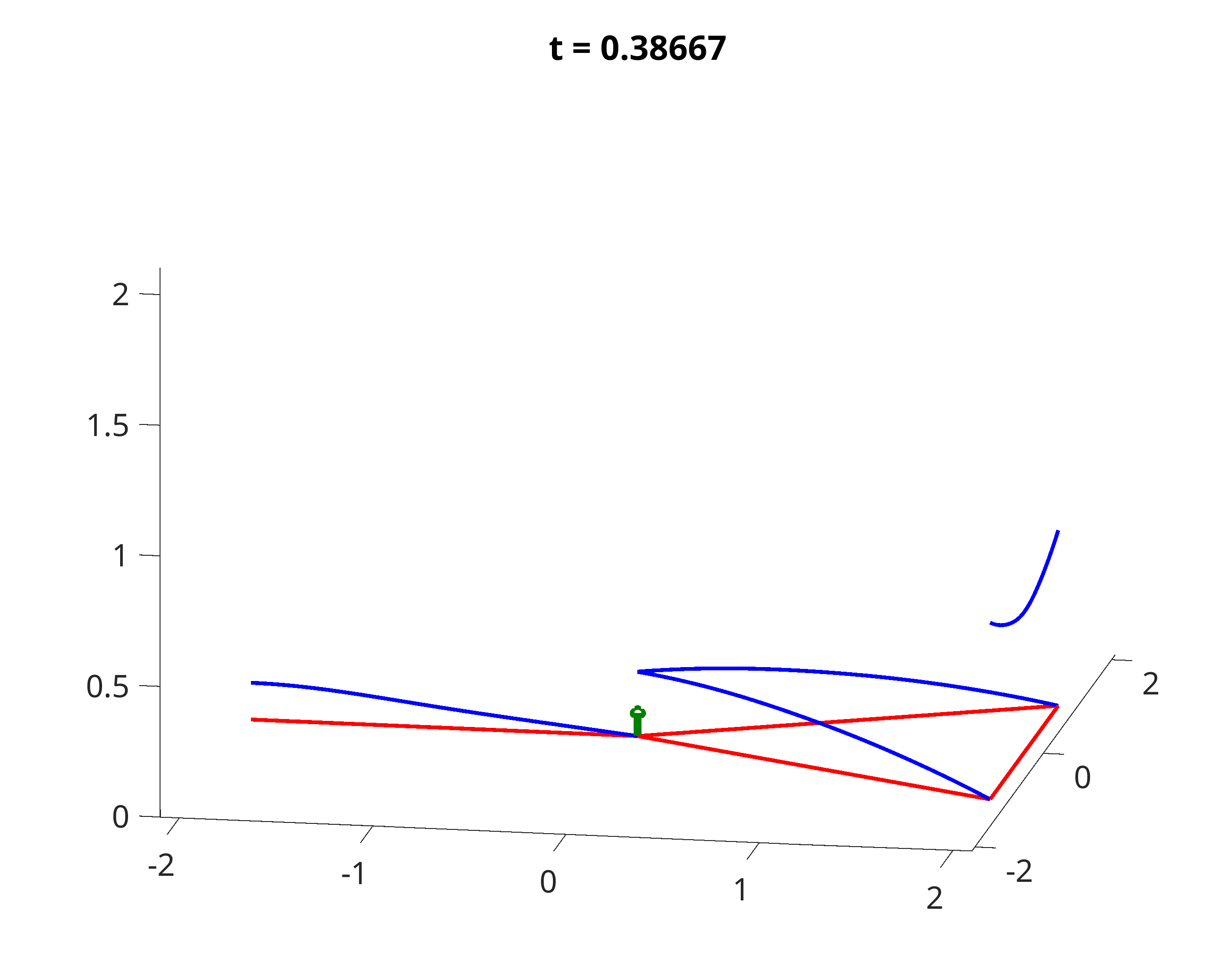}
    \includegraphics[width=.32\textwidth]{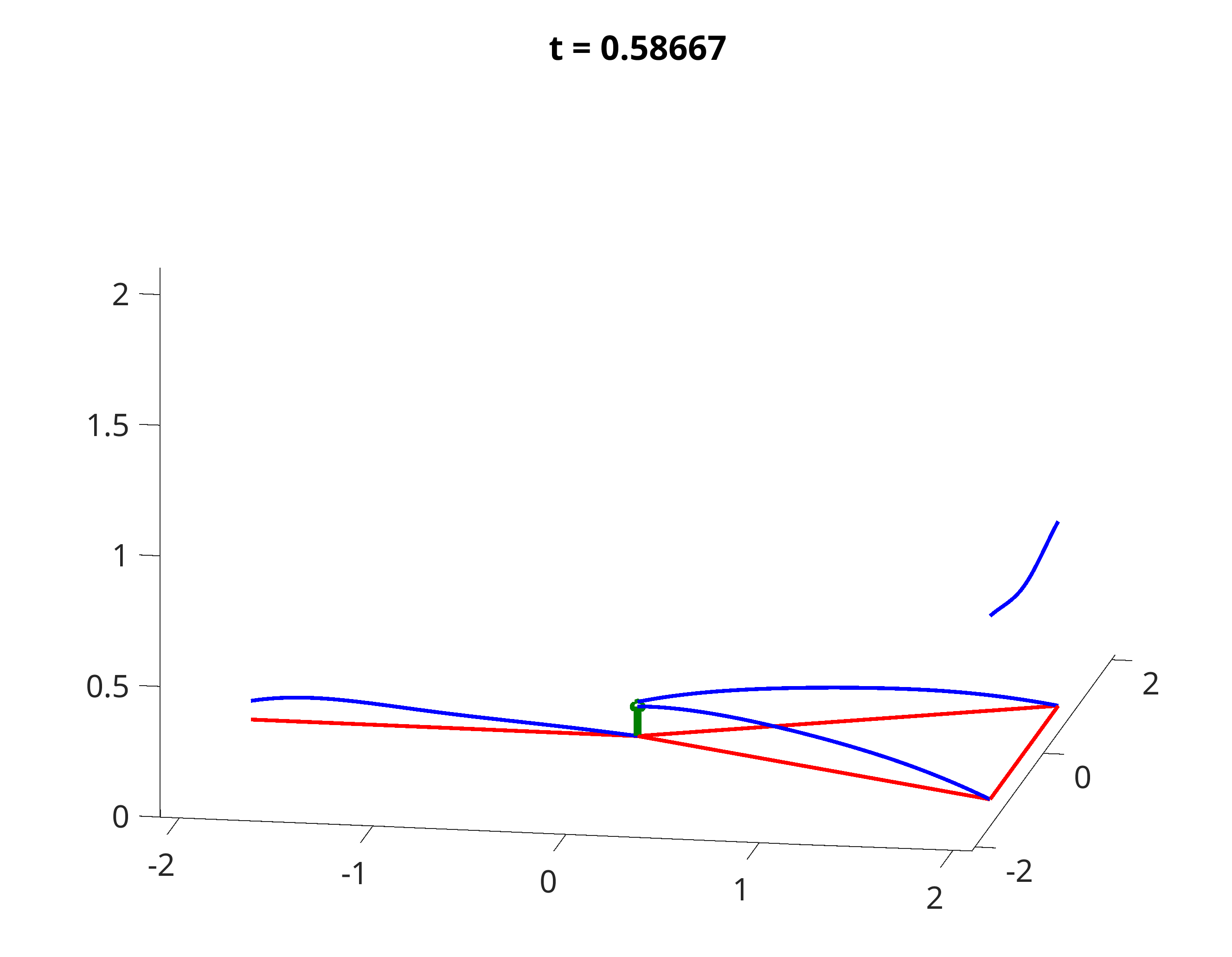}\\
    \includegraphics[width=.32\textwidth]{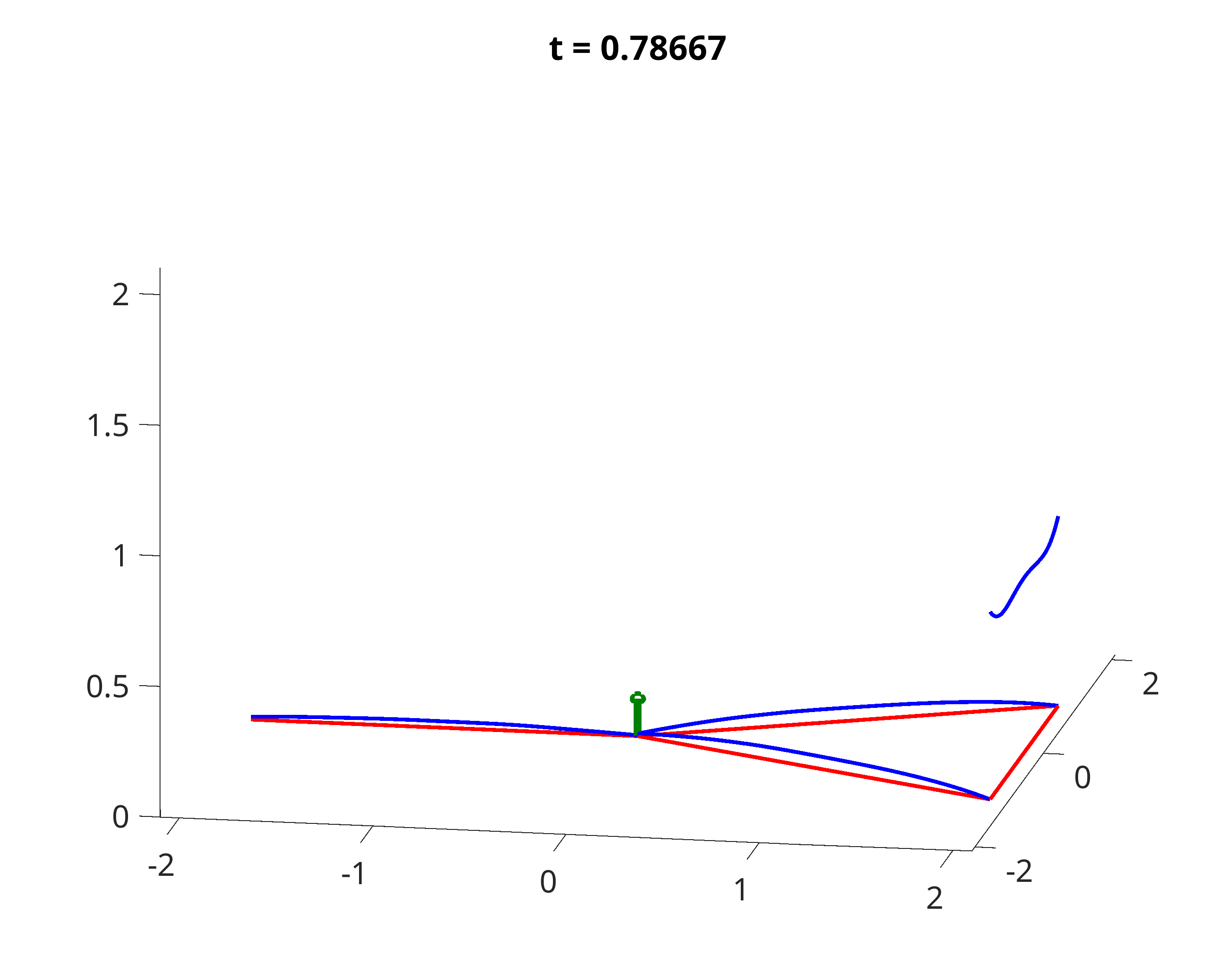}
    \includegraphics[width=.32\textwidth]{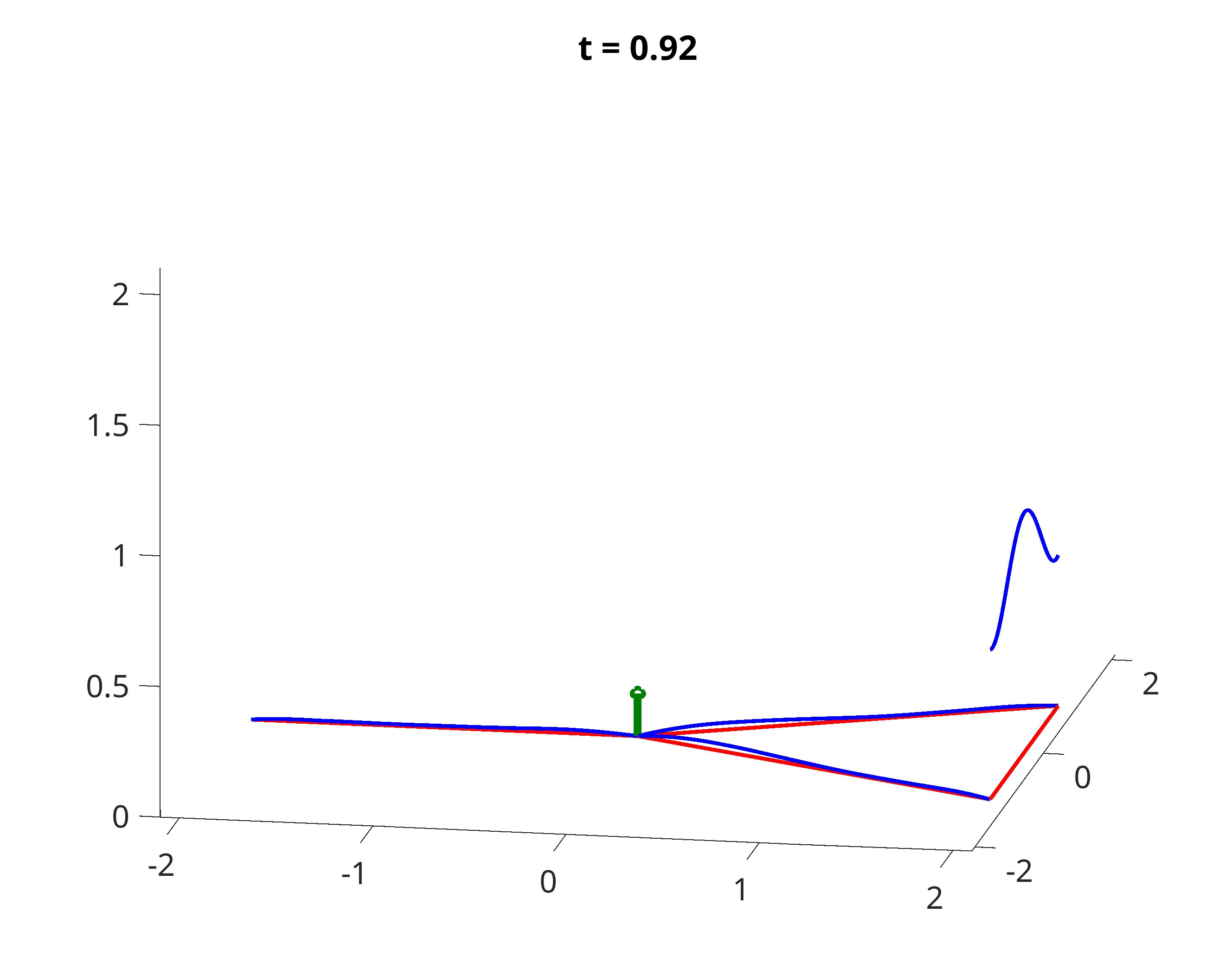}
    \includegraphics[width=.32\textwidth]{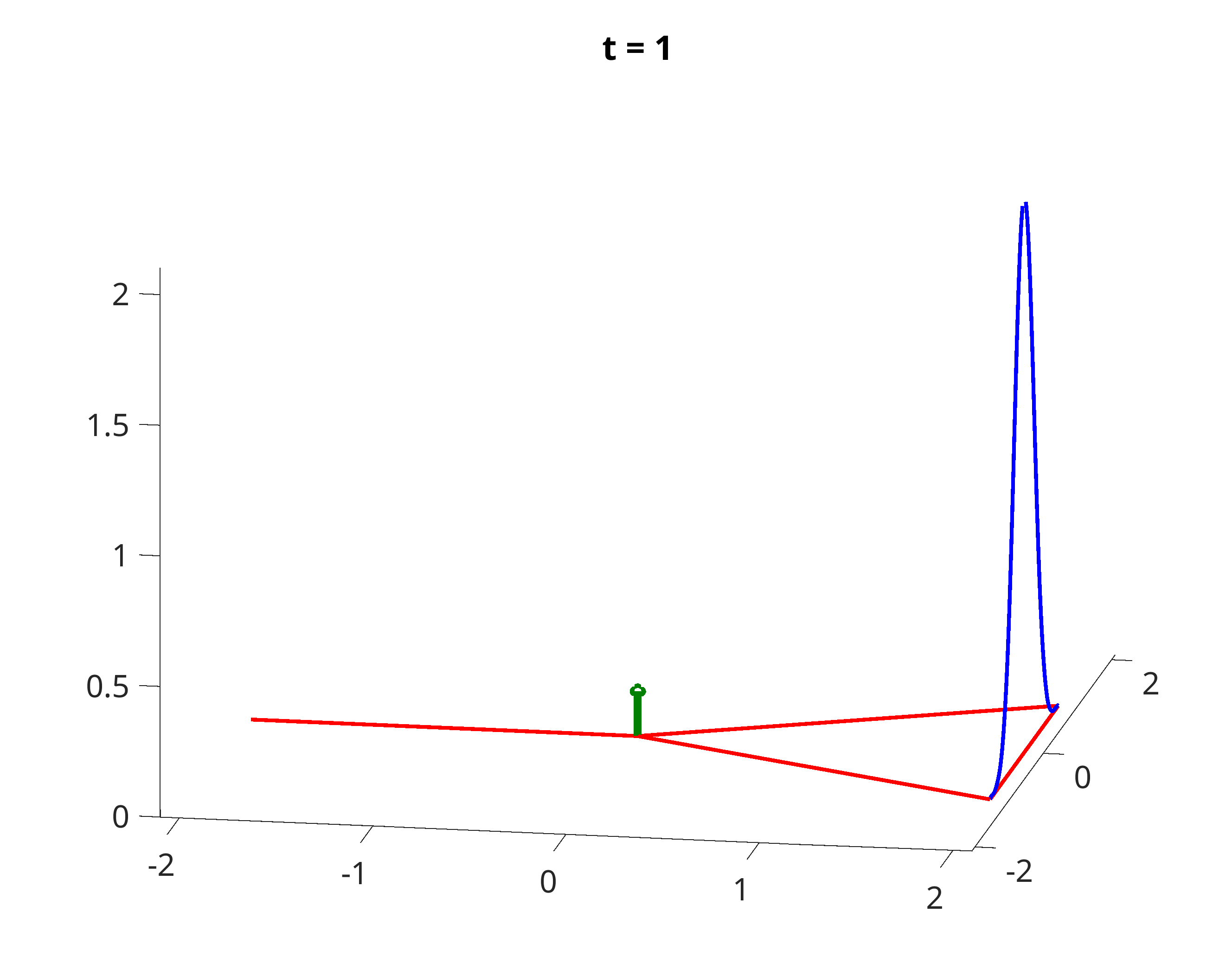}
    \caption{Snapshots of the dynamics of the densities $\rho_e$ and $\gamma_{\nu}$ with symmetric boundary conditions \eqref{eq:bc_in_out_sym} at different times.}
    \label{fig:in_out_sym}
\end{figure}

\begin{figure}
    \centering
    \includegraphics[width=.32\textwidth]{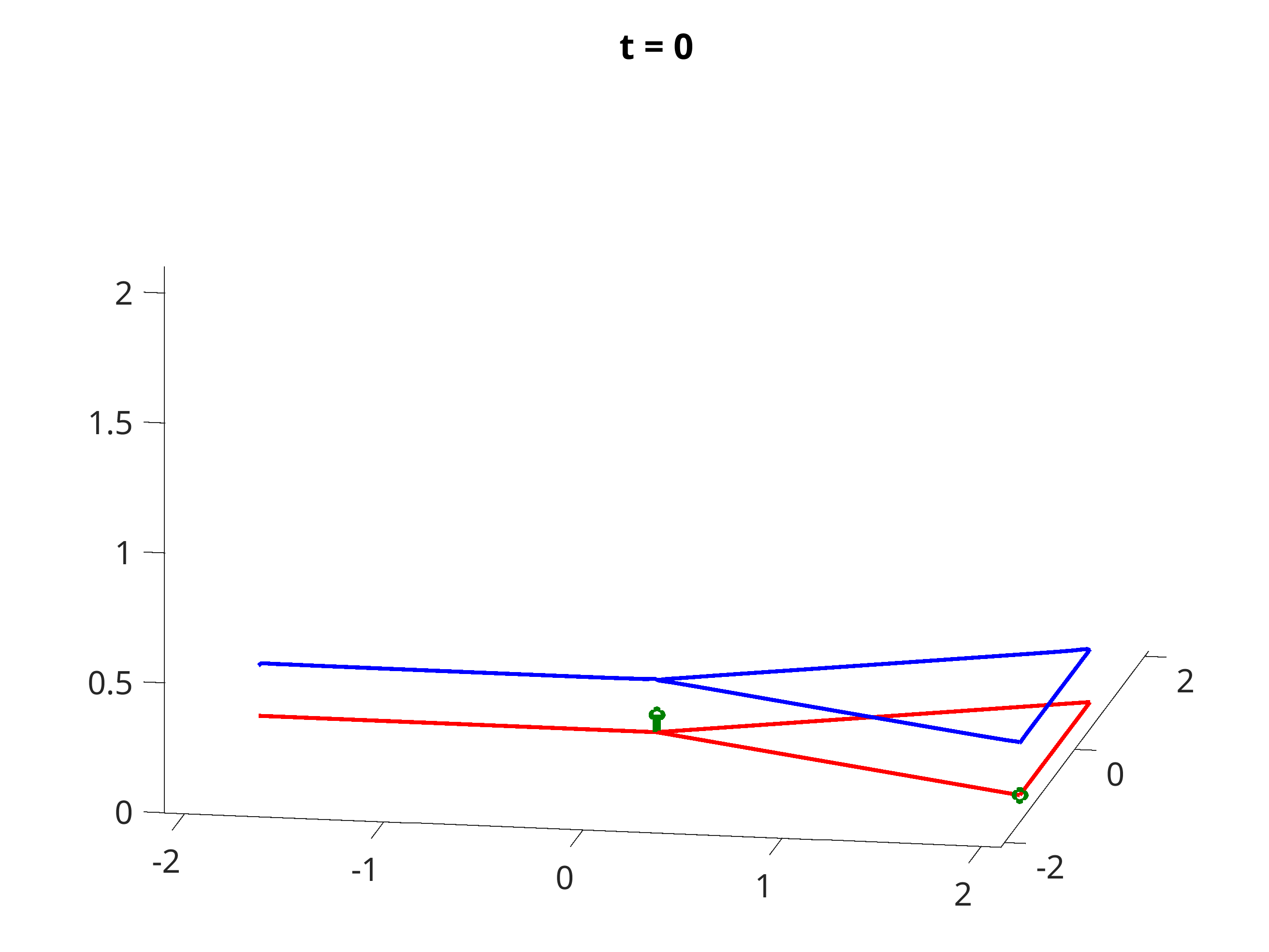}
    \includegraphics[width=.32\textwidth]{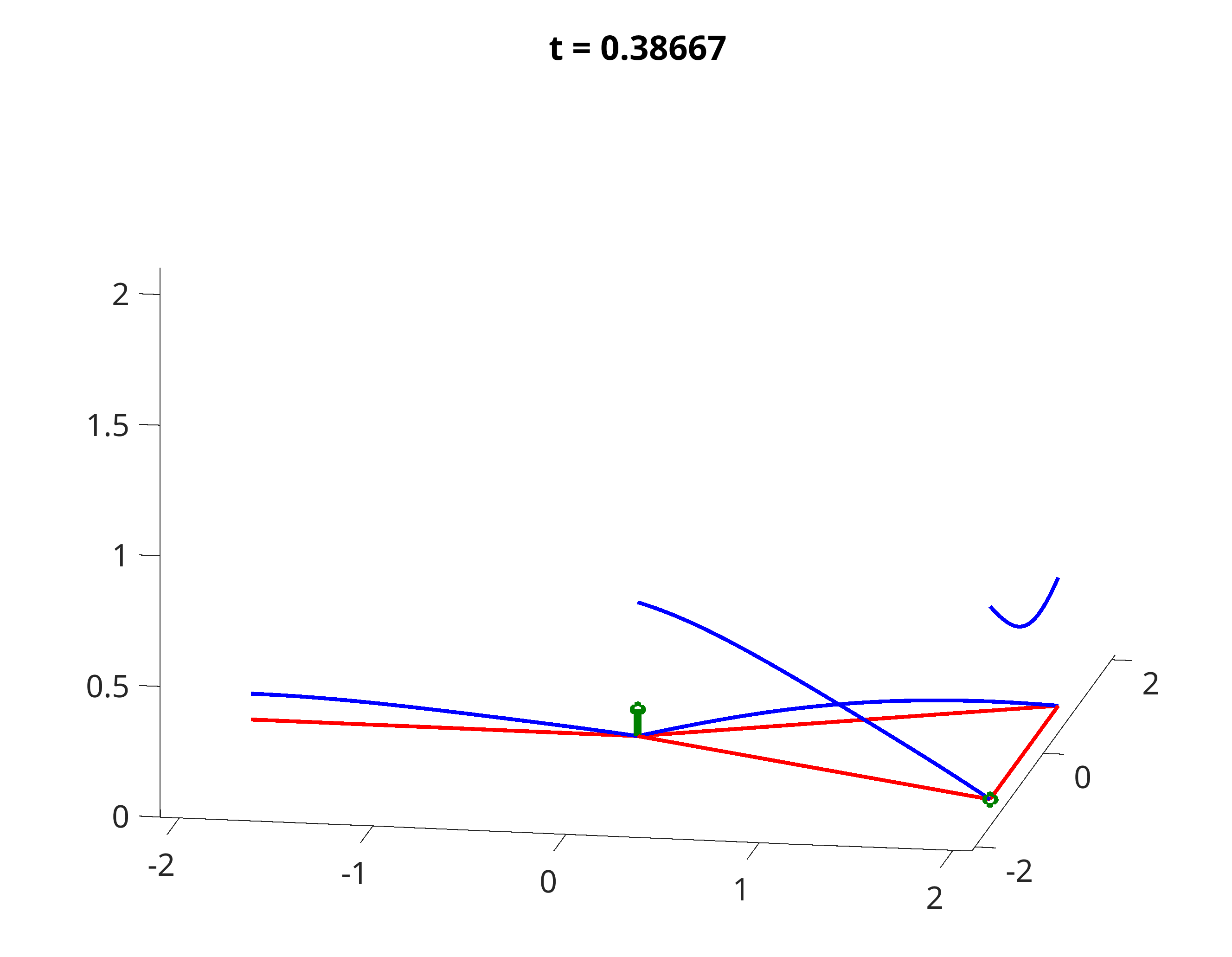}
    \includegraphics[width=.32\textwidth]{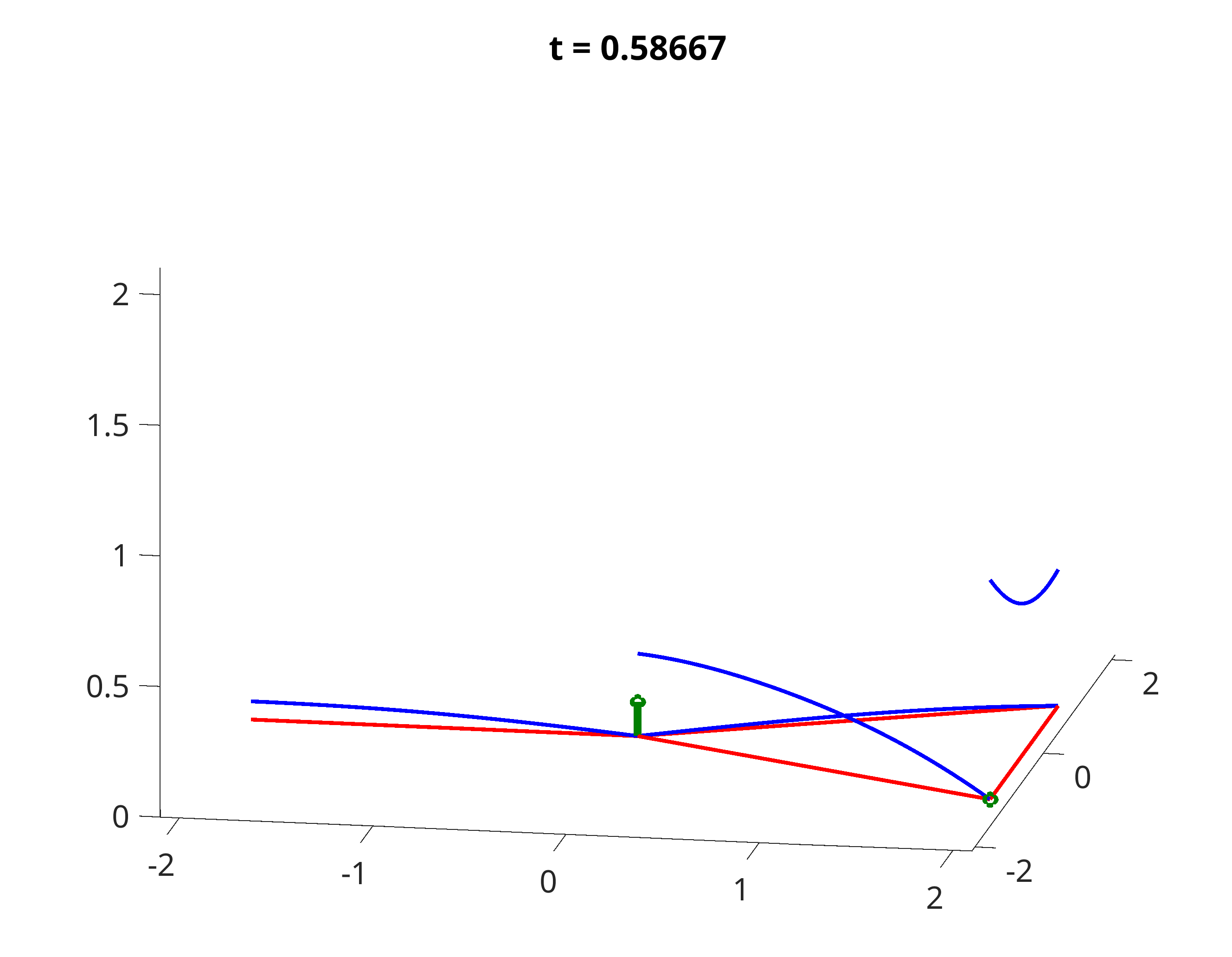}\\
    \includegraphics[width=.32\textwidth]{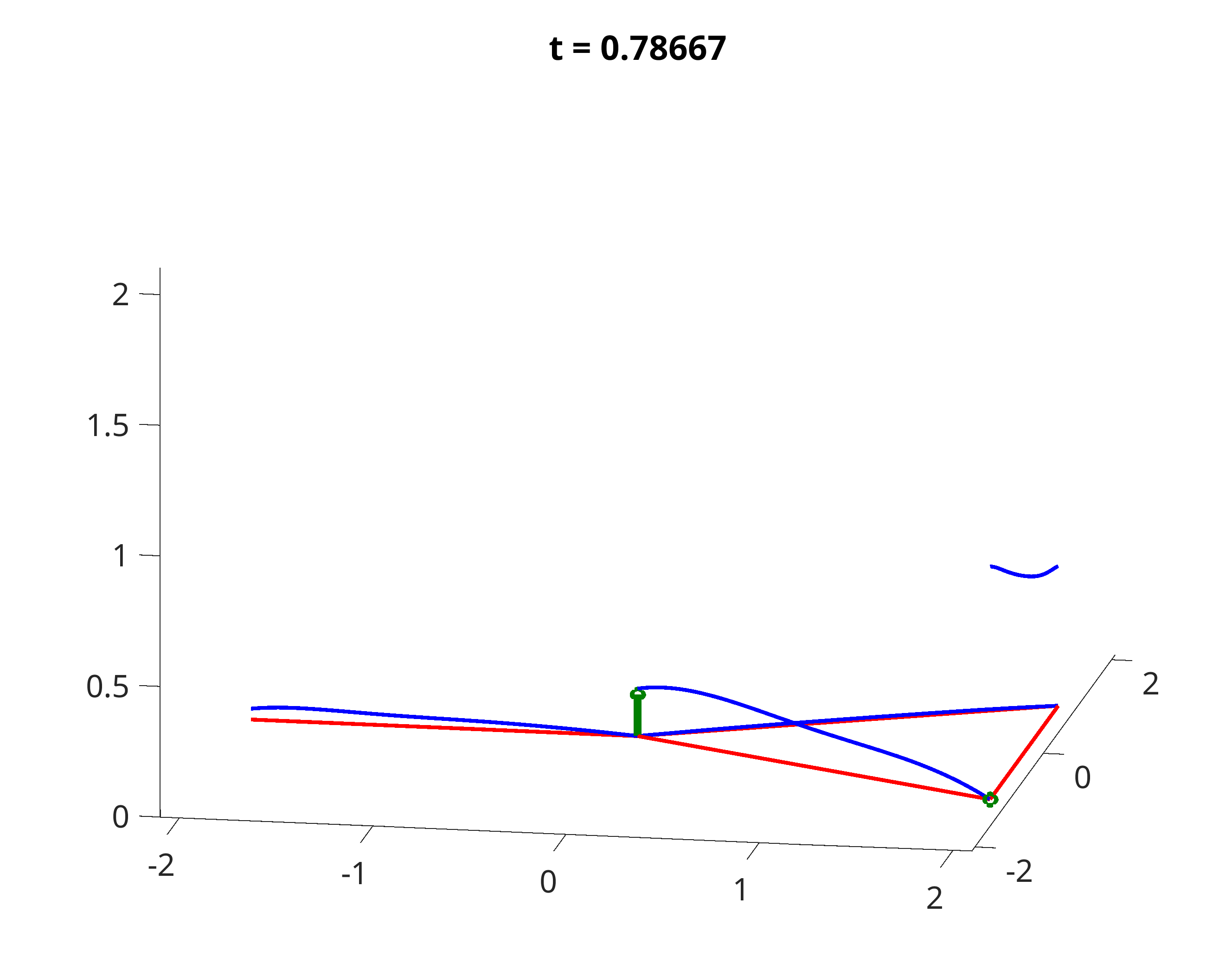}
    \includegraphics[width=.32\textwidth]{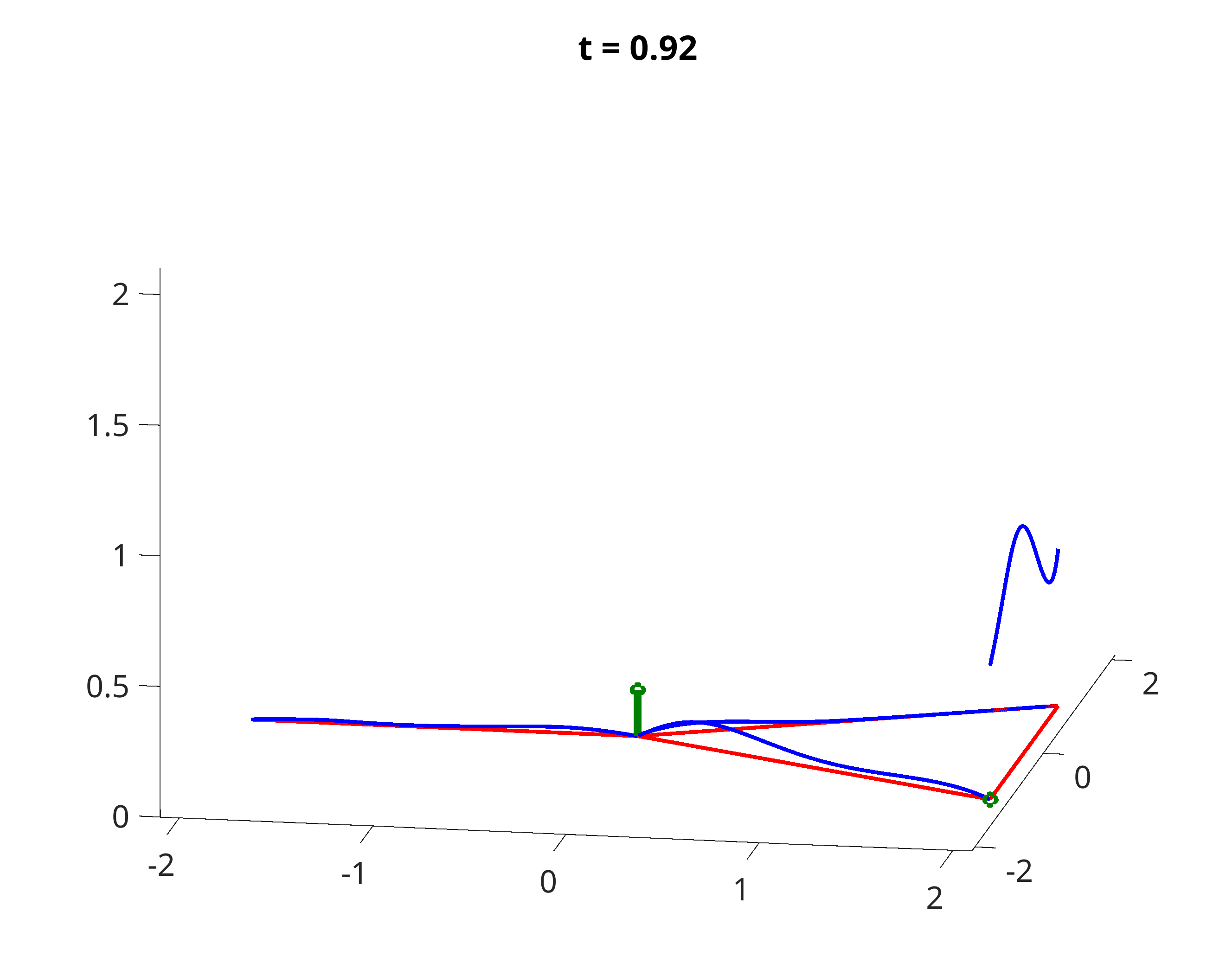}
    \includegraphics[width=.32\textwidth]{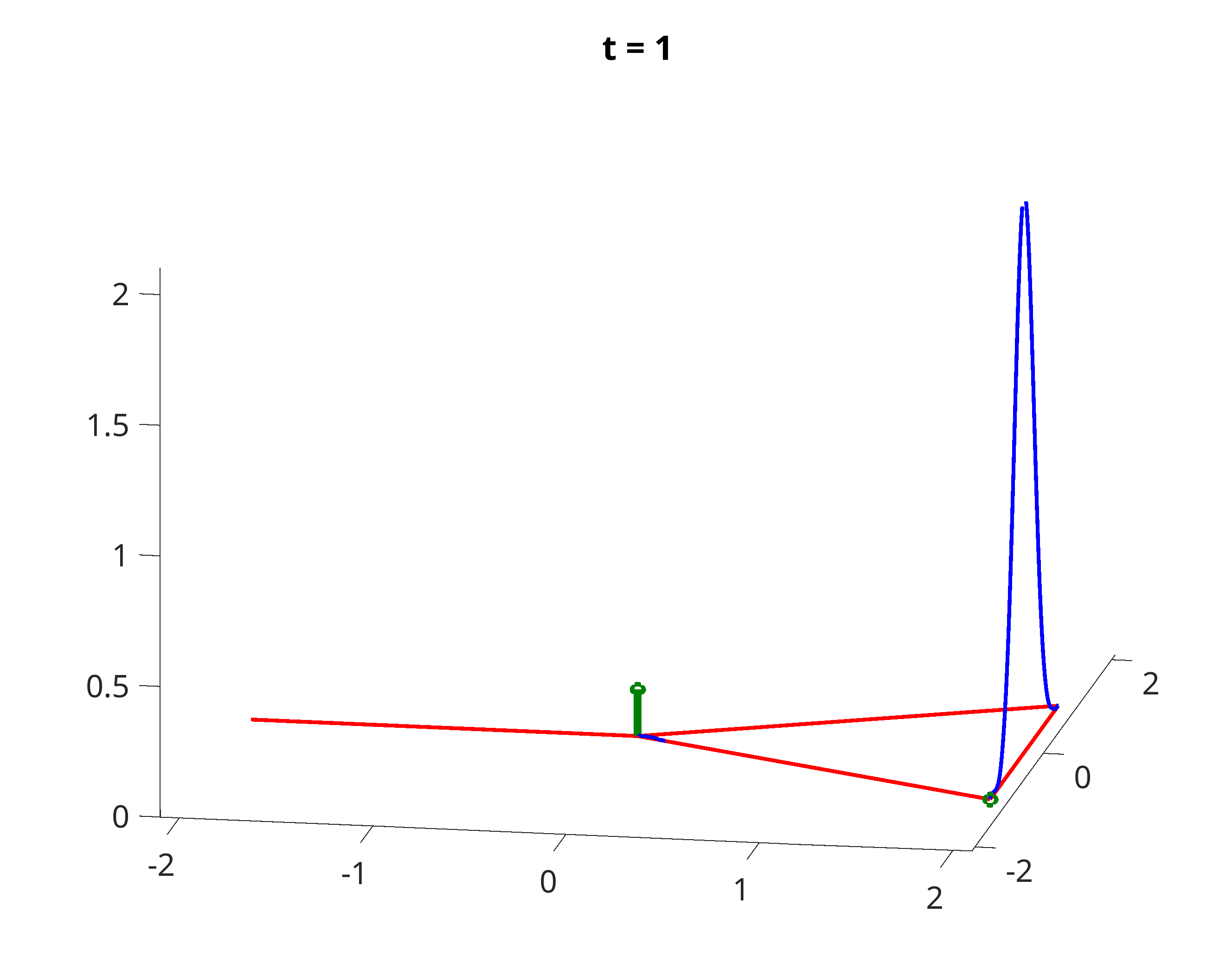}
    \caption{Snapshots of the dynamics of the densities $\rho_e$ and $\gamma_{\nu}$ with asymmetric boundary conditions \eqref{eq:bc_in_out_asym} at different times.}
    \label{fig:in_out_asym}
\end{figure}

In this example, we study the influence of time-dependent boundary conditions on the solutions of the resulting transport problem with vertex dynamics given in definition \ref{def:OT_problems}, combined with (\ref{TD BC 1}) and (\ref{TD BC 2}). Therefore, we fix the graph depicted in figure \ref{fig:graph_in_out} and set
$$
\partial^+ \V = \{\nu_1\},\quad \partial^- \V = \{\nu_3,\,\nu_4\},\,\text{ and }\,\mathring{\V} = \{\nu_2\}.
$$
\begin{figure}[ht]
    \centering
    \begin{minipage}{.4\textwidth}
    \begin{tikzpicture}[scale=1, transform shape]
    \pgfdeclarelayer{background}
    \pgfsetlayers{background,main}
    \begin{scope}[every node/.style={draw,circle}]
            \node (v1) at (-2,0) {$\nu_1$};
            \node (v2) at (0,0) {$\nu_2$};
            \node (v3) at (2,-1) {$\nu_3$};
            \node (v4) at (2,1) {$\nu_4$};
    \end{scope}
        \begin{scope}[>={Stealth[bluegray]},
            every edge/.style={draw=bluegray,line width=1pt}]
            \path [-] (v1) edge node[bluegray,above] {$e_1$} (v2);
            \path [-] (v2) edge node[bluegray,below left] {$e_2$} (v3);
            \path [-] (v2) edge node[bluegray,above left] {$e_3$} (v4);
            \path [-] (v4) edge node[bluegray,right] {$e_4$} (v3);
    \end{scope}
    \end{tikzpicture}
    \end{minipage}
    \begin{minipage}{.5\textwidth}
        \begin{align*}
            & \delta^S\left(e_1\right) = \nu_1 \qquad \qquad \delta^S\left(e_2\right) = \nu_2 \\
            & \delta^E\left(e_1\right) = \nu_2 \qquad \qquad \delta^E\left(e_2\right) = \nu_3 \qquad \\ \vspace{0.5cm}
            & \delta^S\left(e_3\right) = \nu_2 \qquad \qquad \delta^S\left(e_4\right) = \nu_3 \\
            & \delta^E\left(e_3\right) = \nu_4 \qquad \qquad \delta^E\left(e_4\right) = \nu_4
        \end{align*}
    \end{minipage}
    \caption{Sketch of the graph used in the second example. We set $\partial^+ \V = \{\nu_1\}$, $\partial^- \V = \{\nu_3,\,\nu_4\}$ and $\mathring{\V} = \{\nu_2\}$.}
    \label{fig:graph_in_out}
\end{figure}
For initial and final data we chose
\begin{align*}
\rho_{e_1}^0 = \rho_{e_2}^0 =\rho_{e_3}^0 = \rho_{e_4}^0 = 0.225\quad\text{ as well as }\quad\gamma_{\nu_2}^0 = 0.1,\,\gamma^1_{\nu_3} =0,
\end{align*}
and 
\begin{align*}
\rho_{e_1}^1 = \rho_{e_2}^1 =\rho_{e_3}^1 = 0,\, \rho_{e_4}^1(x) = 
c_1e^{-\frac{(x-1/2)^2}{(0.2)^2}}\quad\text{ as well as }\quad\gamma_{\nu_2}^1 = 0.4,\,\gamma^1_{\nu_3} =0,
\end{align*}
where $c_1$ is chosen such that the total mass is normalized to $1$.
For in- and out-flow we consider two cases:
\begin{align}
s^{G,{\nu_1}}(t) = t\quad \text{ and } \quad  d^{G,{\nu_3}}(t) = d^{G,{\nu_4}}(t) = \frac{1}{2} t,\tag{InOutSym}\label{eq:bc_in_out_sym}\\
\intertext{as well as}
s^{G,{\nu_1}}(t) = t,\quad d^{G,{\nu_3}}(t) = 0  \quad \text{ and } \quad d^{G,{\nu_4}}(t) =  t.\tag{InOutAsym}\label{eq:bc_in_out_asym}
\end{align}
Note that  the values for $\gamma_{\nu_3}$ are only needed in the asymmetric case. The respective results, evaluated at different time steps, are depicted in Figures \ref{fig:in_out_sym} (for \eqref{eq:bc_in_out_sym}) and \ref{fig:in_out_asym} (for \eqref{eq:bc_in_out_asym}), respectively. \jfp{The results show that when only one vertex is an outflow vertex, the mass necessary to obtain the final configuration is predominantly transported via the other edge.}

%% file: 8_con.tex
\section{Conclusion}

This paper gives a detailed introduction in the modelling of gas networks as oriented metric graphs, including two different types of coupling conditions at vertices (one which allows for storing gas in vertices, and one which does not) as well as two different kinds of boundary conditions (enabling gas entering and exiting the network as supply and demand). With this setup, we thoroughly investigate mass conservation properties on the network for given initial, final and boundary data, and formulate various transport type problems on metric graphs.\\

\noindent Furthermore, we generalize the dynamic formulation of the $2$-Wasserstein metric in \cite{burger2023dynamic}, to general $p$ and we also include the before mentioned time-dependent and time-independent boundary conditions in the formulations. We also utilize the presented $p$-Wasserstein metrics to derive gradient flows and for the case $p = 3$ we recover the (ISO3) gas model, which is a particularly interesting result. Moreover, we highlight some difficulties in defining appropriate potentials on metric graphs, when going from a single edge to a simple connected graph, which naturally occur in the study of vanishing diffusion limits. \\

\noindent In the last section of the paper we also present some numerical results based on a space-time discretization called primal-dual gradient scheme, which allows us to compute solutions of the presented optimal transport problems for different coupling and boundary conditions. These examples give insights about how gas storage at interior vertices as well as vertices responsible for the supply of the network or for meeting demands, affect the dynamics of the network and thus the optimal transport solution. \\

\noindent Concerning future work, we are aiming at generalizing the results of \cite{burger2023dynamic} and \cite{erbar2021gradient}, which proof the existence of minimizers for the $p$-Wasserstein metrics in the case of $p = 2$, no boundary vertices and no gas storage at interior vertices. Furthermore, an extension to mixture models (for instance a mix of natural gas and hydrogen) constitutes an interesting research opportunity as well.

\section*{Acknowledgements}

JFP thanks the DFG for support via the Research Unit FOR 5387 POPULAR, Project No. 461909888. AF and MB thank the DFG for support via the SFB TRR 154. MB acknowledges support from DESY (Hamburg, Germany),
a member of the Helmholtz Association HGF.

%% file: appendix.tex
\section{Appendix} \labeltext{appendix}{appendix}

\subsection*{Coupled measures}

A detailed introduction of coupled measures on a subset of edges or subset of vertices on the graph is given by the following definition.

\begin{definition}[Coupled measures on graphs] ~ \\
    \noindent On a graph $\mathcal{G} = \left(\mathcal{V}, \mathcal{E}\right)$, we define coupled mass measures, either on all edges or a subset of or all vertices, as
    \begin{equation*}
        \mathcal{M}_+\left(\mathcal{E}\right) := \mathcal{M}_+\left(e_1\right) \times \mathcal{M}_+\left(e_2\right) \times \dots \times \mathcal{M}_+\left(e_m\right)
    \end{equation*}
    and
    \begin{align*}
        \mathcal{M}_+\left(\partial^+ \mathcal{V}\right) := \, & \mathcal{M}_+\left(\nu_1\right) \times \mathcal{M}_+\left(\nu_2\right) \times \dots \times \mathcal{M}_+\left(\nu_o\right), \\
        \mathcal{M}_+\left(\partial^- \mathcal{V}\right) := \, & \mathcal{M}_+\left(\nu_{o + 1}\right) \times \mathcal{M}_+\left(\nu_{o + 2}\right) \times \dots \times \mathcal{M}_+\left(\nu_{\theta}\right), \\
        \mathcal{M}_+\left(\mathring{\mathcal{V}}\right) := \, & \mathcal{M}_+\left(\nu_{\theta + 1}\right) \times \mathcal{M}_+\left(\nu_{\theta + 2}\right) \times \dots \times \mathcal{M}_+\left(\nu_n\right).
    \end{align*}

    \noindent Moreover, for a fixed time point $t \in \left[0, T\right]$, by $\mathcal{M}_+^t\left(\mathcal{E}\right)$ we define
    \begin{equation*}
        \mathcal{M}_+^t\left(\mathcal{E}\right) := \mathcal{M}_+^t\left(e_1\right) \times \mathcal{M}_+^t\left(e_2\right) \times \dots \times \mathcal{M}_+^t\left(e_m\right)
    \end{equation*}
    and 
    \begin{align*}
        \mathcal{M}_+^t\left(\partial^+ \mathcal{V}\right) := & \, \mathcal{M}_+^t\left(\nu_1\right) \times \mathcal{M}_+^t\left(\nu_2\right) \times \dots \times \mathcal{M}_+^t\left(\nu_o\right),\\
        \mathcal{M}_+^t\left(\partial^- \mathcal{V}\right) := & \, \mathcal{M}_+^t\left(\nu_{o + 1}\right) \times \mathcal{M}_+^t\left(\nu_{o + 2}\right) \times \dots \times \mathcal{M}_+^t\left(\nu_{\theta}\right),\\
        \mathcal{M}_+^t\left(\mathring{\mathcal{V}}\right) := & \, \mathcal{M}_+^t\left(\nu_{\theta + 1}\right) \times \mathcal{M}_+^t\left(\nu_{\theta + 2}\right) \times \dots \times \mathcal{M}_+^t\left(\nu_n\right).
    \end{align*}

    \noindent In a similar manner, for measures without non-negativity constraints, we define
    \begin{equation*}
        \mathcal{M}\left(\mathcal{E}\right) = \mathcal{M}\left(e_1\right) \times \mathcal{M}\left(e_2\right) \times \dots \times \mathcal{M}\left(e_m\right)
    \end{equation*}
    and
    \begin{equation*}
        \mathcal{M}\left(\mathring{\mathcal{V}}\right) = \mathcal{M}\left(\nu_{\theta + 1}\right) \times \mathcal{M}\left(\nu_{\theta + 2}\right) \times \dots \times \mathcal{M}\left(\nu_n\right).
    \end{equation*}
\end{definition}

\noindent Therefore, the initial and final densities for a given optimal transport problem are given as coupled measures:
\begin{itemize}
    \item Mass densities:
    \begin{flalign*}
        \rho_0 & := \left(\left.\rho_{e_1}\right|_{t = 0}, \left.\rho_{e_2}\right|_{t = 0}, \dots \left.\rho_{e_m}\right|_{t = 0}\right) \in \mathcal{M}_+^0\left(\mathcal{E}\right) \\
        \rho_T & := \left(\left.\rho_{e_1}\right|_{t = T}, \left.\rho_{e_2}\right|_{t = T}, \dots \left.\rho_{e_m}\right|_{t = T}\right) \in \mathcal{M}_+^T\left(\mathcal{E}\right)
    \end{flalign*}
    \item vertex mass densities:
    \begin{flalign*}
        \gamma_0 & := \left(\left.\gamma_{\nu_{\theta + 1}}\right|_{t = 0}, \left.\gamma_{\nu_{\theta + 2}}\right|_{t = 0}, \dots \left.\gamma_{\nu_n}\right|_{t = 0}\right) \in \mathcal{M}_+^0\left(\mathring{\mathcal{V}}\right) \\
        \gamma_T & := \left(\left.\gamma_{\nu_{\theta + 1}}\right|_{t = T}, \left.\gamma_{\nu_{\theta + 2}}\right|_{t = T}, \dots \left.\gamma_{\nu_n}\right|_{t = T}\right) \in \mathcal{M}_+^T\left(\mathring{\mathcal{V}}\right)
    \end{flalign*}
    \item Source and sink vertex mass densities:
    \begin{flalign*}
        S_0 & := \left(\left.S_{\nu_1}\right|_{t = 0}, \left.S_{\nu_2}\right|_{t = 0}, \dots \left.S_{\nu_o}\right|_{t = 0}\right) \in \mathcal{M}_+^0\left(\partial^+ \mathcal{V}\right) \\
        S_T & := \left(\left.S_{\nu_1}\right|_{t = T}, \left.S_{\nu_2}\right|_{t = T}, \dots \left.S_{\nu_o}\right|_{t = T}\right) \in \mathcal{M}_+^T\left(\partial^+ \mathcal{V}\right) \\
        D_0 & := \left(\left.D_{\nu_{o + 1}}\right|_{t = 0}, \left.D_{\nu_{o + 2}}\right|_{t = 0}, \dots \left.D_{\nu_{\theta}}\right|_{t = 0}\right) \in \mathcal{M}_+^0\left(\partial^- \mathcal{V}\right) \\
        D_T & := \left(\left.D_{\nu_{o + 1}}\right|_{t = T}, \left.D_{\nu_{o + 2}}\right|_{t = T}, \dots \left.D_{\nu_{\theta}}\right|_{t = T}\right) \in \mathcal{M}_+^T\left(\partial^- \mathcal{V}\right)
    \end{flalign*}
\end{itemize}